\documentclass[a4paper,11pt]{article}
\usepackage[english]{babel}

\usepackage[a4paper,margin=2.0cm]{geometry}

\usepackage{amsmath,amssymb}

\usepackage{graphicx}
\usepackage{float}
\usepackage{placeins}
\usepackage{multirow}
\usepackage{authblk}

\usepackage[sort&compress,numbers]{natbib}

\usepackage[colorlinks=true, allcolors=blue]{hyperref}

\usepackage{soul}
\usepackage{xcolor}
\definecolor{mint}{RGB}{180,255,220}
\definecolor{yellow}{RGB}{255,217,102}

\usepackage{lineno}
\title{Identifiability-aware neural ordinary differential equations for parsimonious and reliable dynamic modelling}
\date{}

\author[1,2]{Núria Campo-Manzanares}
\author[1]{Eva Balsa-Canto\thanks{Corresponding author: eva.balsa@csic.es}}

\affil[1]{Biosystems and Bioprocess Engineering Group, IIM-CSIC, Vigo, Spain}
\affil[2]{Applied Mathematics II, University of Vigo, Vigo, Spain}

\begin{document}
\maketitle
\vspace{-3em}

\begin{abstract}

Neural ordinary differential equations (NODE) and hybrid NODE models provide flexible continuous-time representations of complex dynamic systems, but their expressive capacity can exceed the information content of the available data. Consequently, these models may reproduce observed trajectories while retaining weakly identifiable parameters, poorly constrained neural components, and unreliable extrapolation.

Here we introduce identifiability-aware neural ordinary differential equations (iNODE), a framework that incorporates practical identifiability into neural differential equation design. iNODE models embed neural components as explicit analytic functions within the governing equations, enabling direct sensitivity analysis, Fisher-information-based confidence intervals, and identifiability-aware architecture selection. Candidate architectures are generated under data-support constraints, jointly calibrated, and ranked according to predictive accuracy, parsimony, and parameter identifiability.

We evaluate the complete iNODE workflow against conventional NODE and hybrid NODE workflows representative of current practice using four controlled ground-truth benchmarks spanning fully data-driven and hybrid formulations, latent time-varying parameters, partial observability, and sparse or noisy measurements. Using the same training data and evaluation scenarios, the iNODE workflow selected more compact architectures, reduced parameter uncertainty, and improved extrapolation and recovery from latent-dynamics. These results establish practical identifiability as a model-design principle for parsimonious and reliable neural differential equations.
\end{abstract}

\section*{Introduction}

Dynamic-system models are central to understanding, predicting and controlling complex processes in physics, engineering, and the life sciences. Recent advances in machine learning have expanded dynamic modelling by allowing neural components to be embedded directly within differential equations, as part of the broader development of scientific machine learning and physics-informed learning \citep{BruntonNoack2020,Karniadakis2021PIML,Cuomo2022ScientificML}. This idea underpins a broad family of approaches, including physics-informed neural networks (PINN) \citep{raissi2019physics}, neural ordinary differential equations (NODE) \citep{chen2018neural}, hybrid neural ODE (HNODE) \citep{alber2019integrating,shah2022deep,noordijk2024rise,schweidtmann2024review} and universal differential equations (UDE) \citep{rackauckas2020universal}. These methods can represent unknown interactions, latent variables or missing mechanisms directly from data while retaining a continuous-time dynamic formulation.

This flexibility has made neural differential equations increasingly attractive across scientific and engineering domains, including metabolic and systems biology \citep{Faure2023NeuralMechanistic}, microbial community dynamics \citep{Thompson2026Physics}, environmental and water-system modelling \citep{quaghebeur2022hybrid}, weather and climate modelling \citep{Kashinath_etal:weatherclimate}, and neural surrogates for mechanism-based biological models \citep{Wang2019MechanismNN}. Yet, their growing use raises a fundamental question: which parts of the learnt neural structure are actually supported by the available data? In practice, NODE and HNODE models are commonly formulated through workflows in which neural architectures are proposed, tuned, and trained primarily based on predictive performance. Such workflows can improve trajectory fitting, but they do not provide a principled means of determining whether the learnt neural component is genuinely constrained by the data. Consequently, models may reproduce observed trajectories while relying on weakly supported, highly correlated or unnecessarily complex parameterisations. In such cases, a good fit can mask poor extrapolation, unstable parameter estimates and unreliable uncertainty quantification. These limitations become especially important when neural differential equations are used not only for interpolation but also for mechanistic inference, model discrimination, uncertainty quantification or experimental design.

These difficulties are closely related to practical identifiability, defined as the extent to which model parameters can be uniquely estimated from available data in the presence of noise and limited experimental information \citep{walter-pronzato:97}. Inverse problems in non-linear dynamic systems have long demonstrated that limited information content can induce strong parameter correlations, practical non-identifiability and poor predictive reliability, even when data fitting appears satisfactory \citep{Ljung1994,walter1996identifiability,Chis_etal:2016}. This issue is increasingly evident in scientific machine learning. Kharazmi \emph{et al.} showed, in the context of physics-informed neural networks, that identifiability can strongly affect parameter recovery and predictive reliability despite accurate trajectory fits \citep{kharazmi2021identifiability}. Similarly, Giampiccolo \emph{et al.} observed that neural flexibility in HNODE models can compromise the identifiability of the mechanistic component \citep{giampiccolo2024robust}. Recent reviews of universal and hybrid differential equations also identify reliable training with sparse and noisy data, uncertainty quantification and interpretability as open challenges for these models \citep{philipps2025current,schmid2025assessment}.

A more fundamental limitation is that the neural component in standard NODE and HNODE workflows is typically treated as a black-box differentiable module during both training and model selection. Although practical-identifiability analyses can be performed for mechanistic parameters in hybrid formulations, neural-network weights and biases are generally excluded from the identifiability assessment. Consequently, architecture selection remains driven largely by predictive loss or validation performance. Heuristic strategies such as architecture tuning, regularisation and pruning can reduce model complexity or improve training robustness \citep{giampiccolo2024robust,DERESENDEOLIVEIRA2024_RegularizationNN,liu2025automatic}, but they do not resolve the disconnect between neural architecture and practical identifiability.

To address this gap, we introduce identifiability-aware neural ordinary differential equations (iNODE), a framework that integrates practical identifiability directly into the design, calibration and selection of neural differential equations. Rather than treating the neural component solely as a black-box function approximator, iNODE models embed neural components as explicit analytic functions within the governing ordinary differential equations, thereby making neural-network weights and biases explicit parameters of the dynamic model. This representation enables neural and mechanistic parameters to be analysed within a common estimation framework, allowing direct sensitivity analysis, Fisher-information-based confidence intervals and identifiability-informed architecture selection. In doing so, iNODE establishes practical identifiability as a criterion for neural architecture design, linking model complexity to the information content of the available data and addressing a key limitation of conventional NODE and HNODE workflows.

The iNODE workflow comprises architecture generation, joint calibration of mechanistic and neural parameters, practical-identifiability analysis and identifiability-aware model selection. Candidate architectures are generated according to constraints imposed by both the modelling objective and the information content of the available data. The resulting models are then calibrated and ranked according to predictive accuracy, parsimony and parameter uncertainty, with weakly supported components pruned where appropriate. In this way, model complexity is aligned with the information provided by the observations rather than being determined by predictive fit alone.

We evaluate the complete iNODE workflow using four controlled ground-truth systems spanning fully data-driven and hybrid formulations, latent time-varying parameters, full and partial observability, and sparse or noisy measurements. This benchmark design is particularly valuable because it allows direct assessment of whether a model recovers the governing vector field, mechanistic parameters, latent functions and out-of-sample dynamics, rather than simply reproducing observed trajectories. The benchmarks also enable the effects of noise, sampling density, observability and experimental design to be examined under controlled conditions.

Conventional NODE and HNODE baselines were implemented using training and model-selection procedures representative of current practice, including architecture selection driven primarily by predictive performance and calibration via backpropagation through the ODE solver with local gradient-based optimisation \citep{Rumelhart1986,lecun2002efficient,KingmaBa2015Adam,chen2018neural}. Both workflows were evaluated using the same training data and validation scenarios, while preserving their respective calibration and architecture-selection strategies. Across all four benchmarks, iNODE consistently selected more parsimonious architectures, produced better-constrained parameter estimates and improved extrapolation and latent-dynamics recovery. These findings support practical identifiability as a guiding principle for the design of neural differential equations and highlight the benefits of incorporating identifiability directly into model construction and selection.

\section*{Results}

\begin{figure}[h!]
    \centering
    \includegraphics{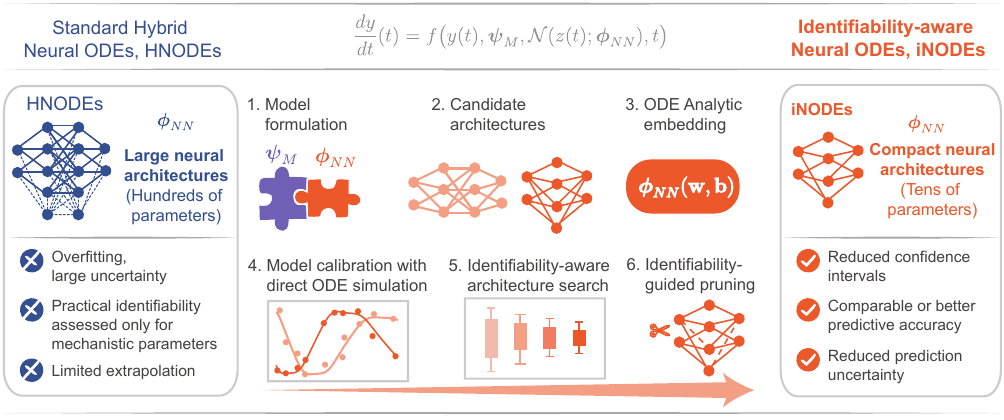}
    \caption{\textbf{The iNODE workflow links architecture search to identifiability.} Conventional HNODE workflows primarily select neural architectures according to predictive performance. In the existing identifiability-oriented literature, practical-identifiability analysis has focused on mechanistic parameters. In the iNODE framework, candidate architectures are systematically generated, analytically embedded into the governing equations, calibrated and ranked through direct practical identifiability analysis. This workflow links neural model complexity to the information content of the data, yielding compact neural differential equations with fewer parameters, narrower confidence intervals, reduced prediction uncertainty, and predictive performance comparable to or better than that of conventional HNODE models.}
    \label{fig:pipeline}
\end{figure}

\subsection*{The iNODE workflow links architecture search to identifiability}

iNODE models integrate neural architecture search, analytic embedding, parameter estimation and model ranking within a single identifiability-aware workflow. Figure~\ref{fig:pipeline} summarises this logic and contrasts it with conventional NODE and HNODE workflows, where neural architectures are commonly selected mainly based on predictive fit.

The central feature of the iNODE framework is that practical identifiability is incorporated into the model formulation. By analytically embedding neural components within the governing equations, the framework makes mechanistic parameters and neural parameters (weights and biases) accessible to the same sensitivity and uncertainty analysis. Candidate architectures can therefore be assessed not only by predictive accuracy, but also by parameter uncertainty and structural parsimony.

\textbf{\textit{Step 1: Model formulation.}}
We consider neural differential equations of the form
\begin{equation}
\frac{dy}{dt}(t)=f\big(y(t),\boldsymbol{\psi}_M,\mathcal{NN}(z(t);\boldsymbol{\phi}_{NN}),t\big),
\qquad y(t_0)=y_0,
\end{equation}
where \(y(t)\in\mathbb{R}^n\) denotes the state of the system, \(\boldsymbol{\psi}_M\) the mechanistic parameters, and \(\mathcal{NN}(\cdot;\boldsymbol{\phi}_{NN})\) a neural component parameterised by weights and biases collected in \(\boldsymbol{\phi}_{NN}\). The input of the network \(z(t)\) may include time, state variables, or both, depending on the modelling context.

This formulation covers several modelling scenarios considered here, including latent parameter dynamics, in which neural networks represent time-varying parameters; state closure, in which neural terms complement mechanistic equations; and mechanistic NODE, in which neural networks replace unknown mechanistic equations while retaining explicit state dependence. The corresponding mathematical formulations are detailed in the Supplementary Information, Section~S.1.

\textbf{\textit{Step 2: Candidate neural architectures.}}
Once the model structure is defined, candidate neural architectures are generated by varying the number of layers, the number of neurons and the activation functions. To limit overparameterisation at this stage, we retain only architectures that satisfy a minimum ratio of four observations per trainable parameter, with \(\boldsymbol{\theta}=[\boldsymbol{\psi}_M,\boldsymbol{\phi}_{NN}]\). This conservative filter excludes clearly underdetermined neural parameterisations during architecture generation, whereas practical identifiability is evaluated explicitly after calibration.

\textbf{\textit{Step 3: Analytic embedding.}}
Each candidate network is written as an explicit analytic expression and embedded directly into the differential equations. For a feedforward neural network with \(L\) hidden layers, the neural mapping is defined recursively as
\begin{equation}
\mathbf{h}^{(1)}=\sigma^{(1)}\!\left(\mathbf{W}^{(1)}\mathbf{z}+\mathbf{b}^{(1)}\right),
\end{equation}
\begin{equation}
\mathbf{h}^{(l)}=\sigma^{(l)}\!\left(\mathbf{W}^{(l)}\mathbf{h}^{(l-1)}+\mathbf{b}^{(l)}\right),
\qquad l=2,\dots,L,
\end{equation}
and
\begin{equation}
\mathcal{NN}(\mathbf{z};\boldsymbol{\phi}_{NN})=
\sigma^{(\mathrm{out})}\!\left(\mathbf{W}^{(\mathrm{out})}\mathbf{h}^{(L)}+\mathbf{b}^{(\mathrm{out})}\right),
\end{equation}
where \(\mathbf{W}^{(l)}\), \(\mathbf{b}^{(l)}\) and \(\sigma^{(l)}\) denote the weights, biases and activation function of layer \(l\).

Substitution of this expression into the original model yields the analytically embedded neural system
\begin{equation}
\frac{dy}{dt}(t)=\tilde{f}\big(y(t),\boldsymbol{\psi}_M,\boldsymbol{\phi}_{NN},t\big),
\qquad y(t_0)=y_0,
\end{equation}
in which mechanistic and neural parameters appear explicitly within a unified ODE model.

\textbf{\textit{Step 4: Calibration and identifiability analysis.}}
The embedded neural model is calibrated by estimating mechanistic and neural parameters to minimise the discrepancy between model predictions and the available data. Because the neural component remains explicit within the ODE system, calibration and practical identifiability analysis are performed within the same formulation. The fitted model can therefore be analysed through parameter sensitivities and Fisher-information-based diagnostics, enabling confidence-interval estimation and direct assessment of which parameters are supported by the data. Full details of the calibration objective, optimisation procedure and identifiability analysis are provided in the Methods section.

\textbf{\textit{Step 5: Identifiability-aware architecture selection.}}
Each calibrated candidate is evaluated using predictive accuracy, quantified by the normalised root mean square error (NRMSE), model parsimony, quantified by the Akaike information criterion (AIC), and parameter uncertainty, summarised by maximum and median confidence intervals. These criteria are combined into the Model Quality Ranking Index (MQRI), a scale-independent composite score that ranks candidate architectures according to their balance between fit quality, complexity and identifiability. Lower MQRI values indicate models with the most favourable trade-off, enabling automated and reproducible architecture selection. Details of the MQRI computation are provided in the Methods section.

\textbf{\textit{Step 6: Identifiability-guided pruning.}}
After architecture selection, an optional pruning phase can further reduce redundancy and improve parameter robustness. Neural parameters with negligible estimated contributions and large relative confidence intervals are removed, and the remaining parameters are re-estimated. This targeted pruning improves the conditioning of the estimation problem and reduces uncertainty without compromising predictive accuracy.

We benchmarked the workflow using four controlled systems with known generating dynamics. This design enables direct assessment of state prediction, vector-field and latent-function recovery, parameter uncertainty and extrapolation under controlled changes in noise, sampling density, observability and experimental design. Further details are provided in the Supplementary Information, Section~S.2.

\subsection*{A compact iNODE improves extrapolation and phase-space recovery in a fully data-driven cubic damped oscillator}

The cubic damped oscillator has recently been used to assess whether NODE models can learn nonlinear continuous-time dynamics from observed trajectories~\citep{goyal2023neural,stepaniants2024discovering}. Here, we use a two-dimensional cubic damped oscillator with strong nonlinearities and sensitivity to initial conditions to test whether identifiability-aware architecture design improves robustness, phase-space recovery and extrapolation in fully data-driven and hybrid settings. The ground-truth system and the corresponding NODE formulations are described in Supplementary Information, Section~S.3.

We first considered a fully data-driven NODE formulation following Goyal and Benner~\citep{goyal2023neural}. In this setting, the system dynamics are represented entirely by a neural network. Synthetic training data were generated from the ground-truth model using the initial conditions \((x_1(0),x_2(0))=(2,0)\,\mathrm{m}\), yielding 1250 measurements per state variable. The data were normalised before training, as this is standard practice and often essential for stable optimisation in conventional NODE workflows.

\begin{figure}[]
\centering
\includegraphics{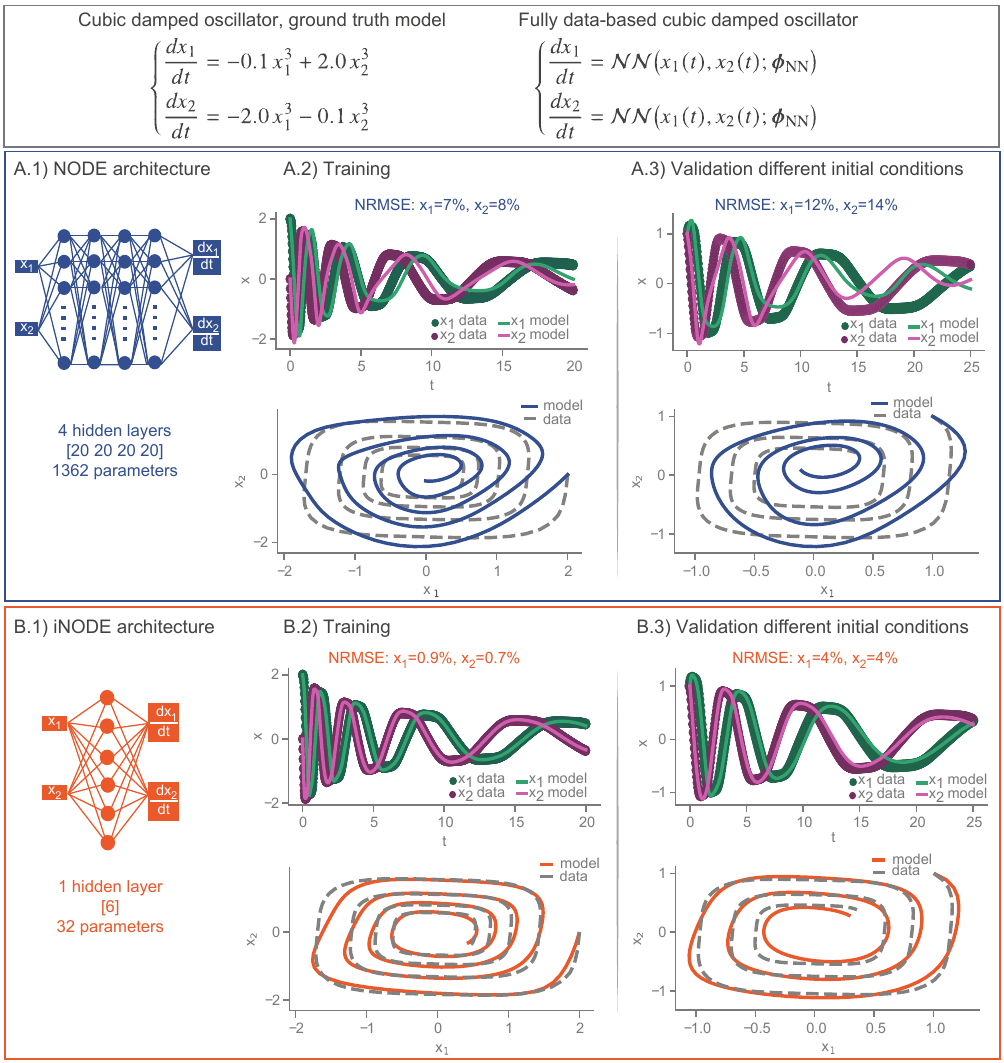}
\caption{\textbf{A compact iNODE improves extrapolation and phase-space recovery compared to the conventional NODE model in a fully data-driven cubic damped oscillator.} \textbf{A.1,} Architecture of the conventional NODE baseline. \textbf{A.2--A.3,} Conventional NODE predictions for the training and validation experiments. \textbf{B.1,} Architecture of the selected iNODE. \textbf{B.2--B.3,} Corresponding iNODE predictions for the training and validation scenarios. Trajectory plots and phase portraits are shown for each experiment and neural architecture. Although both models fit the training data, the conventional NODE extrapolates poorly and distorts the underlying vector field, whereas the more compact identifiability-aware iNODE preserves the phase-space structure and generalises more reliably.}
\label{fig:NODEs_cdm_models}
\end{figure}

As a reference, we trained the conventional NODE model using the architecture proposed in~\citep{goyal2023neural}, namely four hidden layers with 20 neurons each, corresponding to 1362 trainable parameters (Figure~\ref{fig:NODEs_cdm_models}, Panel~\textbf{A.1}). Generalisation was assessed using a validation experiment with unseen initial conditions \((x_{1,\mathrm{val}}(0),x_{2,\mathrm{val}}(0))=(1,1)\,\mathrm{m}\). The model reproduced the training data with NRMSE values of 7\% and 8\% for the two states (Panel~\textbf{A.2}), but the reconstructed phase portrait deviated noticeably from the ground truth. These discrepancies became more pronounced in validation (Panel~\textbf{A.3}), where errors increased to approximately 12\% and 14\% and the learnt dynamics produced oscillations with artificially shortened periods. Thus, despite fitting the training trajectory reasonably well, the conventional NODE model failed to recover the correct phase-space geometry and extrapolated poorly beyond the calibration regime.

We then applied the iNODE workflow to the same problem and dataset. Although the data could in principle support architectures with up to 625 trainable parameters under the heuristic of at least four observations per parameter, we restricted the analysis to the first 20 candidate architectures of increasing complexity. Model selection was based on the identifiability-aware MQRI criterion (Supplementary Information, Table~T.4). A clear trade-off emerged between model complexity and parameter identifiability: more complex models often achieved excellent fits, but showed stronger parameter correlations and poorer identifiability. Single-hidden-layer architectures provided the best compromise between predictive accuracy and parameter robustness. The selected iNODE model consisted of one hidden layer with six neurons, corresponding to 32 trainable parameters (Figure~\ref{fig:NODEs_cdm_models}, Panel~\textbf{B.1}). Detailed parameter values, sensitivity and identifiability analyses are reported in Supplementary Information, Section~S.3.2.

Using the same training data and validation initial conditions, the selected iNODE model outperformed conventional NODE on both training and validation sets. In the training experiment, the iNODE model achieved NRMSE values below 1\% for both state variables and reproduced the phase portrait faithfully (Figure~\ref{fig:NODEs_cdm_models}, Panel~\textbf{B.2}). In validation, the model remained robust, with errors of approximately 4\% for both states and smooth trajectories that captured both oscillation amplitude and period (Panel~\textbf{B.3}).

Importantly, identifiability of the conventional NODE model was assessed using the iNODE framework. After training, the selected neural architecture was analytically embedded in the ODE system, allowing sensitivity- and Fisher-information-based diagnostics. This analysis showed that the Fisher information matrix was non-invertible, preventing confidence-interval estimation. 

Parametric sensitivity analysis revealed values spanning several orders of magnitude, with many parameters showing weak or strongly overlapping sensitivity profiles. The resulting sensitivity matrix was rank-deficient, indicating redundant and poorly informed parameter directions (Supplementary Information, Section~S.3.1).

\subsection*{Experimental diversity improves identifiability and recovery in hybrid iNODE models in a hybrid cubic oscillator}

We next considered a hybrid mechanistic-neural formulation of the cubic damped oscillator, in which the dynamics of \(x_1\) are specified mechanistically and \(x_2\) is learnt from data, i.e. the time derivative of \(x_2\) is represented by a neural network. We use this example to examine how identifiability depends not only on the choice of architecture but also on the diversity of available experiments. We therefore used the iNODE framework to analyse the same hybrid architecture across different experimental scenarios with equal data budgets but varying initial-condition diversity and then compared the best-performing configuration against a conventional HNODE baseline.

\begin{figure}[]
\centering
\includegraphics{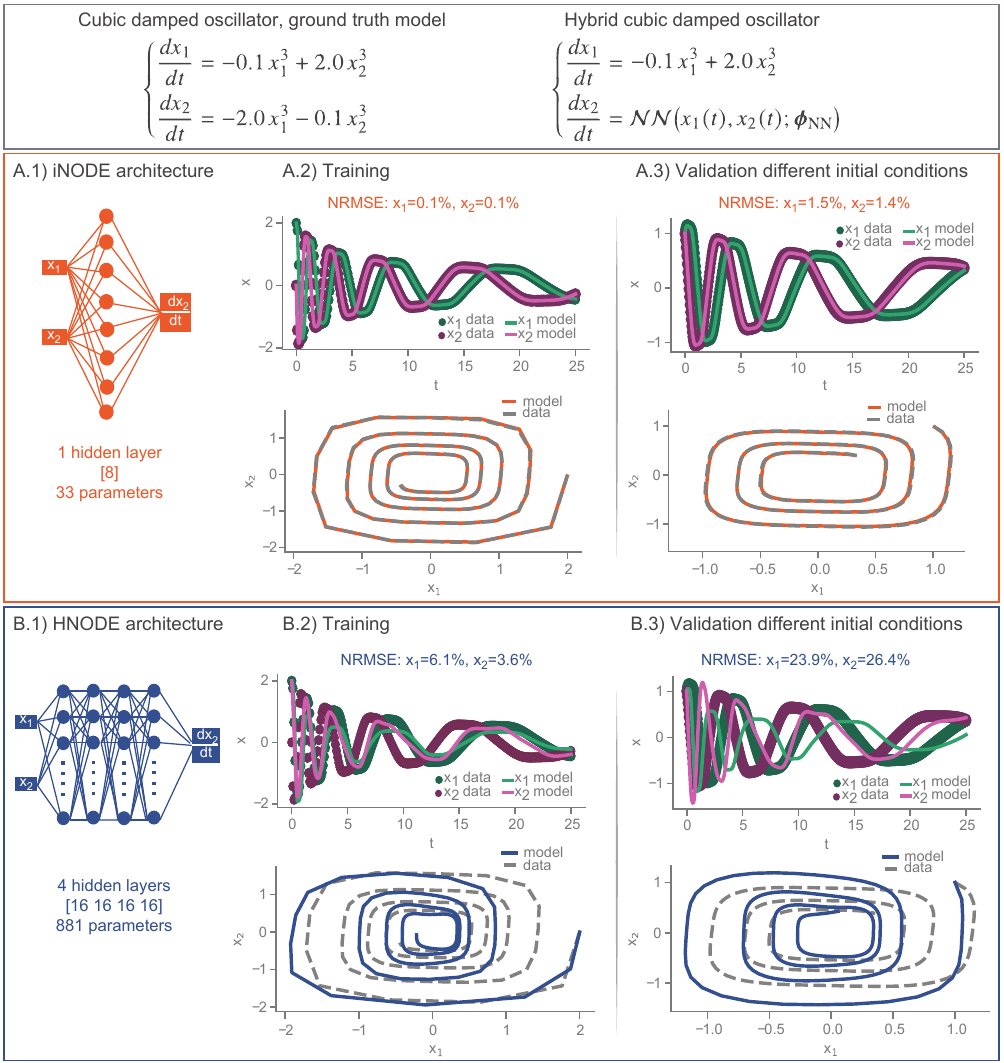}
\caption{\textbf{Experimental diversity and identifiability-aware pruning improve recovery and generalisation with a hybrid iNODE model relative to a conventional HNODE model in a cubic damped oscillator.} \textbf{A.1,} Architecture of the selected hybrid iNODE model. \textbf{A.2-A.3,} Hybrid iNODE results for training and the validation experiments with different initial conditions. \textbf{B.1,} Architecture of the conventional HNODE baseline. \textbf{B.2-B.3,} Corresponding conventional HNODE predictions for the training and validation scenarios. Trajectory plots and phase portraits are shown for each experiment and neural architecture. The selected hybrid iNODE model, obtained after analysing the effect of experimental diversity and identifiability-guided pruning, accurately recovers the oscillatory dynamics and preserves the phase-space structure in training and validation. The much larger conventional HNODE model exhibits substantially larger errors, converges towards overly damped dynamics and fails to recover the correct phase-space portrait.}
\label{fig:hybrid_cdm}
\end{figure}

To select the neural architecture associated with the data-driven component, we used a dataset comprising 500 data points per state from a single experiment with initial conditions \((x_1(0), x_2(0))=(2,0)\,\mathrm{m}\) and applied the iNODE workflow. Models with up to 250 parameters were admissible, given the ratio between the amount of data and the number of parameters. A total of 12 architectures were evaluated and single-hidden-layer networks consistently outperformed deeper alternatives, resulting in lower errors and tighter confidence intervals (Supplementary Information, Table~T.7). The selected hybrid iNODE architecture consisted of a neural network with a hidden layer of eight neurons and 33 trainable parameters (Figure~\ref{fig:hybrid_cdm}, Panel~\textbf{A.1}).

We then examined how parameter identifiability and predictive performance changed in two experimental scenarios. In the single-experiment setting (\texttt{1exp}), the selected model accurately reconstructed the observed trajectories but exhibited practical identifiability limitations, with a maximum confidence interval of 881\%. We next considered a two-experiment design (\texttt{2exp}) with distinct initial conditions but the same amount of data to test whether increased dynamic diversity could improve parameter confidence. In the \texttt{2exp} scenario, the maximum confidence interval was 278\%.

Finally, identifiability-guided pruning was applied by removing parameters with negligible magnitude (\(|\theta|<0.1\)) and poor identifiability (CI \(>100\%\)). Pruning further reduced uncertainty, yielding a maximum confidence interval of 84\% without compromising predictive performance. These results show that data diversity provides the main gain in practical identifiability, while pruning further improves the conditioning of the selected hybrid model. More details on training, explored architectures, and performance metrics are provided in Supplementary Information, Section~S.4.1.

To benchmark predictive performance, we compared the selected two-experiment pruned hybrid iNODE model with a conventional HNODE baseline trained on the same data set. The baseline architecture was selected by grid search and consisted of four hidden layers with 16 neurons each, corresponding to 881 trainable parameters, compared to only 33 trainable parameters in the selected hybrid iNODE model (Figure~\ref{fig:hybrid_cdm}, Panels~\textbf{B.1} and~\textbf{A.1}). Although this deeper HNODE architecture achieved the lowest training error among the conventional candidates explored, it remained substantially larger and less constrained than the identifiability-aware alternative, making it a stringent baseline for comparison.

Comparison of the two models reveals a clear performance gap. Despite the much larger architecture, the conventional HNODE model exhibited relatively large training errors (NRMSE of 6.1\% and 3.6\%; Figure~\ref{fig:hybrid_cdm}, Panel~\textbf{B.2}) and performed poorly in validation, reaching 23.9\% and 26.4\% (Panel~\textbf{B.3}). The selected hybrid iNODE model achieved training NRMSE below 0.12\% and validation errors below 1.5\% while preserving the oscillatory phase-space structure (Panels~\textbf{A.2--A.3}). The conventional HNODE model, in contrast, converged towards overdamped dynamics and distorted the attractor geometry. Further details can be found in the Supplementary  Section~S.4.2.

\subsection*{A compact hybrid iNODE model robustly recovers latent time-dependent interactions from noisy, sparse data in a Lotka--Volterra model}

We next considered a class of hybrid neural models in which the neural component represents latent time-dependent parameters embedded within an otherwise mechanistic system. This setting is especially relevant for the microbial communities described by generalised Lotka--Volterra (gLV) models, where interaction coefficients may vary over time due to physiological changes, resource fluctuations, or higher-order effects~\citep{momeni2017lotka,mustri2025accuracy}. Here, we used a two-species Lotka--Volterra benchmark in which the dynamics of both species are observed and the latent interaction functions are inferred from the data. The ground-truth system and the corresponding hybrid formulation are shown in Figure~\ref{fig:LV-noisycase} and are described in detail in the Supplementary Information, Section~S.5.

We focus on a realistic and experimentally challenging regime with sparse measurements and 10\% Gaussian noise. Three training experiments with different initial conditions were generated, yielding a total of 60 observations. Model calibration was performed on base-10 log-transformed observations to improve parameter estimation under noisy, data-limited conditions.

Within the hybrid iNODE framework, candidate neural architectures were generated under the adopted data-to-parameter constraint and calibrated using the identifiability-aware workflow, resulting in five admissible models (Supplementary Information, Table~T.11). The MQRI-selected architecture consisted of a single hidden layer with two neurons and achieved a favourable compromise between predictive precision, parsimony, and practical identifiability, with improved recovery of \(\mu_2\) and \(k_2\) and reduced uncertainty in neural parameters (maximum CI: 83\%). For comparison, a conventional hybrid NODE baseline was implemented in Python and selected by grid search on the same log-transformed data set using NRMSE (Supplementary Information, Table~T.13).

As shown in Figure~\ref{fig:LV-noisycase}, the two approaches converged to markedly different architectures. The conventional hybrid NODE selected a substantially larger network with 630 trainable parameters (Figure~\ref{fig:LV-noisycase}, Panel~\textbf{A.1}), whereas the hybrid iNODE selected a compact model with only 14 trainable parameters (Panel~\textbf{B.1}). Despite this pronounced difference in complexity, the hybrid iNODE model achieved similar training accuracy and consistently better generalisation.

\begin{figure}[]
    \centering
    \includegraphics{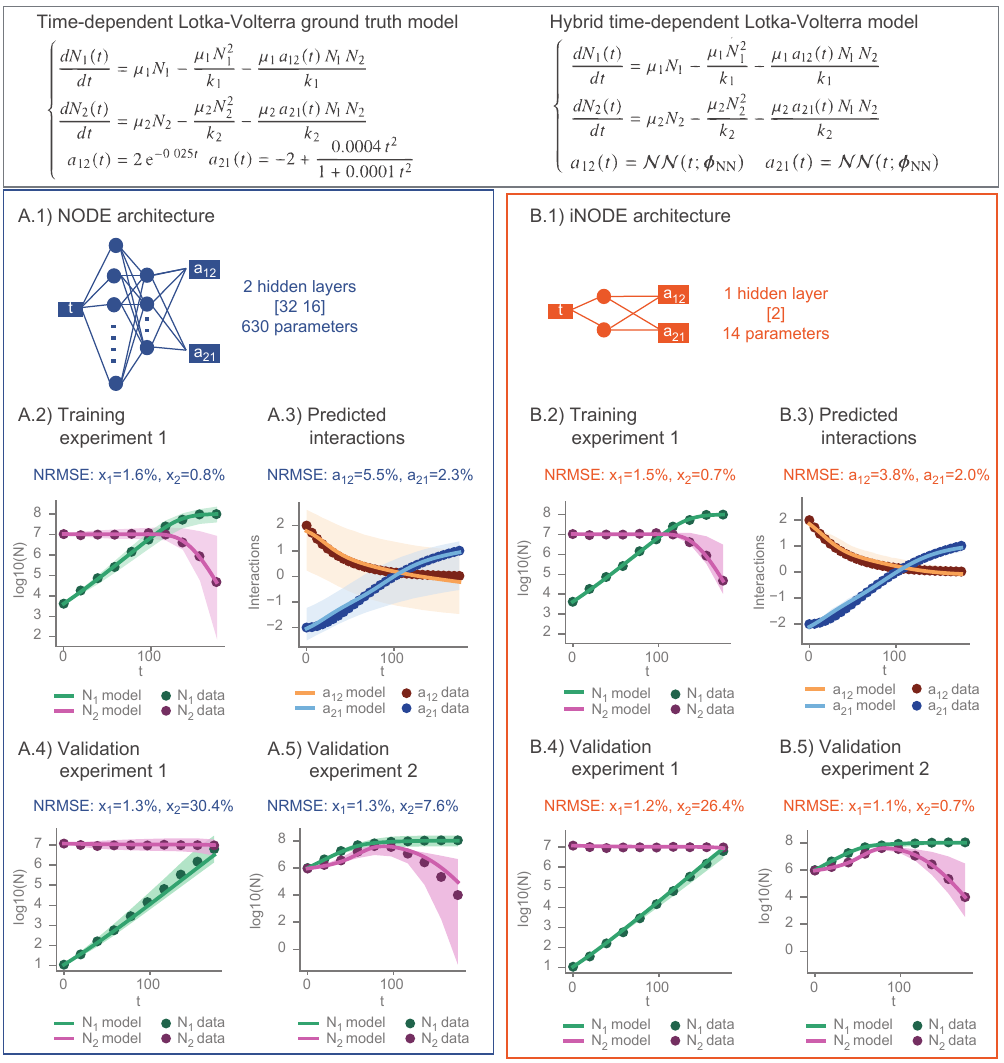}
    \caption{\textbf{A hybrid iNODE model recovers latent interaction dynamics more accurately than a substantially larger conventional hybrid NODE baseline.}
    \textbf{A.1,} Architecture of the conventional hybrid NODE baseline. \textbf{B.1,} Architecture of the selected hybrid iNODE. \textbf{A.2-A.5,} Baseline predictions and inferred interaction functions for training and validation experiments. \textbf{B.2-B.5,} Corresponding results for the hybrid iNODE model. Points denote observations, solid lines model predictions, and shaded regions uncertainty bands. Despite its compact architecture, hybrid iNODE achieves comparable or superior trajectory prediction and more accurate recovery of latent interaction dynamics under sparse and noisy conditions.}
    \label{fig:LV-noisycase}
\end{figure}

Generalisation was evaluated in two validation scenarios with unseen initial conditions. In the first, the hybrid iNODE model yielded a lower prediction error for species~2 than the baseline conventional hybrid NODE (NRMSE=26.4\% versus 30.4\%; Figure~\ref{fig:LV-noisycase}, Panels~\textbf{A.4} and~\textbf{B.4}). In the second, the difference became much larger, and the hybrid iNODE model achieved an NRMSE of 0.7\% for species~2 compared to 7.6\% for the baseline (Panels~\textbf{A.5} and~\textbf{B.5}). The recovery of latent interaction functions further highlighted this advantage: the hybrid iNODE model reconstructed \(a_{12}(t)\) and \(a_{21}(t)\) with NRMSE values of 3.8\% and 2.0\%, compared to 5.5\% and 2.3\% for conventional hybrid NODE (Panels~\textbf{A.3} and~\textbf{B.3}). The larger reconstruction error in \(a_{12}(t)\) obtained with the baseline explains the poorer predictive performance observed for species~2 in both validation scenarios.


Uncertainty quantification further differentiated the two approaches. For the hybrid iNODE model, confidence intervals were obtained directly through the identifiability-aware pipeline, allowing uncertainty bands to be propagated to both state trajectories and latent interaction functions (Figure~\ref{fig:LV-noisycase}, Panels~\textbf{B.2--B.5}). In contrast, sensitivity analysis of the conventional HNODE model revealed strong parameter correlations and a severely rank-deficient sensitivity matrix, precluding direct confidence-interval estimation through inversion of the Fisher Information Matrix (Supplementary Information, Section~S.5.2). To compute confidence intervals in this case, we used a bootstrap resampling method. Subsequently, these confidence intervals were used to generate prediction uncertainty bands for both the state trajectories and the inferred interaction functions of the HNODE model (Panels~\textbf{A.2--A.5}) following the ensemble-based uncertainty propagation methodology described in the Methods section. Despite achieving a comparable fit to the data, the conventional HNODE model exhibited broader uncertainty and less constrained predictions than the hybrid iNODE formulation.

\subsection*{A compact hybrid iNODE model recovers latent time-dependent transmission under partial observability in an SIR model}

We next considered a partially observable latent-parameter setting based on the susceptible-infectious-recovered (SIR) framework. SIR models provide a mechanistic basis for describing transmission dynamics in epidemics. However, the common assumption of a constant transmission rate is rarely valid in real outbreaks, as behavioural changes, interventions, and population heterogeneity can induce pronounced temporal variation~\citep{osi2024parameter}. We considered a scenario in which, given observations of the infected population, the time-varying transmission rate must be inferred from the data using a hybrid NODE formulation. The ground truth model and the corresponding hybrid formulation are shown in Figure~\ref{fig:SIR} and are described in detail in Supplementary Information, Section~S.6. Only the infected population \(I(t)\) was used for calibration; the susceptible and recovered trajectories were retained from the synthetic ground truth only for the validation of state recovery.

We first used a hybrid iNODE model to examine how the experimental design affects parameter identifiability under partial observability. The models were calibrated and evaluated using 60 synthetic measurements of the infected population. Candidate architectures were generated under the adopted data-to-parameter constraint, and all admissible models were calibrated and ranked using the identifiability-aware workflow (Supplementary Information, Table~T.15). The selected configuration consisted of a single hidden layer with two neurons and only 8 trainable parameters (Figure~\ref{fig:SIR}, Panel~\textbf{A.1}). Slightly larger networks yielded only marginal gains in fit quality while showing worse AIC values and substantially poorer practical identifiability.

To assess the role of data richness and experimental diversity, the selected architecture was calibrated under three scenarios: i) a short observation window with 60 measurements and $(S^{(1)}(0), I^{(1)}(0), R^{(1)}(0))$ $=(100000, 100, 0)$ cases, ii) a longer observation window with 80 measurements and iii) a two-experiment scenario combining 40 measurements obtained under two different initial conditions $(S^{(2)}(0),$ $I^{(2)}(0), R^{(2)}(0)) = (50000, 500, 0)$, with 40 measurements per experiment. The corresponding parameter estimates and confidence intervals are reported in the Supplementary Information, Table~T.16. 

\begin{figure}[]
    \centering
    \includegraphics{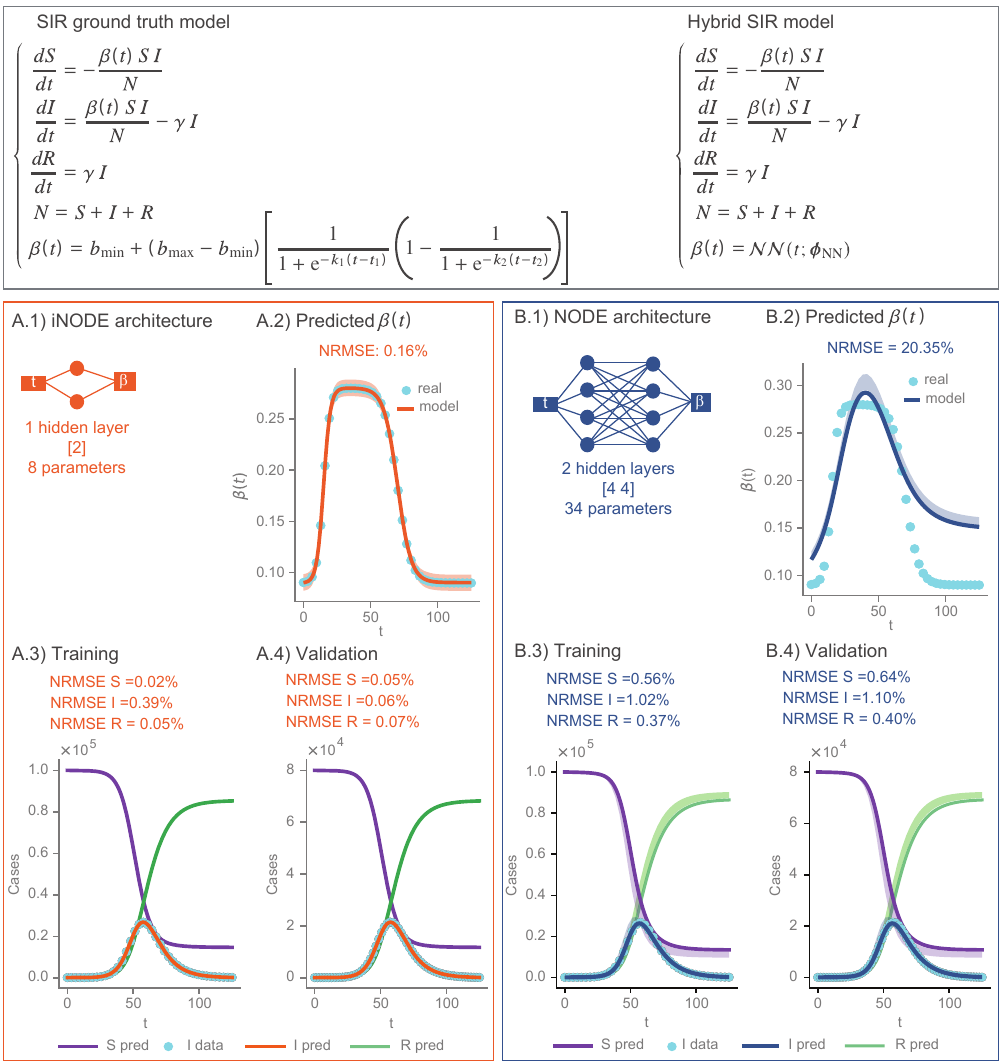}
   \caption{\textbf{The hybrid iNODE model recovers time-varying transmission dynamics more accurately than a larger conventional HNODE baseline under partial observability.} \textbf{A.1,} Architecture of the selected hybrid iNODE. \textbf{A.2,} Recovered transmission-rate profile \(\beta(t)\) for the hybrid iNODE model. \textbf{A.3-A.4,} Hybrid iNODE predictions for the training and validation experiments. \textbf{B.1,} Architecture of the conventional HNODE baseline. \textbf{B.2,} Recovered transmission-rate profile \(\beta(t)\) for the conventional HNODE model. \textbf{B.3-B.4,} Corresponding conventional HNODE predictions for the training and validation experiments. Points denote observations of the infected population, solid lines denote model predictions, and shaded regions indicate uncertainty bands. Despite its compact architecture, hybrid iNODE accurately reconstructs the latent transmission dynamics and yields lower prediction errors and tighter uncertainty bands than the larger conventional HNODE baseline.}
    \label{fig:SIR}
\end{figure}

Practical identifiability was strongly dependent on the experimental design. The single-experiment model showed severe parametric uncertainty, with maximum confidence intervals reaching 1757\%, and extending the observation window reduced this only partially (maximum CI: 474\%). The two-experiment design reduced the uncertainty by two to three orders of magnitude, with a maximum CI of only 17\%, and yielded a substantially more robust reconstruction of the states and transmission rate. These results show that under partial observability, experimental diversity is essential for constraining latent neural parameters and reliably recovering the unobserved epidemic dynamics.

We then compared the selected two-experiment hybrid iNODE model with a conventional HNODE baseline trained on the same data set. The baseline architecture was selected by grid search (Supplementary Information, Table~T.17) and consisted of two hidden layers with four neurons each, corresponding to 34 trainable parameters, more than four times the size of the selected hybrid iNODE model (Figure~\ref{fig:SIR}, Panels~\textbf{A.1} and~\textbf{B.1}). Generalisation was evaluated using a validation scenario with distinct initial conditions.

As shown in Figure~\ref{fig:SIR}, the hybrid iNODE model reconstructed the transmission-rate profile with very high precision, achieving an NRMSE of 0.16\% (Panel~\textbf{A.2}). The NRMSE values in training and validation are well below 0.5\% (Panels~\textbf{A.3-A.4}). The baseline of HNODE still fits the observed epidemic trajectories reasonably well (NRMSE$\leq$1.10\%, Panels~\textbf{B.3-B.4}), yet it fails to recover the underlying transmission-rate dynamics (NRMSE of 20.35\%, Panel~\textbf{B.2}), indicating that the fit of the trajectory-level alone is insufficient to guarantee reliable inference of latent variables.

These differences between models were also reflected in uncertainty quantification. The hybrid iNODE model yielded narrow uncertainty bands consistent with the tight confidence intervals obtained from the identifiability-aware pipeline, whereas conventional HNODE showed wider uncertainty due to larger parameter correlations and poorer practical identifiability (Supplementary Information, Section~S.6.2).

\section*{Discussion}

We have introduced identifiability-aware neural ordinary differential equations (iNODE), a framework that incorporates practical identifiability directly into the design of neural differential equations. By analytically embedding neural components within the governing equations, the framework exposes mechanistic and neural parameters to a common calibration, sensitivity and uncertainty analysis. This representation allows model complexity to be determined by the information content of the data rather than by predictive fit alone, establishing a principled connection between neural architecture, parameter uncertainty and predictive reliability.

This advance addresses a central challenge in scientific machine learning: flexible models can reproduce observed trajectories without reliably recovering the underlying system dynamics. Previous studies have shown that identifiability can strongly affect parameter recovery and prediction in physics-informed neural networks, even when trajectory fitting is accurate \citep{kharazmi2021identifiability}. Neural flexibility can similarly compromise the identifiability of mechanistic components in HNODE models \citep{giampiccolo2024robust}, while reliable training under sparse and noisy data, uncertainty quantification and interpretability remain major challenges for UDE formulations in systems biology \citep{philipps2025current}. The iNODE framework extends these observations by incorporating identifiability directly into model construction and selection, rather than treating it solely as a post hoc diagnostic tool.

Across four ground-truth benchmarks, iNODE consistently selected compact architectures, produced better-constrained parameter estimates and achieved superior performance under unseen conditions compared with conventional neural differential-equation workflows. The controlled benchmark design allowed trajectory agreement to be distinguished from recovery of the generating dynamics. Because the governing vector fields, mechanistic parameters and latent functions were known, it was possible to determine directly whether accurate fits corresponded to faithful dynamic inference. The benchmark systems spanned fully data-driven and hybrid formulations, latent time-varying processes, partial observability, and sparse or noisy measurements. The largest gains were observed in settings where flexible dynamic models are particularly difficult to constrain, including incomplete state measurements, indirectly observed latent processes and extrapolation beyond the calibration domain. Taken together, these results identify identifiability as a key determinant of whether a neural differential-equation model can support reliable inference and prediction.

A central finding of this study is that predictive accuracy alone provides an insufficient basis for architecture selection. Large conventional NODE and HNODE models often reproduced observed trajectories while retaining strong parameter correlations, rank-deficient sensitivity matrices, large uncertainty or poor recovery of the governing vector field and latent dynamics. In contrast, the iNODE workflow favoured substantially smaller architectures that preserved or improved predictive performance while remaining better constrained by the data. The framework therefore differs fundamentally from generic complexity regularisation because it links architectural choice directly to the parameter combinations supported by the available experiments.

This distinction is important because architecture size and predictive loss do not uniquely determine dynamic reliability. Two parameterisations may achieve similar trajectory errors yet represent markedly different vector fields, phase-space structures or latent functions. These differences become particularly evident during extrapolation, when models are evaluated outside the region directly constrained by the training observations. By combining predictive performance with parameter-level uncertainty, iNODE favours architectures that capture the observed dynamics without relying on weakly informed directions in parameter space. This balance between fit, parsimony and identifiability helps explain the improved extrapolation and latent-function recovery observed across the benchmark systems.

Benchmarking against conventional NODE and HNODE workflows also highlights the interaction between model representation, architecture selection and calibration. Conventional baselines were trained using backpropagation through the ODE solver, local gradient-based optimisation and established stabilisation procedures \citep{Rumelhart1986,lecun2002efficient,KingmaBa2015Adam,chen2018neural}. These approaches underpin much of current neural differential-equation practice, including many HNODE, UDE and PINN applications. The resulting optimisation problems are strongly non-convex and often contain flat directions and ill-conditioned regions, particularly when neural architectures are only weakly constrained by the measurements. Repeated initialisation, learning-rate scheduling, gradient clipping and multiple shooting can improve numerical robustness, but the resulting models may still retain redundant parameterisations that reproduce trajectories without recovering the underlying dynamics.

The case studies illustrate how these factors interact. In the hybrid cubic-oscillator benchmark, the same nominal MLP(8) architecture produced substantially different outcomes across the two workflows. In the partially observed SIR benchmark, the MLP(2) architecture selected by iNODE was also among the candidate architectures considered in the conventional hybrid NODE workflow, yet the latter did not recover the latent transmission dynamics with comparable accuracy. These findings demonstrate that architecture size alone cannot explain the observed performance differences. Rather, analytic embedding, calibration and identifiability-aware selection act together to determine whether a model is both dynamically informative and statistically supported.

This perspective complements recent efforts to replace manually tuned neural structures with more principled architecture-discovery strategies in scientific machine learning, including physics-informed distillation approaches for physics-informed neural networks \citep{liu2025automatic}. For dynamic inference, the present results suggest an additional requirement: candidate architectures should be evaluated according to whether their parameters are supported by the available observations. Identifiability provides a direct criterion for this purpose and connects architecture selection to established principles of inverse problems and dynamic-system identification.

The rank-deficient sensitivity matrices observed for several conventional models reveal the consequences of ignoring this information structure. Large numbers of measurements did not necessarily eliminate parameter redundancy, even in densely sampled scenarios. Measurement count alone is therefore a poor surrogate for inferential information. Parameter identifiability depends instead on whether the experiments excite sufficiently distinct parameter directions and generate observations capable of distinguishing among competing dynamic mechanisms. This principle is well established in systems identification and model-based experimental design, where parameter uncertainty is governed by experimental informativeness rather than sample size alone \citep{Banga-Balsa:2008,Franceschini2008,BalsaCanto2025Uncertainty}.

The hybrid cubic-oscillator and partially observed SIR benchmarks further demonstrate the importance of dynamic diversity. Combining experiments performed under distinct initial conditions or perturbations improved parameter identifiability more effectively than extending a single trajectory or increasing neural flexibility. Complementary experiments probe different regions of the state space and activate different parameter sensitivities, thereby reducing correlations and improving uncertainty estimates. Because the iNODE formulation provides sensitivities for both mechanistic and neural parameters, it establishes a direct connection between neural differential-equation modelling and model-based optimal experimental design. The same information used to select an architecture can therefore help identify experimental conditions that are most informative for refining it.

Uncertainty quantification represents a further advantage of the analytic formulation. Neural components and mechanistic parameters are incorporated into a common parameter vector, allowing confidence intervals to be estimated directly whenever supported by the data. This capability is particularly relevant for HNODE and UDE models, where overparameterised neural components, non-linear dynamics, symmetries and multimodal objective functions complicate parameter-level uncertainty analysis \citep{schmid2025assessment}. Across the benchmark systems, singular or ill-conditioned Fisher information matrices revealed models whose parameters could not be reliably constrained despite apparently satisfactory trajectory fits. Incorporating this information into model selection prevented predictive accuracy from masking uncertainty and redundancy.

Explicit uncertainty estimates also provide a principled basis for model pruning. Conventional pruning methods primarily target network compression, computational efficiency or memory reduction and typically rely on parameter magnitude, sparsity or saliency \citep{Molchanov2017VariationalDropout,Frankle2019LotteryTicket,Cheng2024PruningSurvey}. In the iNODE workflow, pruning instead follows calibration and identifiability analysis. Parameters are removed when they are weakly supported by the observations and contribute negligibly to the inferred dynamics. Pruning therefore becomes an inferential operation that improves conditioning and parsimony while preserving the dynamic behaviour supported by the data.

The ability to infer latent dynamics is particularly important in systems where key processes cannot be observed directly but nevertheless determine collective behaviour. In microbial communities, recovering time-varying interaction coefficients can reveal how ecological relationships shape coexistence, succession and stability, while supporting the rational design of defined consortia with desired biotechnological functions \citep{venturelli2018deciphering,clark2021design}. In epidemic models, reconstructing a time-varying transmission rate can reveal changes in transmissibility and inform the adaptation of public-health interventions \citep{thompson2019improved}. Reliable latent-function inference therefore extends neural differential equations beyond trajectory prediction towards interpretable and decision-relevant dynamic modelling.

These results support a broader perspective on neural differential equations. Their scientific value depends not only on their ability to fit complex trajectories, but also on whether the inferred representation is constrained by experimental information and remains reliable under changes in initial conditions, inputs or observation regimes. Identifiability provides a common language for assessing these requirements in both fully data-driven and mechanistic–neural models, linking architecture design, calibration, uncertainty quantification, pruning and experimental design within a unified framework. By combining analytic embedding, global calibration, explicit sensitivity analysis and identifiability-aware architecture selection, iNODE elevates identifiability from a diagnostic property to a model-design principle and provides a principled route towards neural differential-equation models that are accurate, parsimonious and reliable beyond interpolation. More broadly, the framework connects the information content of experimental data to the complexity of neural dynamic models, supporting their use for mechanistic inference and the reliable reconstruction of otherwise inaccessible dynamic processes.
\section*{Methods}

\subsection*{General iNODE formulation}

We consider dynamic systems described by ordinary differential equations of the form
\begin{equation}
\frac{dy}{dt}(t)=\tilde f\big(y(t),\boldsymbol{\theta},t\big), \qquad y(t_0)=y_0,
\label{eq:general_inode}
\end{equation}
where \(y(t)\in\mathbb{R}^n\) denotes the vector of state variables and \(\boldsymbol{\theta}=[\boldsymbol{\psi}_M,\boldsymbol{\phi}_{NN}]\) collects both the mechanistic parameters \(\boldsymbol{\psi}_M\) and the neural-network parameters \(\boldsymbol{\phi}_{NN}\). In iNODE, neural components are written explicitly as analytical functions and embedded directly within the governing equations, so that both mechanistic and neural parameters appear in the same ODE system and can be analysed within a unified calibration and uncertainty-quantification framework. Depending on the application, the neural component may represent fully data-driven dynamics, latent time-varying parameters, or unknown terms that complement a mechanistic model. The specific formulations used in the case studies are provided in the Supplementary Information, Section~S.1.

\subsection*{Pseudo-data generation}

Pseudo-experimental data were generated in AMIGO2 \citep{balsa-etal:2016AMIGO2} by simulating each ground-truth dynamic model with predefined nominal parameter values. The resulting trajectories were sampled at the time points specified for each benchmark to obtain synthetic observations. For a measured observable \(y_o(t)\), the pseudo-data at sampling time \(t_s\) were generated as
\begin{equation}
y_{s,o}^{\mathrm{meas}} =
y_o(t_s;\boldsymbol{\theta}^{\mathrm{true}}) + \varepsilon_{s,o},
\qquad
\varepsilon_{s,o} \sim \mathcal{N}(0,\sigma_{s,o}^2),
\label{eq:pseudodata}
\end{equation}
where \(y_o(t_s;\boldsymbol{\theta}^{\mathrm{true}})\) denotes the simulated ground-truth output and \(\varepsilon_{s,o}\) is an additive Gaussian measurement error. Noise terms were sampled independently across observables, experiments and sampling times. For the noise-free oscillator benchmarks, \(\sigma_{s,o}=0\). For the Lotka--Volterra and SIR benchmarks, Gaussian noise was added to emulate more realistic experimental conditions. The initial conditions, sampling schemes and noise levels used in each case study are reported in the Supplementary Information.

\subsection*{Model calibration with direct ODE simulation}

To assess the framework, we considered four benchmark case studies based on pseudo-experimental data generated from ground-truth dynamic systems. This setup provides complete knowledge of the true model structure and parameter values, enabling the evaluation of predictive accuracy, parameter recovery, uncertainty quantification, and identifiability. Depending on the case study, pseudo-data were generated under ideal or realistic conditions, including sparse sampling, partial observability, and additive measurement noise.

We considered an experimental scheme comprising \(n_e\) experiments, which may differ in initial conditions, duration, and sampling times. For experiment \(e\), let \(n_o^e\) denote the number of measured observables, and let \(n_s^{e,o}\) be the number of sampling times associated with observable \(o\). The corresponding discrete measurements are represented by \(y_{s}^{e,o,\mathrm{meas}}\), and the model predictions by \(y_{s}^{e,o,\mathrm{pred}}(\boldsymbol{\theta})\).

Model calibration was formulated as a nonlinear parameter-estimation problem in which the mechanistic and neural parameters are jointly estimated by minimising a weighted least-squares objective:
\begin{equation}
J(\boldsymbol{\theta})=
\sum_{e=1}^{n_e}\sum_{o=1}^{n_o^e}\sum_{s=1}^{n_s^{e,o}}
\left[
\frac{y_{s}^{e,o,\mathrm{pred}}(\boldsymbol{\theta})-y_{s}^{e,o,\mathrm{meas}}}
{\sigma_{s}^{e,o}}
\right]^2,
\label{eq:wlsq}
\end{equation}
where \(\sigma_{s}^{e,o}\) denotes the standard deviation associated with each observation. Under the assumption of Gaussian measurement noise with known or constant variance, this objective corresponds to maximum-likelihood estimation \cite{walter-pronzato:97}. The estimation of parameters was subject to the dynamics of the system in Eq.~\ref{eq:general_inode} and to parameter limits that define the admissible search space \(\theta_L \le \theta \le \theta_U\).

\subsection*{Sensitivity analysis}

Sensitivity analysis was used to quantify how changes in the estimated parameters affect the model outputs over the experimental conditions considered. For experiment \(e\), observable \(o\), sampling time \(t_s^{e,o}\), and parameter \(\theta_p\), the local parametric sensitivity was defined as
\begin{equation}
S_p^{e,o}(t_s^{e,o})=
\frac{\partial y^{e,o}}{\partial \theta_p}(t_s^{e,o}),
\qquad p=1,\dots,n_\theta .
\label{eq:sens}
\end{equation}
Because the neural component is analytically embedded into the ODE system, these sensitivities can be computed for both mechanistic and neural parameters within the same simulation-estimation framework.

The sensitivities were assembled into the sensitivity matrix
\begin{equation}
S_{ij} =
\frac{\partial y_i}{\partial \theta_j},
\qquad i=1,\dots,N,\quad j=1,\dots,n_\theta,
\label{eq:sens_matrix}
\end{equation}
where \(N\) is the total number of measured data points and \(n_\theta\) is the number of estimated parameters. The rank of this matrix provides information on the number of independent parameter directions that are locally informed by the data. If
\begin{equation}
\mathrm{rank}(S) < n_\theta,
\label{eq:rank_def}
\end{equation}
then at least some parameter directions are locally linearly dependent or indistinguishable in the neighbourhood of the estimated solution. This rank deficiency indicates that the available data do not provide sufficient independent information to identify all parameters locally.

In addition to the rank analysis, we used the $\delta^{msqr}$ metric to summarise the relative influence of individual parameters on the model outputs, particularly in large conventional NODE and HNODE models where visual inspection of all sensitivity profiles is impractical. For each parameter, $\delta^{msqr}$ was computed from the local parametric sensitivities as
\begin{equation}
\delta^{\mathrm{msqr}}_p =
\frac{1}{N}
\sqrt{
\sum_{e,o,s}
\left(
\frac{\partial y^{e,o}}
{\partial \theta_p}
\left(t_s^{e,o}\right)
\right)^2
},
\label{eq:dmsqr}
\end{equation}
where \(N\) is the total number of measured data points. This quantity provides a scaled aggregate measure of the magnitude of the output sensitivities associated with parameter \(\theta_p\). Parameters with larger \(\mathrm{dmsqr}\) values exert a stronger local influence on the simulated outputs under the experimental conditions considered, whereas parameters with very small values have little effect on the model trajectories.

\subsection*{Practical identifiability analysis}

Practical identifiability analysis was used to assess whether the available data contain sufficient information to estimate model parameters with finite and reliable uncertainty in the presence of experimental noise. Parameter uncertainty was quantified locally using the Fisher Information Matrix (FIM), computed from the parametric sensitivities evaluated at the estimated parameter vector \(\boldsymbol{\theta}^*\).

Under the weighted least-squares formulation used for model calibration, the FIM can be approximated as
\begin{equation}
F \approx S^\top W S,
\label{eq:fim_sens}
\end{equation}
where \(S\) is the sensitivity matrix and \(W\) is a diagonal matrix containing the observation weights. The FIM provides a local measure of how strongly the parameters are constrained by the available data: well-informed parameters increase the information content of the experiment, whereas weakly sensitive or highly correlated parameters lead to poor conditioning and large uncertainty.

When the FIM is non-singular and sufficiently well conditioned, parameter confidence intervals were derived from the covariance matrix approximation given by the Cramér--Rao bound,
\begin{equation}
C \ge F^{-1}(\boldsymbol{\theta}^*),
\label{eq:crb}
\end{equation}
where \(C\) denotes the parameter covariance matrix. For parameter \(\theta_i\), the corresponding confidence interval half-width was computed as
\begin{equation}
\rho_i =
t_{\alpha/2}^{\gamma}\sqrt{C_{ii}},
\label{eq:ci}
\end{equation}
where \(t_{\alpha/2}^{\gamma}\) is the Student's \(t\)-value for confidence level \((1-\alpha)\times100\%\), and \(\gamma\) denotes the degrees of freedom. Relative confidence intervals were reported as
\begin{equation}
\mathrm{CI}_i(\%) =
100\,\frac{\rho_i}{|\theta_i^*|}.
\label{eq:relative_ci}
\end{equation}

When the FIM was singular or nearly singular, confidence intervals could not be reliably estimated from its inverse. In these cases, the model was considered locally practically non-identifiable under the corresponding experimental design, and the source of non-identifiability was further examined through the sensitivity matrix, parameter correlations and \(\mathrm{dmsqr}\) values. Rank-deficient sensitivity matrices indicate that some parameter combinations are redundant or poorly informed by the data, whereas strong parameter correlations indicate compensatory relationships among parameters.  In addition, confidence intervals were computed using a bootstrap resampling method \cite{joshi-seidel-morgenstern-kremling:2006, Balsa-Canto-etal:2010}. The underlying idea is to simulate the feasibility of performing hundreds of repetitions of the same experimental scheme, given a specified experimental error. The model calibration problem is solved for each data realisation, and the cloud of solutions is recorded in a matrix. The confidence intervals obtained from this analysis were subsequently used to generate the parameter ensembles and corresponding prediction uncertainty bands.

\subsection*{Identifiability-aware neural architecture search}

For each case study, candidate neural architectures were generated by varying the number of hidden layers, the number of neurons per layer, and the activation functions. To reduce overparameterisation, we imposed a preliminary data-support constraint based on the total number of available observations. Specifically, only architectures that met a minimum ratio of 4 observations per trainable parameter were retained. The rationale is that architectures with too many free parameters relative to the number of measurements are more likely to be underdetermined, highly correlated, and weakly identifiable, even before formal calibration. This rule was used as a pragmatic first filter to exclude clearly underdetermined neural parameterisations during architecture generation; then practical identifiability was assessed explicitly after calibration.

All admissible candidate architectures were calibrated and compared using complementary criteria reflecting predictive accuracy, model complexity, and parameter uncertainty. Candidate models were evaluated sequentially, from simpler to more complex structures, and ranked using the following metrics.

\textbf{Normalised root-mean-square error.} Predictive accuracy was quantified through the normalised root-mean-square error (NRMSE), computed across all experiments, observables, and sampling times as
\begin{equation}
\mathrm{NRMSE}=
\frac{
\sqrt{
\frac{1}{N}
\sum_{e=1}^{n_e}\sum_{o=1}^{n_o^e}\sum_{s=1}^{n_s^{e,o}}
\left(y_{s}^{e,o,\mathrm{pred}}-y_{s}^{e,o,\mathrm{meas}}\right)^2
}
}{
y_{\max}^{\mathrm{meas}}-y_{\min}^{\mathrm{meas}}
},
\label{eq:nrmse}
\end{equation}
where \(N\) denotes the total number of measurements and the denominator corresponds to the range of the measured data.

\textbf{Akaike Information Criterion.} To compare candidate models with different numbers of parameters, we used the Akaike Information Criterion (AIC, \cite{BanksJoyner2017}), which balances goodness of fit
against model complexity: 
\begin{equation}
\mathrm{AIC}=N\ln\left(\frac{J^*}{N}\right)+2(n_\theta+1),
\label{eq:aic}
\end{equation}
where \(J^*\) is the minimum value of the weighted least-squares objective and \(n_\theta\) is the number of estimated parameters. Lower AIC values indicate a more favourable compromise between predictive accuracy and parameter parsimony, and were therefore used to penalise unnecessarily complex architectures during model selection.

\textbf{Model Quality Ranking Index.} To combine predictive accuracy, parsimony, and practical identifiability in a single criterion, we introduced the Model Quality Ranking Index (MQRI). For each candidate model, the values of NRMSE, AIC, maximum confidence-interval width, and median confidence-interval width were ranked throughout the model candidate set, assigning rank 1 to the best-performing value in each category. The MQRI of model \(m_i\) was then defined as
\begin{equation}
\mathrm{MQRI}(m_i)=\sum_{j=1}^{4}R_i^{(j)}.
\label{eq:mqri}
\end{equation}
Lower MQRI values identify architectures that achieve the most favourable balance between predictive accuracy, model complexity, and parameter uncertainty, and were therefore used to guide automated architecture selection.

Following architecture selection, in some case studies, an optional identifiability-guided pruning step was applied to reduce redundancy and improve parameter robustness. Parameters with negligible estimated contributions and large relative confidence intervals were removed, and the remaining parameters were re-estimated.

\subsection*{Uncertainty quantification}

To visualise how parameter uncertainty propagates to model predictions, we employed an ensemble-based uncertainty quantification approach. Starting from the calibrated parameter vector $\hat{\boldsymbol{\theta}}$, repeated simulations were performed using perturbed parameter sets sampled within the corresponding confidence intervals, generating an ensemble of model trajectories from which prediction bands were derived. A total of $10 \times n_{\theta}$ simulations were carried out to propagate parameter uncertainty through all training and validation experiments.

Perturbed parameter vectors were generated by sampling uniformly within the $95\%$ confidence interval of each parameter. For a parameter $\theta_i$ with estimated standard deviation $\sigma_i$, the interval was defined as

\[
\left[\theta_i - 1.96\,\sigma_i,\;\; \theta_i + 1.96\,\sigma_i\right].
\]

At each Monte Carlo iteration, a new parameter vector was obtained by drawing independent uniform samples inside these intervals, ensuring exploration of all values compatible with the confidence bounds without assuming any particular distributional shape. For mechanistic parameters that must remain strictly positive, the lower bound was truncated to prevent sampling values that would be physically or biologically meaningless.

Each sampled parameter vector was then used to simulate all training and validation experiments, producing ensembles of predicted trajectories. Collecting these trajectories across all Monte Carlo samples yields, for each observable and each time point, an empirical distribution of predicted values.

The uncertainty bands shown in the figures correspond to the interquartile range of the simulated trajectories, namely the region between the 25th and 75th percentiles at each time point, while the central line represents the median prediction. These percentiles were computed directly from the ensemble of simulated trajectories. This provides a simple and consistent way to visualise how parameter uncertainty propagates into model predictions, without imposing any parametric assumptions on the distribution of the outputs.

\subsection*{Numerical methods and computational tools for iNODE}

The iNODE models were calibrated using the AMIGO2 toolbox \cite{balsa-etal:2016AMIGO2}, which provides a unified environment for dynamic simulation, parameter estimation, sensitivity analysis, identifiability assessment and uncertainty quantification. Ordinary differential equations were solved using the CVODES solver in the SUNDIALS suite \citep{hindmarsh2005sundials}, which provides adaptive time integration and forward sensitivity computation. The numerical accuracy was controlled through user-defined relative and absolute tolerances. Model calibration was performed using global optimisation via the enhanced scatter search algorithm (eSS, \citep{Egea2009eSS}).

\subsection*{Numerical methods and computational tools for conventional NODE and Hybrid NODE}

Conventional NODE and HNODE baselines were implemented in Python using back-propagation through the ODE solver and local gradient-based optimisation \citep{Rumelhart1986,chen2018neural,KingmaBa2015Adam}. Depending on the case study, stabilisation strategies included Xavier-type parameter initialisation \citep{GlorotBengio2010}, gradient-norm clipping \citep{Pascanu2013RNN}, adaptive learning-rate scheduling, and multi-shooting strategies \citep{Bock1984} to improve numerical robustness over long training horizons.

PyTorch \citep{Paszke2019PyTorch} served as the core computational framework for all models, providing GPU-accelerated tensor operations, automatic differentiation, and modular neural network construction. Neural components were implemented as fully connected multilayer perceptrons (MLPs) with smooth activation functions, such as $\tanh$, which are well suited for continuous-time dynamics. Weights were initialised using Xavier-type schemes to promote numerical stability, and model optimisation relied primarily on the Adam optimiser, combined with gradient-norm clipping and adaptive learning-rate schedulers to ensure stable convergence.

Across all formulations, system dynamics were expressed as ordinary differential equations whose vector fields depend on the states and time. In conventional NODE, the full vector field was learnt directly by a neural network from data. In hybrid NODE, neural modules were embedded within explicitly defined mechanistic equations, capturing unknown or time-dependent components—such as interaction functions or latent dynamics—while preserving the known physical or biological structure. Mechanistic and neural parameters were optimised jointly, with structural constraints such as positivity enforced through smooth transformations (e.g., Softplus or exponential mappings).

Time integration was performed using a custom fourth-order Runge–Kutta (RK4) solver implemented entirely in PyTorch, ensuring full differentiability of the numerical scheme. Fixed or adaptive internal substepping was used to balance numerical accuracy and computational efficiency, depending on the temporal resolution of the data and the training stage. In computationally demanding settings, neural evaluations were vectorised and reused across substeps to reduce overhead, and state clipping was applied in selected models to prevent numerical divergence.

In addition to these strategies, alternative training designs were explored. In particular, trajectory subsampling was applied in selected formulations, whereby a different subset of time points was used at each optimisation step. This results in a sequence of stochastic training problems—typically one subsampled trajectory per epoch—which introduces variability in the gradients and perturbs the optimisation process without requiring explicit multi-start strategies. This mechanism acts as an implicit regularisation and can help mitigate overfitting to specific trajectories while improving exploration of the parameter space. In some configurations, multiple-shooting was also incorporated by periodically reintegrating short trajectory segments from intermediate states, introducing local consistency constraints that further stabilise long-horizon training. When multiple experimental trajectories were available, they were fitted simultaneously, enforcing global consistency across different initial conditions and acting as an additional form of structural regularisation.

Model training minimised the mean squared error between simulated trajectories and experimental observations, with the loss weighted by the temporal increments $\Delta t$ to account for non-uniform sampling. In addition to the training objective used for backpropagation, model performance was quantified using normalised root mean square error (NRMSE) metrics for trajectory reconstruction and, when applicable, for the recovery of time-dependent parameters against theoretical references.

\section*{Funding}
This work has received funding from MICIU/AEI/10.13039/501100011033 and FEDER, UE Grant numbers: CPP2024-011608 and PID2024-161806OB-C22 and from Xunta de Galicia grant IN607B.

\section*{Data availability}
The synthetic datasets generated and analysed in this study are available at: \\
\href{}{https://doi.org/10.5281/zenodo.21526685}

\section*{Code availability}
The code used to generate the models, simulations and analyses is available at:\\
\href{}{https://doi.org/10.5281/zenodo.21526685}

\section*{Competing interests}
The authors declare no competing interests.

\section*{Author contributions}
N.C.-M. implemented and performed the computations, contributed to the analysis and interpretation of the results, prepared visualisations, and drafted the manuscript.
E.B.-C. conceived and designed the study, contributed to the analysis and interpretation of the results, prepared visualisations, revised and edited the manuscript, and acquired funding.

\bibliographystyle{unsrt}
\bibliography{biblio}

\end{document}


\section*{Supplementary Information File}

\section*{Identifiability-aware neural ordinary differential equations for parsimonious and reliable dynamic modelling}

\noindent \textbf{Authors:} N. Campo-Manzanares, E. Balsa-Canto

\vspace*{0.25cm}

This Supplementary Information provides the mathematical details and extended results supporting the main manuscript. It includes the general model formulations considered in this work, a summary of the benchmark case studies, and additional analyses of calibration, identifiability, uncertainty, and model comparison that complement the results presented in the main text.

\tableofcontents
\listoffigures
\listoftables
\clearpage

\section*{S.1 General model formulation}
\addcontentsline{toc}{section}{S.1 General model formulation} 
We consider a unified formulation for modelling dynamical systems based on ordinary differential equations (ODE) that systematically integrates prior mechanistic knowledge with data-driven components. This formulation subsumes purely mechanistic models, fully data-driven Neural ODE, and intermediate hybrid approaches within a single mathematical framework. Specifically, we define the following general model form:

\begin{equation}
\left\{
\begin{aligned}
\frac{d y}{d t}(t)
&=
f\!\left(
c_1\, (y(t),\,\boldsymbol{\psi}_M),\,
c_2\, w(t),\,
c_3\, \mathcal{NN}\big(z(t);\boldsymbol{\phi}_{NN}\big),\,
t
\right)\\
w(t)
&=
\mathcal{NN}_1\big(z_1(t);\boldsymbol{\phi}_{NN}\big)\\
y(t_0) &= y_0
\end{aligned}
\right.
\end{equation}

where $y(t)\in\mathbb{R}^n$ denotes the system state vector and $\boldsymbol{\psi}_M$ the set of mechanistic parameters. The latent variable $w(t)$ represents a time-dependent parameter inferred from data through the neural network $\mathcal{NN}_1(\cdot;\boldsymbol{\phi}_{NN})$, parametrised by weights and biases $\boldsymbol{\phi}_{NN}$ and driven by inputs $z_1(t)$. The function $\mathcal{NN}(\cdot;\boldsymbol{\phi}_{NN})$ denotes an additional neural component embedded directly in the system dynamics to account for unknown or unmodelled effects. The inputs $z(t)$ and $z_1(t)$ may depend on time, the state variables, or both.

The binary coefficients $c_1,c_2,c_3\in\{0,1\}$ allow for selective activation of the mechanistic component, the latent parameter dependent on time, and the embedded neural correction, respectively. By appropriately choosing these coefficients, the formulation recovers a broad class of mechanistic, hybrid, and fully data-driven models as special cases.

\subsection*{Mechanistic Models}
\addcontentsline{toc}{subsection}{Mechanistic Models} 

Classical mechanistic models are recovered by activating only the mechanistic component, that is,

\[
c_1=1,\qquad c_2=0,\qquad c_3=0.
\]

The resulting system reads:

\begin{equation}
\frac{d y}{d t}(t)
=
f\!\left(
y(t),\,
\boldsymbol{\psi}_M,\,
t
\right),
\qquad
y(t_0)=y_0
\end{equation}

In this setting, the system dynamics are fully prescribed by prior domain knowledge, such as physical laws, conservation principles, or established phenomenological relationships. All interactions and causal dependencies are explicitly encoded in the functional form of $f$.


\subsection*{Latent Parameter Dynamic Model}
\addcontentsline{toc}{subsection}{Latent Parameter Dynamic Model} 
The first hybrid extension allows selected parameters to vary over time, while preserving the mechanistic structure of the governing equations. This is obtained by choosing

\[
c_1=1,\qquad c_2=1,\qquad c_3=0,
\]

which yields the latent parameter dynamics formulation:

\begin{equation}
\left\{
\begin{aligned}
\frac{d y}{d t}(t)
&=
f\!\left(
y(t),\,\boldsymbol{\psi}_M,\,
w(t),\,
t
\right)\\
w(t)
&=
\mathcal{NN}_1\big(z_1(t);\boldsymbol{\phi}_{NN}\big)\\
y(t_0)&=y_0
\end{aligned}
\right.
\end{equation}

Here, the neural network $\mathcal{NN}_1$ represents unknown, time-varying, or poorly characterised parameters, while the mechanistic backbone of the model is retained. 

This formulation offers a good compromise between flexibility and interpretability and is particularly appropriate when parameters are expected to evolve over time, but no explicit evolution law is available. 

\subsection*{Hybrid Neural ODE}
\addcontentsline{toc}{subsection}{Hybrid Neural ODE} 
An alternative hybrid strategy embeds neural networks directly into the system dynamics to capture unresolved mechanisms or structural model mismatch. Setting

\[
c_1=1,\qquad c_2=0,\qquad c_3=1,
\]

the model becomes:

\begin{equation}
\frac{d y}{d t}(t)
=
f\!\left(
y(t),\,
\boldsymbol{\psi}_M,\,
\mathcal{NN}\big(z(t);\boldsymbol{\phi}_{NN}\big),\,
t
\right),
\qquad
y(t_0)=y_0
\end{equation}

In this formulation, mechanistic knowledge provides the structural backbone, while the neural component serves as a flexible correction term that accounts for partially unknown or neglected dynamics. Models of this type are commonly referred to as Mechanistic Neural ODE, Universal Differential Equations (UDE), or Hybrid Neural ODE (HNODE).

This hybridisation can improve robustness and generalisation—particularly in data-scarce regimes—by combining mechanistic inductive bias with neural expressiveness. 

\subsection*{Neural ODE}
\addcontentsline{toc}{subsection}{Neural ODE} 

Finally, the unified framework includes purely data-driven approaches. By selecting
\[
c_1=0,\qquad c_2=0,\qquad c_3=1,
\]

the system reduces to standard Neural ODE:

\begin{equation}
\frac{d y}{d t}(t)
=
\mathcal{NN}\big(z(t);\boldsymbol{\phi}_{NN}\big),
\qquad
y(t_0)=y_0
\end{equation}

Here, the neural network directly parameterises the full state dynamics without any explicit mechanistic structure. This approach offers high expressive capacity and minimal reliance on prior assumptions, allowing it to capture complex nonlinear dynamics given sufficient data. 

\pagebreak
\section*{S.2 Benchmarking framework}
\addcontentsline{toc}{section}{S.2 Benchmarking framework} 

We evaluate the iNODE framework using four controlled case studies with known generating dynamics. The benchmarks span three established dynamic systems and vary systematically in mechanistic knowledge, data availability, observability, sampling density, measurement noise and model complexity.

Each case study addresses a distinct inference problem, ranging from fully data-driven vector-field recovery to mechanistic-neural closure and the reconstruction of latent time-varying parameters under full or partial observability. Because the states, parameters and latent functions of the generating systems are known, the benchmarks enable direct assessment of prediction accuracy, dynamic-structure recovery, practical identifiability, uncertainty and extrapolation. This distinction is essential because accurate reproduction of observed trajectories does not by itself establish recovery of the underlying dynamics.

In each case, the iNODE workflow is compared with the corresponding conventional NODE or hybrid NODE workflow using the same training data and evaluation scenarios. Table~\ref{tab:case_studies_summary} summarises the modelling challenge, experimental design and data characteristics of the four benchmarks.

\begin{table}[h!]
\centering
\caption[Benchmark case studies used to evaluate the iNODE framework.]{
\textbf{Benchmark case studies used to evaluate the iNODE framework.}
The four case studies span fully data-driven and hybrid neural ODE formulations, latent time-varying parameters, and varying levels of observability, sampling density and measurement noise. Together, they provide controlled ground-truth systems for assessing prediction accuracy, parameter identifiability, latent-function recovery and extrapolation.
}
\vspace{6pt}
\label{tab:case_studies_summary}
\begin{tabular}{p{3.5cm} p{5.6cm} p{6.8cm}}
\toprule
\textbf{Case study} & \textbf{Type of scenario} & \textbf{Type of data} \\
\midrule 
\texttt{Data-CDO: }Fully data-driven cubic damped oscillator
& Neural ODE without mechanistic prior knowledge; benchmark for extrapolation and phase-space recovery 
& Synthetic pseudo-experimental data from a single trajectory; both states observed; densely sampled; noise-free \\
\addlinespace
\texttt{Hybrid-CDO: }Hybrid mechanistic-neural cubic damped oscillator  
& Hybrid / mechanistic Neural ODE with known equation for one state and neural representation of the other 
& Synthetic pseudo-experimental data; both states observed; one- and two-experiment designs; densely sampled; noise-free \\
\addlinespace
\texttt{tLV: }
Lotka-Volterra system with time-varying interactions 
& Hybrid latent-parameter dynamics; fully observable system with neural representation of time-varying interaction coefficients 
& Synthetic pseudo-experimental data from multiple experiments; both species observed; sparse sampling; analyses under realistic noisy conditions (10\% Gaussian noise) \\
\addlinespace
\texttt{tSIR: }
SIR model with time-varying transmission rate 
& Hybrid latent-parameter dynamics under partial observability; neural representation of the time-varying transmission rate 
& Synthetic pseudo-experimental data; only infected population observed; single- and multi-experiment designs; sparse sampling; noise-free observations \\
\bottomrule
\end{tabular}
\end{table}

Table~\ref{tab:supplementary_findings_summary} summarises the main findings obtained across the four benchmark case studies. For each scenario, it highlights the modelling challenge, the comparison performed between conventional NODE/HNODE and iNODE approaches, and the key conclusions regarding dynamic recovery, predictive performance and practical identifiability.

\begin{table}[H]
\caption[Summary of Supplementary findings across benchmark case studies.]{
\textbf{Summary of Supplementary findings across benchmark case studies.}
The table summarises the modelling challenge, the main comparison and the key finding for each benchmark, including the limitations observed for conventional NODE or HNODE workflows.
}
\vspace{6pt}
\label{tab:supplementary_findings_summary}
\begin{tabular}{p{2cm}p{4.0cm}p{3.6cm}p{5.6cm}}
\toprule
Case study & Main challenge & Main comparison & Main finding \\
\midrule
\texttt{Data-CDO} &
Fully data-driven vector-field recovery &
NODE vs iNODE &
The conventional NODE fits the training trajectory but has a rank-deficient sensitivity matrix, indicating practical non-identifiability. A compact iNODE improves phase-space recovery, extrapolation and parameter identifiability. \\

\texttt{Hybrid-CDO} &
Hybrid closure and role of experimental diversity &
HNODE vs hybrid iNODE; 1exp vs 2exp &
The conventional HNODE requires a much larger architecture and has a rank-deficient sensitivity matrix, indicating practical non-identifiability. A compact hybrid iNODE recovers the system dynamics, while experimental diversity and identifiability-guided pruning reduce uncertainty and parameter correlations. \\

\texttt{tLV} &
Noisy, sparse data and recovery of latent interactions &
HNODE vs hybrid iNODE &
Conventional HNODE models can fit the observed trajectories but provide weaker latent-interaction recovery and higher uncertainty. The selected hybrid iNODE recovers the time-varying interactions with reduced uncertainty. \\

\texttt{tSIR} &
Partial observability and recovery of the latent transmission rate &
1exp vs extended 1exp vs 2exp; HNODE vs hybrid iNODE &
Single-experiment configurations remain weakly informative for the latent transmission rate. The hybrid iNODE recovers the latent transmission dynamics, with the two-experiment design providing the strongest improvement in uncertainty and parameter identifiability. \\
\bottomrule
\end{tabular}
\end{table}

\section*{S.3~\texttt{Data-CDO: }Fully data-driven cubic damped oscillator}
\addcontentsline{toc}{section}{S.3~\texttt{Data-CDO: }Fully data-driven cubic damped oscillator} 
The cubic damped oscillator is a nonlinear dynamic system that has been extensively used in physics and engineering as a benchmark for oscillatory dynamics with dissipation and nonlinear restoring forces. Its combination of a compact mathematical formulation and rich nonlinear behaviour makes it well suited for evaluating analytical, numerical, data-driven, and hybrid modelling approaches. In this work, the mechanistic model is used as a ground-truth system to generate dense synthetic pseudo-experimental data, which are then used to assess whether neural differential equation frameworks can recover the underlying dynamics from observations alone. The system is formulated in terms of two interacting state variables and is governed by the following nonlinear ordinary differential equations:

\begin{equation}
\left\{
\begin{aligned}
\frac{d x_1}{d t} &= -0.1\, x_1^{3} + 2.0\, x_2^{3}\\
\frac{d x_2}{d t} &= -2.0\, x_1^{3} - 0.1\, x_2^{3}
\end{aligned}
\right.
\end{equation}

Here, \(x_1(t)\) and \(x_2(t)\) denote the system states. This formulation preserves the characteristic cubic nonlinearities of damped oscillatory systems while introducing cross-coupled interactions between the states. Its nonlinear nature can lead to strong sensitivity to initial conditions and makes accurate recovery of the governing vector field non-trivial. These features make the cubic damped oscillator a useful benchmark for analysing the trade-offs between model flexibility, parameter identifiability, and predictive robustness in neural differential equation frameworks.

Synthetic data were generated following the benchmark configuration proposed by Goyal and Benner~\cite{goyal2023neural}. The dataset is noise-free and densely sampled, thereby providing a favourable setting for evaluating model expressiveness, calibration stability, and practical identifiability under idealised conditions. Table~\ref{tab:data_CDO_scenarios} summarises the datasets used for both training and validation analyses.

\begin{table}[h!]
\centering
\caption[\texttt{Data-CDO: }Training and validation scenarios for the fully data-driven cubic damped oscillator.]{
\textbf{\texttt{Data-CDO: }Training and validation scenarios for the fully data-driven cubic damped oscillator.}
The scenarios use distinct initial conditions to assess whether models fitted to one trajectory recover the underlying phase-space structure and generalise to an unseen trajectory.}
\vspace{6pt}
\label{tab:data_CDO_scenarios}
\begin{tabular}{ccccc}
\hline
\textbf{Dataset} & \textbf{Scenario} & \textbf{Experiment(s)}                & \textbf{Measurements} & \textbf{Total data} \\
 &  &  &  \textbf{per state} & \textbf{points}\\ \hline
Train            & \texttt{1exp}             & \((x^{(1)}_1(0),x^{(1)}_2(0))=(2,0)\) & 1250                            & 2500                       \\[3mm]
Validation     & \texttt{1exp}             & \((x^{(1)}_1(0),x^{(1)}_2(0))=(1,1)\) & 500                             & 1000                       \\ \hline
\end{tabular}
\end{table}

\subsection*{S.3.1 \texttt{Data-CDO: }Implementation with the conventional NODE framework}
\addcontentsline{toc}{subsection}{S.3.1 \texttt{Data-CDO: }Implementation with the conventional NODE framework} 
We first consider a fully data-driven modelling approach based on conventional Neural Ordinary Differential Equations (NODE). In this formulation, the system dynamics are represented entirely by a neural network, without explicitly incorporating mechanistic structure or prior physical knowledge. Accordingly, the cubic damped oscillator is modelled as

\begin{equation}
\left\{
\begin{aligned}
\frac{d x_1}{d t} &= \mathcal{NN}\big(x_1(t), x_2(t); \boldsymbol{\phi}_{\mathrm{NN}}\big)\\
\frac{d x_2}{d t} &= \mathcal{NN}\big(x_1(t), x_2(t); \boldsymbol{\phi}_{\mathrm{NN}}\big)
\end{aligned}
\right.
\end{equation}

where \(\mathcal{NN}(\cdot;\boldsymbol{\phi}_{\mathrm{NN}})\) denotes a neural network parameterised by weights and biases \(\boldsymbol{\phi}_{\mathrm{NN}}\). The network learns the vector field directly from the data, providing a flexible but structurally unconstrained representation.

Synthetic data were generated following the benchmark configuration proposed by Goyal and Benner~\cite{goyal2023neural}. A single training experiment with initial conditions \((x^{(1)}_1(0),x^{(1)}_2(0))=(2,0)\,\mathrm{m}\) was simulated, yielding 1250 measurements per state variable (2500 data points in total). The data are noise-free and densely sampled, providing a favourable setting for evaluating model expressiveness, calibration stability, and practical identifiability under idealised conditions.

The reference architecture considered here is the one proposed by Goyal and Benner~\citep{goyal2023neural}, namely a fully connected neural network with four hidden layers of twenty neurons each: MLP(20,20,20,20). We used hyperbolic tangent activation functions in the hidden layers and a linear activation function in the output layer.

Training was performed using a Python-based NODE pipeline in which the neural network is embedded within a numerical ODE solver and optimised from trajectory-level discrepancies rather than pointwise derivative information, following the standard NODE framework~\citep{chen2018neural}. 

To improve numerical stability and efficiency, the state variables were normalised, and the learnt vector field was integrated using an adaptive fourth-order Runge--Kutta (RK4) scheme with progressively reduced maximum step size during training. The loss function was defined as a mean squared error weighted by temporal increments to account for non-uniform sampling. 

To improve stability and limit long-term error accumulation, a multi-shooting strategy was periodically applied by reintegrating short trajectory segments from intermediate states~\citep{turan2021multipleshooting}. This was performed every three epochs, resulting in approximately 670 multi-shooting updates over 2000 training epochs. Additionally, trajectory subsampling was applied at each epoch, generating one reduced dataset per optimisation step and yielding a total of 2000 subsampled trajectories during training.

After training, we assessed practical identifiability and uncertainty for the optimal NODE model using the analytical embedding in AMIGO2. The resulting Fisher Information Matrix (FIM) was found to be non-invertible, preventing confidence-interval estimation and indicating severe practical identifiability limitations or redundancy in the neural parameterisation. To further investigate the origin of this non-invertibility, parameter sensitivities were analysed using the \texttt{AMIGO\_LRank} function. Due to the large size of the model, which comprises 1362 parameters, this sensitivity analysis is computationally demanding and numerically challenging.

\begin{figure}[h!]
    \centering
    \includegraphics[width=0.90\textwidth]{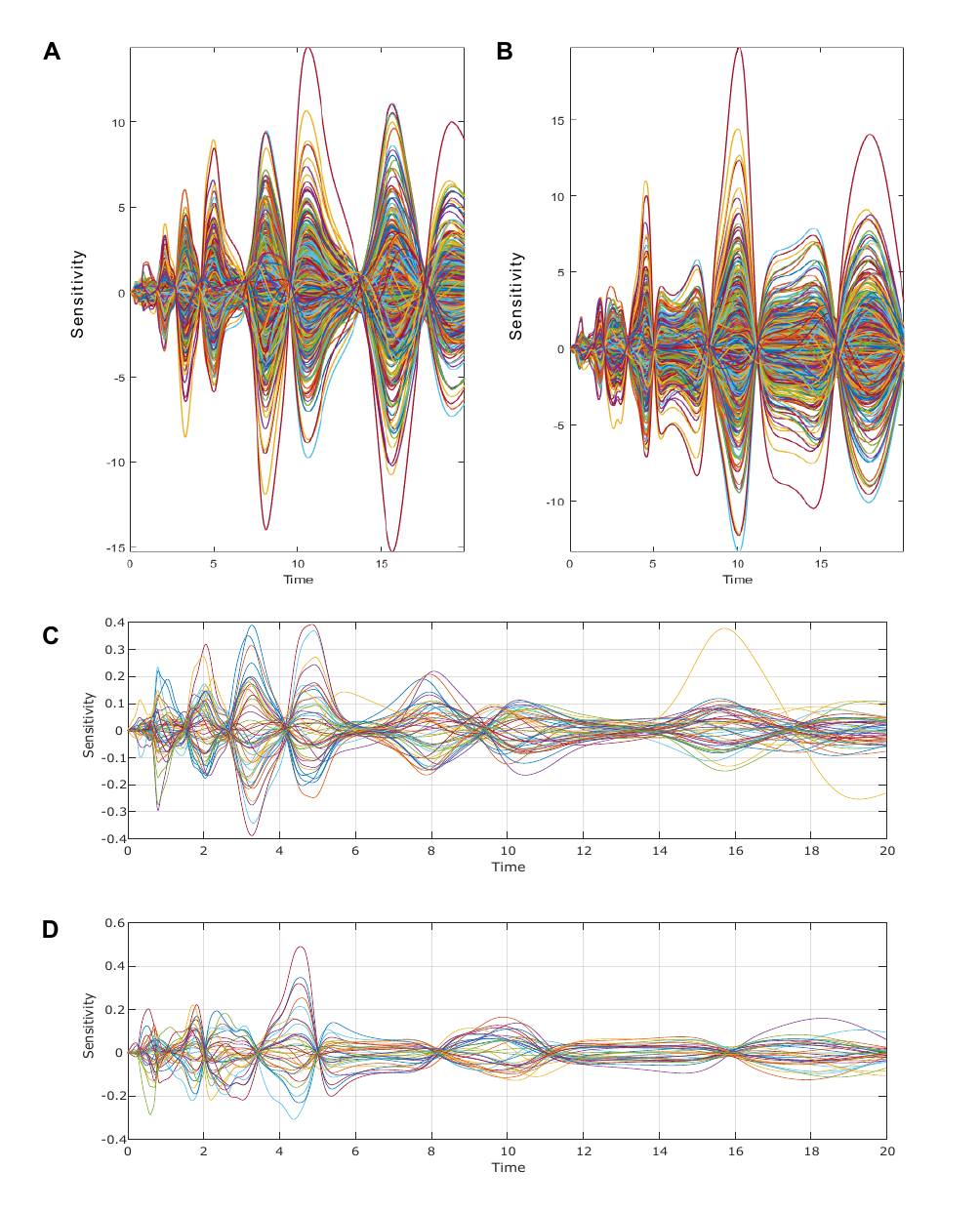}
    \caption[\texttt{Data-CDO: }Sensitivity profiles for the conventional NODE model of the cubic damped oscillator.]{
\textbf{\texttt{Data-CDO: }Sensitivity profiles for the conventional NODE model of the cubic damped oscillator.}
Model trajectories and sensitivities were computed using the training initial conditions \((x^{(1)}_1(0),x^{(1)}_2(0))=(2,0)\).
\textbf{A--B,} Time-dependent sensitivities of the state variables \(x_1\) and \(x_2\) with respect to all model parameters.
\textbf{C--D,} Zoomed views showing near-zero and highly overlapping sensitivity profiles.
These patterns indicate weakly informed and redundant parameter directions, consistent with a non-invertible FIM, a rank-deficient sensitivity matrix and practical non-identifiability.}
    \label{fig:data_CDO_sens_NODEs}
\end{figure}

Figure~\ref{fig:data_CDO_sens_NODEs} shows the time-dependent sensitivities of both state variables with respect to all model parameters. Many parameters exhibit near-overlapping profiles, whereas others remain weakly sensitive throughout the experimental window. These patterns indicate redundant and weakly informed parameter directions, consistent with the non-invertibility of the FIM. The \(\delta_{\mathrm{msqr}}\) values span several orders of magnitude \((2\times10^{-1}, 6\times10^{-4})\), further confirming the highly uneven influence of individual parameters on the model outputs.

The sensitivity matrix has rank=1259 (instead of 1362), confirming that the conventional NODE contains parameter directions that cannot be independently identified from the available data. Thus, despite the large number of observations, the model is practically non-identifiable under this experimental design.

\FloatBarrier

\subsection*{S.3.2 \texttt{Data-CDO: }Implementation with the iNODE framework}
\addcontentsline{toc}{subsection}{S.3.2 \texttt{Data-CDO: }Implementation with the iNODE framework} 
We next formulated the fully data-driven cubic damped oscillator using the iNODE framework. We generated candidate multilayer perceptrons with one and two hidden layers, subject to the minimum data-to-parameter ratio adopted. The admissible architectures were analytically embedded in the ODE system. Candidate models were calibrated, and sensitivity and identifiability analyses were performed using AMIGO2 as described in the main text. Table~\ref{tab:data_CDO_iNODEs_grid} summarises the results.

\begin{table}[h!]
\centering
\caption[\texttt{Data-CDO: }Candidate iNODE architectures for the fully data-driven cubic damped oscillator.]{
\textbf{\texttt{Data-CDO: }Candidate iNODE architectures for the fully data-driven cubic damped oscillator.}
Candidate multilayer perceptron architectures are denoted as MLP(\(n_1,\dots,n_L\)), where each entry indicates the number of neurons in a hidden layer.
For each architecture, we report the best-fit objective value, the total number of estimated parameters (\(p\)), the Akaike information criterion (AIC), the training normalised root-mean-square error (NRMSE, \%), the maximum and median confidence intervals (CI) across all estimated parameters, and the Model Quality Ranking Index (MQRI).
The architecture highlighted in orange has the lowest MQRI and is selected as the best compromise between predictive accuracy, parsimony and practical identifiability.
}
\vspace{6pt}
\label{tab:data_CDO_iNODEs_grid}
\begin{tabular}{cccccccc}
\hline
Architecture & Best fit & p & AIC & $ \text{NRMSE}_{\text{train}} (\%) $ & Max CI (\%) & Median CI (\%) & MQRI \\ \hline
MLP(2)     & 206.77 & 12 & -4466  & 15.80 & $4.33 \times 10^{3}$ & $1.08 \times 10^{3}$  & 59 \\
MLP(2,2)   & 205.34 & 18 & -4464  & 15.70 & $3.41 \times 10^{5}$ & $8.90 \times 10^{3}$  & 66 \\
MLP(4)     & 2.21   & 22 & -7395  & 1.62 & $3.24 \times 10^{2}$ & $8.99$     & 14 \\
MLP(4,2)   & 98.81  & 28 & -4878  & 11.10 & $1.02 \times 10^{4}$ & $6.44 \times 10^{2}$  & 49 \\
MLP(2,4)   & 6.08   & 28 & -9044  & 2.76 & $1.00 \times 10^{3}$ & $99.93$    & 27 \\
\rowcolor[HTML]{FFCE93}
MLP(6)     & 0.84   & 32 & -9817  & 1.04 & $3.81 \times 10^{2}$ & $19.34$    & 12 \\
MLP(6,2)   & 207.98 & 38 & -7205  & 15.87 & $1.60 \times 10^{6}$ & $1.21 \times 10^{4}$& 76 \\
MLP(2,6)   & 6.32   & 38 & -10001 & 2.82 & $1.93 \times 10^{3}$ & $1.13 \times 10^{2}$   & 31 \\
MLP(8)     & 0.39   & 42 & -15432 & 0.68 & $1.27 \times 10^{4}$ & $1.25 \times 10^{2}$   & 20 \\
MLP(4,4)   & 3.42   & 42 & -12182 & 2.05 & $5.38 \times 10^{5}$ & $47.92$    & 32 \\
MLP(8,2)   & 53.56  & 48 & -4554  & 8.51 & $9.78 \times 10^{5}$ & $3.54 \times 10^{3}$  & 60 \\
MLP(2,8)   & 5.84   & 48 & -16447 & 2.63 & $6.87 \times 10^{6}$ & $4.90 \times 10^{2}$   & 47 \\
MLP(10)    & 0.18   & 52 & -18503 & 0.49 & $3.09 \times 10^{3}$ & $1.32 \times 10^{2}$   & 14 \\
MLP(6,4)   & 2.24   & 56 & -12274 & 1.65 & $5.77 \times 10^{4}$ & $3.41\times 10^{2}$   & 34 \\
MLP(4,6)   & 1.82   & 56 & -11960 & 1.52 & $3.45 \times 10^{4}$ & $2.54\times 10^{2}$   & 27 \\
MLP(10,2)  & 137.26 & 58 & -4578  & 12.86 & $6.55 \times 10^{4}$ & $1.01\times 10^{3}$  & 57 \\
MLP(2,10)  & 36.93  & 58 & -18647 & 6.98 & $9.43 \times 10^{6}$ & $1.68\times 10^{3}$  & 61 \\
MLP(12)    & 0.27   & 62 & -12171 & 0.59 & $7.28 \times 10^{5}$ & $3.20\times 10^{2}$   & 27 \\
MLP(12,2)  & 188.47 & 68 & -4534  & 15.14 & $6.56 \times 10^{8}$ & $2.27\times 10^{4}$ & 74 \\
MLP(2,12)  & 8.25   & 68 & -9349  & 3.16 & $1.57 \times 10^{6}$ & $9.79\times 10^{2}$   & 53 \\ \hline
\end{tabular}
\end{table}

Among the evaluated candidates, MLP(6) provides the best MQRI-ranked compromise between fit, complexity and practical identifiability. Although larger single-layer architectures reduce the training error, they also produce substantially wider confidence intervals. Two-layer architectures generally perform worse, showing either poorer fits or markedly larger uncertainty. Thus, additional neural capacity does not translate into more reliable models when the added parameters are weakly constrained by the data.

\begin{figure}[h!]
    \centering
    \includegraphics[width=0.87\textwidth]{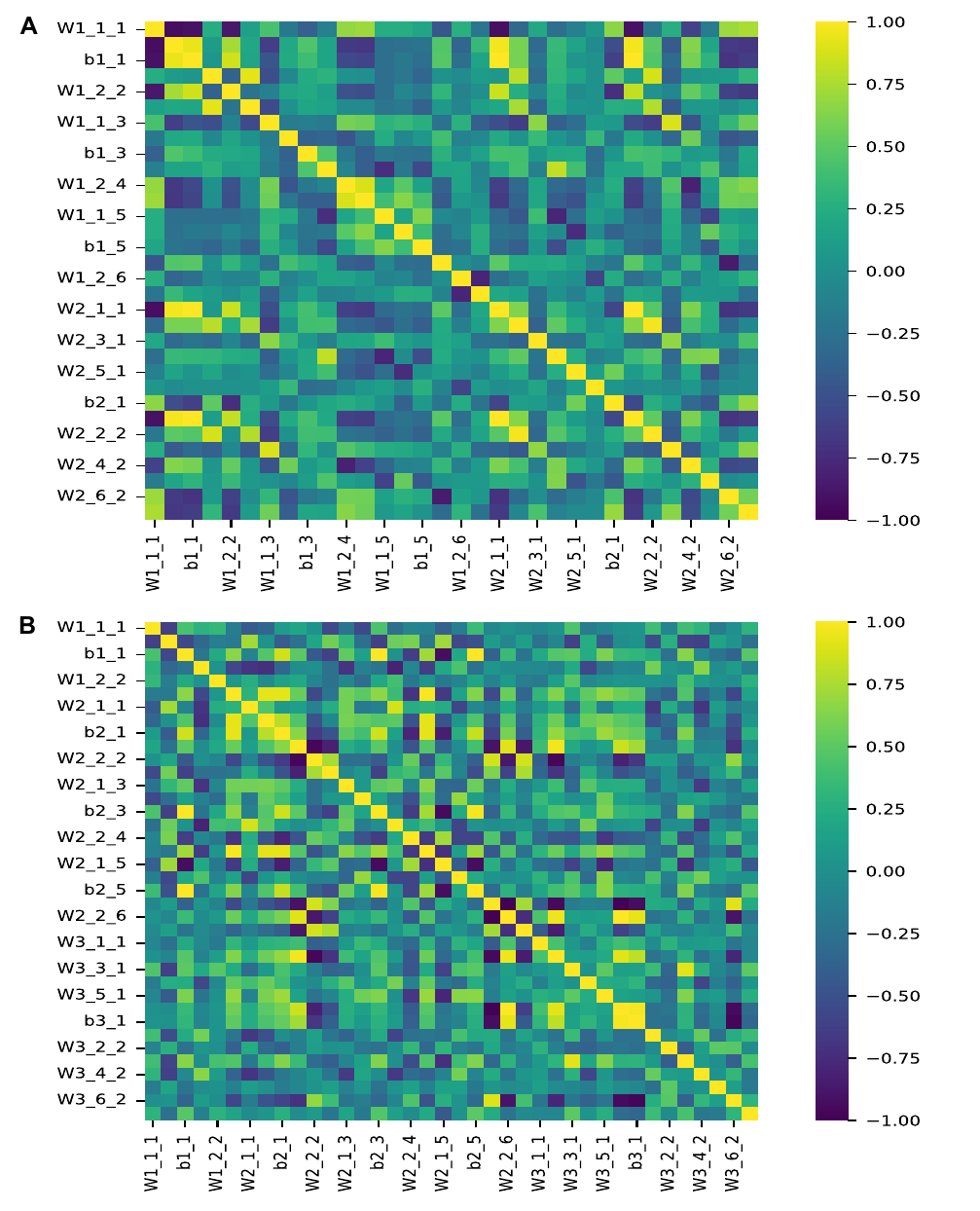}
\caption[\texttt{Data-CDO: }Parameter correlation matrices for two candidate iNODE models.]{
\textbf{\texttt{Data-CDO: }Parameter correlation matrices for two candidate iNODE models of the fully data-driven cubic damped oscillator.}
\textbf{A,} Selected iNODE architecture with one hidden layer of six neurons, MLP(6), corresponding to 32 parameters.
\textbf{B,} Alternative iNODE architecture with two hidden layers, MLP(6,2), corresponding to 38 parameters. Rows and columns correspond to neural-network weights and biases. The colour scale shows pairwise parameter correlation coefficients; values close to \(\pm 1\) indicate strong linear dependencies, whereas values near zero indicate weak correlations.
For readability, only a subset of parameter labels is displayed on each axis.
}
    \label{fig:data_CDO_iNODEs_corr}
\end{figure}

The parameter correlation matrices (Figure~\ref{fig:data_CDO_iNODEs_corr}) provide further information on the markedly different behaviour of single- and two-layer models. For the selected model with moderate complexity MLP(6) (Figure~\ref{fig:data_CDO_iNODEs_corr}, Panel~\textbf{A}), the correlation matrix is largely dominated by diagonal terms, and most of the off-diagonal entries remain close to zero. This indicates weak linear dependencies between parameters, both within the neural network weights and between weights and biases across layers. The generally weaker off-diagonal correlations indicate that the MLP(6) parameterisation is better conditioned than the deeper MLP(6,2) alternative, although several individual parameters remain weakly constrained.

In contrast, the deeper architecture MLP(6,2) (Figure~\ref{fig:data_CDO_iNODEs_corr}, Panel~\textbf{B}) exhibits extended blocks of strong positive and negative correlations, particularly among parameters across different layers. This indicates pronounced compensatory effects whereby changes in one parameter can be offset by others without substantially affecting model output. This correlation is a clear signature of overparameterisation and non-identifiability, directly explaining the excessively large maximum and median confidence intervals observed for this model.

\begin{figure}[h!]
    \centering
    \includegraphics[width=0.90\textwidth]{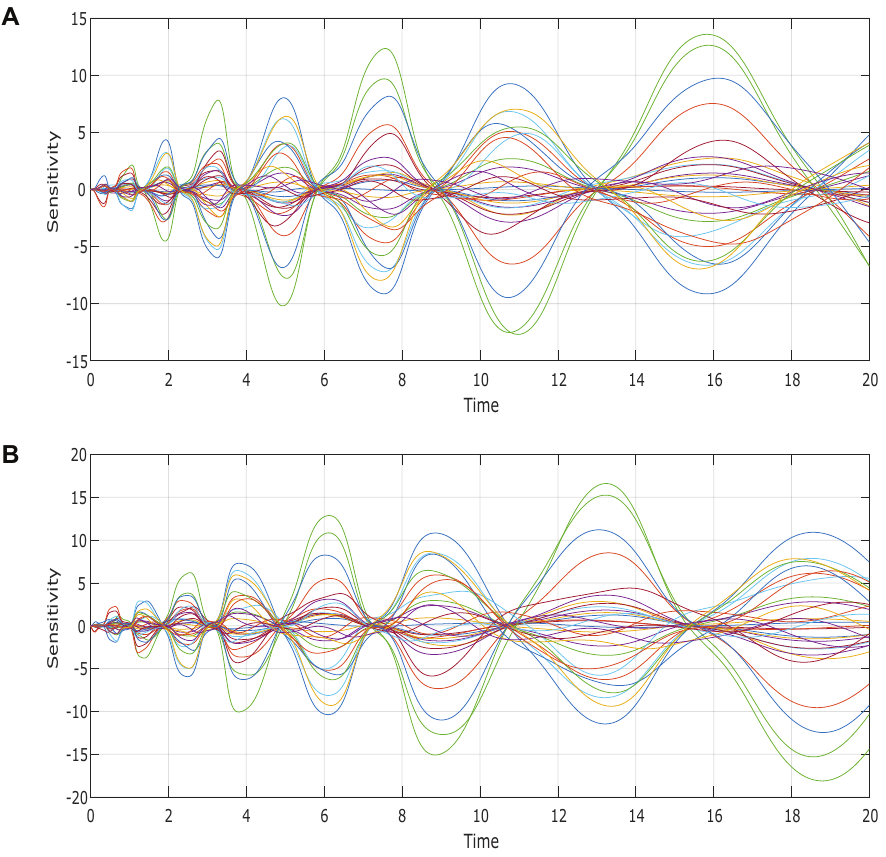}
\caption[\texttt{Data-CDO: }Sensitivity profiles for the selected iNODE model of the cubic damped oscillator.]{
\textbf{\texttt{Data-CDO: }Sensitivity profiles for the selected iNODE model of the cubic damped oscillator.}
Time-dependent sensitivities of the iNODE outputs \(x_1\) (\textbf{A}) and \(x_2\) (\textbf{B}) with respect to the estimated model parameters, computed using the initial conditions of Experiment~1.
}
    \label{fig:data_CDO_sens_iNODEs}
\end{figure}

Sensitivity analysis of the selected MLP(6) model supports this interpretation (Figure~\ref{fig:data_CDO_sens_iNODEs}). Unlike the conventional NODE baseline, the iNODE sensitivities do not collapse towards near-zero values or show widespread overlap. This is consistent with lower parameter correlation, a full-rank sensitivity matrix and an invertible FIM.

This qualitative improvement is further supported by the distribution of $\delta^{msqr}$ values, which ranges from $2 \times 10^{-1}$ to $5 \times 10^{-3}$ with most parameters concentrated in the range $10^{-2}$--$10^{-1}$, reflecting a more balanced sensitivity structure across the parameter space as compared to conventional NODE model. These results indicate a substantially better-conditioned parameterisation than the conventional NODE model, although several parameters remain weakly constrained.

Table~\ref{tab:data_CDO_iNODEs} reports the calibrated parameters and relative confidence intervals for the selected iNODE model. Most weights and biases are well constrained, although a subset has confidence intervals above 100\%. These parameters have non-negligible magnitudes and were therefore not removed, because pruning in this work requires both high uncertainty and limited estimated contribution.
The single-experiment design follows the conventional NODE benchmark configuration used in the reference study.

\begin{table}[H]
\centering
\caption[\texttt{Data-CDO: }Estimated parameters and confidence intervals for the selected iNODE model.]{\textbf{\texttt{Data-CDO: }Estimated parameters and confidence intervals for the selected iNODE model of the cubic damped oscillator.} The table reports the calibrated neural-network weights and biases together with their relative confidence intervals. Parameters with confidence intervals greater than \(100\%\) are highlighted in cyan, indicating weak practical identifiability.
}
\label{tab:data_CDO_iNODEs}
\vspace{6pt}
\begin{tabular}{ccc}
\hline
Parameter & $\theta^*$ & CI (\%) \\ \hline
$W_{111}$ & 1.431  & 98.94  \\
$W_{211}$ & -0.497 & 96.94 \\
\rowcolor[HTML]{96FFFB}
$b_{11}$  & 1.033  & 162.94 \\
$W_{112}$ & -3.139 & 12.46 \\
$W_{212}$ & 0.468  & 26.56  \\
$b_{12}$  & -2.725 & 5.70  \\
$W_{113}$ & 3.231  & 5.58   \\
$W_{213}$ & 0.098  & 53.05  \\
$b_{13}$  & -3.101 & 2.40  \\
$W_{114}$ & -1.949 & 17.92 \\
$W_{214}$ & 2.888  & 11.69  \\
$b_{14}$  & 2.617  & 6.95   \\
$W_{115}$ & -0.695 & 12.45 \\
$W_{215}$ & -2.776 & 3.39  \\
$b_{15}$  & -3.054 & 1.87  \\
$W_{116}$ & -0.423 & 17.05 \\
$W_{216}$ & -3.456 & 2.40  \\
$b_{16}$  & 2.832  & 1.48   \\
\rowcolor[HTML]{96FFFB}
$W_{211}$ & 0.206  & 381.35 \\
$W_{221}$ & 1.644  & 35.47  \\
$W_{231}$ & -0.870 & 22.13 \\
$W_{241}$ & 0.624  & 22.46  \\
$W_{251}$ & -4.742 & 9.45  \\
$W_{261}$ & -2.176 & 3.79  \\
$b_{21}$  & -2.578 & 29.78 \\
\rowcolor[HTML]{96FFFB}
$W_{212}$ & 0.720  & 381.51 \\
$W_{222}$ & 7.910  & 31.25  \\
$W_{232}$ & -10.515 & 20.99 \\
$W_{242}$ & 0.663  & 20.76  \\
$W_{252}$ & 1.624  & 15.95  \\
$W_{262}$ & 0.141  & 46.95  \\
\rowcolor[HTML]{96FFFB}
$b_{22}$  & -2.353 & 139.86 \\ \hline
\end{tabular}
\end{table}

\section*{S.4~\texttt{Hybrid-CDO: }Hybrid mechanistic-neural cubic damped oscillator} 
\addcontentsline{toc}{section}{S.4~\texttt{Hybrid-CDO: }Hybrid mechanistic-neural cubic damped oscillator} 
In the second case study based on the cubic damped oscillator, we investigate a hybrid modelling approach using Hybrid Neural Ordinary Differential Equations (HNODE). As discussed above, this class of models combines explicit mechanistic differential equations with neural network components, allowing prior physical knowledge to be incorporated while retaining the flexibility of data-driven representations. The hybrid formulation considered here is given by

\begin{equation}
\left\{
\begin{aligned}
\frac{d x_1}{d t} &= -0.1\, x_1^{3} + 2.0\, x_2^{3}\\
\frac{d x_2}{d t} &= \mathcal{NN}\big(x_1(t), x_2(t); \boldsymbol{\phi}_{NN}\big)
\end{aligned}
\right.
\end{equation}

where the first equation represents the mechanistically defined component of the system, and $\mathcal{NN}(\cdot;\boldsymbol{\phi}_{NN})$ denotes a neural network parameterised by weights and biases $\boldsymbol{\phi}_{NN}$. This neural component captures the unknown or unmodelled dynamics governing the evolution of the second state variable. This formulation allows assessment of how the incorporation of mechanistic structure affects predictive performance and identifiability, and how the informativeness of the available data further shapes parameter recovery within an HNODE framework.

We considered two experimental scenarios with the same total number of measurements ($1000$) but different levels of dynamic diversity. Full details are reported in Table~\ref{tab:hybrid_CDO_scenarios}. We used both the  Hybrid iNODE and the conventional HNODE frameworks to model the system.

\begin{table}[h!]
\centering
\caption[\texttt{Hybrid-CDO: }Training and validation scenarios for one- and two-experiment hybrid cubic damped oscillator designs.]{
\textbf{\texttt{Hybrid-CDO: }Training and validation scenarios for one- and two-experiment hybrid cubic damped oscillator designs.}
The scenarios define the training and validation trajectories used to assess how experimental diversity affects parameter identifiability, phase-space recovery and extrapolation in the hybrid iNODE formulation.
}
\vspace{6pt}
\label{tab:hybrid_CDO_scenarios}
\begin{tabular}{cccccc}
\hline
\textbf{Dataset} &
  \textbf{Scenario} &
  \textbf{Experiment(s)} &
  \textbf{Time} &
  \textbf{\begin{tabular}[c]{@{}c@{}}Measurements\\  per state\end{tabular}} &
  \textbf{\begin{tabular}[c]{@{}c@{}}Total data \\ points\end{tabular}} \\ \hline
\multirow{3}{*}{Train} & 1exp & \((x^{(1)}_1(0),x^{(1)}_2(0))=(2,0)\)   & {[}0, 20{]}  & 500 & 1000 \\ [3mm]
 &
  \multirow{2}{*}{2exp} &
  \((x^{(1)}_1(0),x^{(1)}_2(0))=(2,0)\) &
  {[}0, 25{]} &
  \multirow{2}{*}{\begin{tabular}[c]{@{}c@{}}250 in each \\ experiment\end{tabular}} &
  \multirow{2}{*}{1000} \\ 
                       &      & \((x^{(2)}_1(0),x^{(2)}_2(0))=(0.5,2)\) & {[}0, 25{]}  &     &      \\ [3mm]
\multirow{2}{*}{Validation} & val1 & \((x^{(1)}_1(0),x^{(1)}_2(0))=(2,0)\)   & {[}20, 40{]} & 100 & 200  \\[3mm]
                            & val2 & \((x^{(1)}_1(0),x^{(1)}_2(0))=(1,1)\)   & {[}0, 25{]}  & 500 & 1000 \\ \hline
\end{tabular}
\end{table}

\subsection*{S.4.1 \texttt{Hybrid-CDO: }Implementation with the hybrid iNODE framework}
\addcontentsline{toc}{subsection}{S.4.1 \texttt{Hybrid-CDO: }Implementation with the hybrid iNODE framework} 
We first used the hybrid iNODE framework under the \texttt{1exp} scenario. We performed a systematic grid search over candidate multilayer perceptron architectures. Table~\ref{tab:hybird_CDO_iNODEs_gridsearch} reports the corresponding calibration, identifiability, and model-selection metrics, including the MQRI used to identify the best overall architecture. After selecting this architecture, we examined how predictive performance and practical identifiability changed in the different experimental schemes. These results are summarised in Table~\ref{tab:hybrid_CDO_iNODEs_models}.

The architecture search again shows that the lowest training error does not necessarily correspond to the most identifiable model. Although MLP(10) gives the smallest training error, MLP(8) achieves the lowest MQRI among the evaluated candidates by balancing predictive accuracy with lower parameter uncertainty. Deeper two-layer architectures generally show larger confidence intervals and higher MQRI values, indicating weaker practical identifiability.

\begin{table}[h!]
\centering
\caption[\texttt{Hybrid-CDO: }Candidate hybrid iNODE architectures for the cubic damped oscillator.]{
\textbf{\texttt{Hybrid-CDO: }Candidate hybrid iNODE architectures for the cubic damped oscillator.}
Candidate multilayer perceptron architectures are denoted as MLP(\(n_1,\dots,n_L\)), where each entry indicates the number of neurons in a hidden layer.
For each architecture, we report the best-fit objective value, the total number of estimated parameters (\(p\)), the Akaike information criterion (AIC), the training normalised root-mean-square error (NRMSE), the maximum and median confidence intervals across all estimated parameters, and the Model Quality Ranking Index (MQRI).
All neural networks use hyperbolic tangent activation functions in the hidden layers and a linear activation function in the output layer.
The architecture highlighted in orange has the lowest MQRI and is selected as the best compromise between predictive accuracy, parsimony and practical identifiability.
}
\vspace{6pt}
\label{tab:hybird_CDO_iNODEs_gridsearch}
\begin{tabular}{cccccccc}
\hline
Architecture & Best fit & p & AIC & $ \text{NRMSE}_{\text{train}} (\%) $ & Max CI (\%) & Median CI (\%) & MQRI \\ \hline
MLP(2)   & 0.884 & 9  & -5427  & $1.69$ & $4.73\times 10^{2}$ & $1.22 \times 10^{2}$ & 26 \\
MLP(2,2) & 0.102 & 15 & -7144  & $0.56$ & $2.10 \times 10^{5}$  & $2.55 \times 10^{3}$  & 33 \\
MLP(4)   & 0.089 & 17 & -7249  & $0.52$ & $1.70 \times 10^{5}$ & $3.18 \times 10^{2}$   & 27 \\
MLP(2,4) & 0.156 & 23 & -6787  & $0.69$ & $6.69  \times 10^{6}$ & $1.34\times 10^{4}$ & 38 \\
MLP(6)   & 0.018 & 25 & -8495  & $0.24$ & $3.1\times 10^{3}$ & 58 & 11 \\
MLP(4,2) & 0.083 & 25 & -7285  & $0.50$ & - & - & -\\
MLP(2,6) & 0.080 & 31 & -7309  & $0.48$ & $5.11  \times 10^{5}$ & $2.95\times 10^{4}$ & 28 \\
\rowcolor[HTML]{FFCE93} 
MLP(8)   & 0.008 & 33 & -9099  & $0.16$ & $8.81 \times 10^{2}$ & 53 & 7  \\
MLP(6,2) & 3.276 & 35 & -4326  & $3.15$ & $3.29 \times 10^{7}$  & $4.51\times 10^{2}$   & 40 \\
MLP(4,4) & 0.046 & 37 & -7734  & $0.38$ & $1.15 \times 10^{5}$  & $1.43\times 10^{3}$  & 20 \\
MLP(2,8) & 0.082 & 39 & -7268  & $0.49$ & - & - & - \\
MLP(10)  & 0.002 & 41 & -10188 & $0.08$ & $7.05 \times 10^{3}$ & $2.32\times 10^{2}$ & 10 \\ \hline
\end{tabular}
\end{table}

As in the \texttt{Data-CDO} case study, the grid search shows that shallow architectures systematically outperform deeper two-layer configurations in both fit quality and practical identifiability. For comparable numbers of trainable parameters, single-hidden-layer models yield lower errors and markedly tighter confidence intervals than two-layer alternatives. In addition, deeper architectures frequently exhibit greater uncertainty and higher MQRI values. 

\pagebreak

The selected MLP(8) model achieved a best-fit value of \(8.45\times10^{-3}\), but several low-magnitude parameters had confidence intervals above 100\%. Removing these weakly informed parameters and recalibrating the model reduced the maximum CI from \(>800\%\) to \(368\%\), without degrading the fit (Table~ \ref{tab:hybrid_CDO_iNODEs_models}).

We then recalibrated the same architecture using two experiments with distinct initial conditions. This multi-experiment design further improved practical identifiability, and pruning the remaining weakly informed parameters reduced the maximum CI to \(84.2\%\). The mean absolute CI decreased from \(97.76\%\) to \(57.07\%\) in the single-experiment case and from \(55.84\%\) to \(22.08\%\) in the two-experiment case after pruning. These results show that pruning improves conditioning, but experimental diversity provides the larger gain in identifiability.

\begin{figure}[H]
    \centering
    \includegraphics[width=0.95\textwidth]{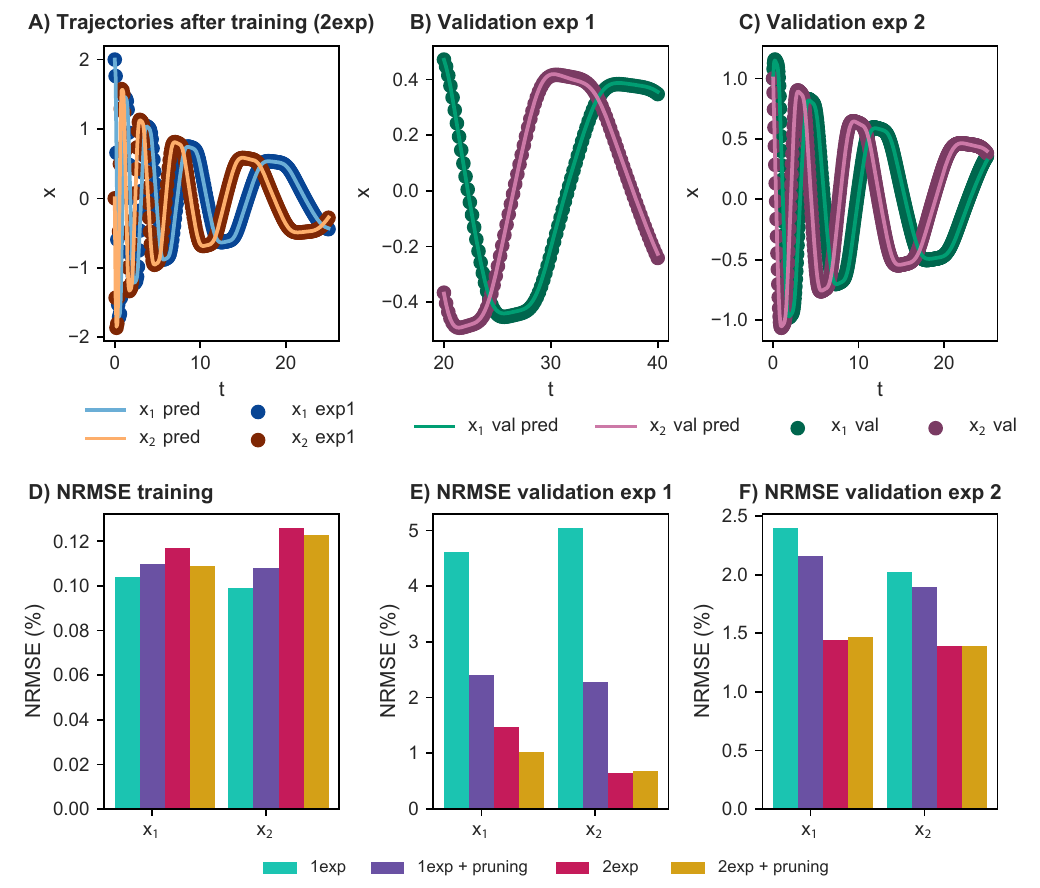}
\caption[\texttt{Hybrid-CDO: }Predictions and NRMSE comparison for the hybrid iNODE cubic damped oscillator.]{
\textbf{\texttt{Hybrid-CDO: }Predictions and NRMSE comparison for the hybrid iNODE cubic damped oscillator.}
\textbf{A--C,} Predicted trajectories of the \texttt{CDO\_2exp\_iprunn} model for the training experiment and two validation experiments, including one extending beyond the training time range. Points denote synthetic observations, solid lines denote model predictions, and shaded regions indicate uncertainty bands.
\textbf{D--F,} NRMSE comparison across the four model configurations for the training and two validation scenarios.
The results show that multi-experiment calibration combined with identifiability-aware pruning improves recovery of the cubic damped oscillator dynamics and provides the best balance between predictive accuracy and practical identifiability.
}
    \label{fig:Hybrid_CDO_iNODEs_models}
\end{figure}

The single-experiment hybrid iNODE achieves low training error but only moderate validation performance, with errors around 5\% in Figure~\ref{fig:Hybrid_CDO_iNODEs_models} Panels~\textbf{D-F}. Generalisation improves in the two-experiment setting, with \(NRMSE_{\mathrm{val1}}=1.06\%\) and \(NRMSE_{\mathrm{val2}}=1.42\%\). The pruned two-experiment model provides the best overall balance, yielding \(NRMSE_{\mathrm{train}}=0.116\%\), \(NRMSE_{\mathrm{val1}}=0.846\%\) and \(NRMSE_{\mathrm{val2}}=1.428\%\), together with improved conditioning.

\begin{table}[H]
\centering
\small
\caption[\texttt{Hybrid-CDO: }Estimated parameters and confidence intervals for the hybrid iNODE configuration of the cubic damped oscillator.]{
\textbf{\texttt{Hybrid-CDO: }Estimated parameters and confidence intervals for the hybrid iNODE configuration of the cubic damped oscillator.}
The table reports the calibrated parameter values \((\theta^*)\) and relative confidence intervals (CI, \%) for four configurations: \texttt{CDO\_1exp} and \texttt{CDO\_2exp}, calibrated using one and two experimental datasets, respectively, and their identifiability-guided pruned counterparts, \texttt{CDO\_1exp\_iprun} and \texttt{CDO\_2exp\_iprun}.
Parameters highlighted in blue have large relative confidence intervals, indicating weak practical identifiability.
Entries marked with ``--'' denote parameters removed during identifiability-guided pruning.
}
\vspace{6pt}
\label{tab:hybrid_CDO_iNODEs_models}
\begin{tabular}{c|cc|cc|cc|cc|c}
\hline
 &
  \multicolumn{2}{c|}{1exp} &
  \multicolumn{2}{c|}{1exp + pruning} &
  \multicolumn{2}{c|}{2exp} &
  \multicolumn{2}{c|}{2exp + pruning} &
   \\[3mm]
Best fit &
  \multicolumn{2}{c|}{$8.45\times10^{-3}$} &
  \multicolumn{2}{c|}{$7.54\times10^{-3}$} &
  \multicolumn{2}{c|}{$1.42\times10^{-2}$} &
  \multicolumn{2}{c|}{$1.34\times10^{-2}$} &
   \\[3mm]
Param. &
  $\theta^*$ &
  CI (\%) &
  $\theta^*$ &
  CI (\%) &
  $\theta^*$ &
  CI (\%) &
  $\theta^*$ &
  CI (\%) &
  Final Param. \\ \hline
$W_{111}$ &
  -2.64 &
  27.03 &
  -2.43 &
  18.16 &
  -2.62 &
  17.88 &
  -2.30 &
  10.51 &
  $W_{111}$ \\
$W_{211}$ &
  \cellcolor[HTML]{96FFFB}-0.12 &
  \cellcolor[HTML]{96FFFB}$1.53\times10^{2}$ &
  \cellcolor[HTML]{96FFFB}-0.08 &
  \cellcolor[HTML]{96FFFB}$1.22\times10^{2}$ &
  \cellcolor[HTML]{96FFFB}-0.07 &
  \cellcolor[HTML]{96FFFB}$2.19\times10^{2}$ &
  - &
  - \\
$b_{11}$ &
  0.20 &
  45.58 &
  0.19 &
  32.04 &
  0.23 &
  27.90 &
  0.23 &
  16.83 &
  $b_{11}$\\
$W_{112}$ &
  2.98 &
  40.47 &
  2.79 &
  28.65 &
  3.05 &
  12.28 &
  2.75 &
  11.00 &
  $W_{112}$\\
$W_{212}$ &
  \cellcolor[HTML]{96FFFB}-0.01 &
  \cellcolor[HTML]{96FFFB}$8.81\times10^{2}$ &
  - &
  - &
  \cellcolor[HTML]{96FFFB}-0.06 &
  \cellcolor[HTML]{96FFFB}$1.03\times10^{2}$ &
  - &
  - \\
$b_{12}$ &
  2.70 &
  24.79 &
  2.56 &
  13.26 &
  2.92 &
  7.77 &
  2.71 &
  5.38 &
  $b_{12}$\\
$W_{113}$ &
  -3.98 &
  36.11 &
  -4.16 &
  41.89 &
  -3.66 &
  40.37 &
  -3.52 &
  21.67 &
  $W_{113}$\\
$W_{213}$ &
  \cellcolor[HTML]{96FFFB}0.14 &
  \cellcolor[HTML]{96FFFB}$2.76\times10^{2}$ &
  \cellcolor[HTML]{96FFFB}0.15 &
  \cellcolor[HTML]{96FFFB}$3.68\times10^{2}$ &
  \cellcolor[HTML]{96FFFB}-0.01 &
  \cellcolor[HTML]{96FFFB}$2.78\times10^{2}$ &
  - &
  - \\
$b_{13}$ &
  7.55 &
  26.30 &
  7.96 &
  27.16 &
  7.09 &
  35.15 &
  6.98 &
  17.87 &
  $b_{13}$\\
$W_{114}$ &
  -1.32 &
  39.62 &
  -1.09 &
  24.44 &
  -1.35 &
  32.90 &
  -1.26 &
  11.40 &
  $W_{114}$\\
$W_{214}$ &
  -0.54 &
  60.55 &
  -0.49 &
  30.30 &
  -0.54 &
  42.91 &
  -0.53 &
  19.22 &
  $W_{214}$\\
$b_{14}$ &
  \cellcolor[HTML]{96FFFB}0.07 &
  \cellcolor[HTML]{96FFFB}$1.57\times10^{2}$ &
  - &
  - &
  \cellcolor[HTML]{96FFFB}0.07 &
  \cellcolor[HTML]{96FFFB}$1.45\times10^{2}$ &
  - &
  - \\
$W_{115}$ &
  0.44 &
  53.20 &
  0.41 &
  50.30 &
  1.05 &
  29.24 &
  0.97 &
  17.11 &
  $W_{115}$\\
$W_{215}$ &
  0.96 &
  65.54 &
  0.87 &
  27.68 &
  0.80 &
  31.25 &
  0.82 &
  16.67 &
  $W_{215}$\\
$b_{15}$ &
  1.74 &
  44.56 &
  1.76 &
  15.62 &
  1.75 &
  20.98 &
  1.87 &
  8.94 &
  $b_{15}$\\
$W_{116}$ &
  \cellcolor[HTML]{96FFFB}-0.10 &
  \cellcolor[HTML]{96FFFB}$1.66\times10^{2}$ &
  - &
  - &
  \cellcolor[HTML]{96FFFB}0.08 &
  \cellcolor[HTML]{96FFFB}$1.98\times10^{2}$ &
  - &
  - \\
$W_{216}$ &
  1.40 &
  18.58 &
  1.30 &
  14.36 &
  1.28 &
  12.35 &
  1.25 &
  10.96 &
  $W_{216}$\\
$b_{16}$ &
  0.55 &
  14.03 &
  0.43 &
  15.02 &
  0.29 &
  23.75 &
  0.26 &
  21.72 &
  $b_{16}$\\
$W_{117}$ &
  1.74 &
  10.50 &
  1.64 &
  11.70 &
  1.67 &
  6.83 &
  1.59 &
  4.16 &
  $W_{117}$\\
$W_{217}$ &
  0.08 &
  21.32 &
  0.07 &
  18.44 &
  0.10 &
  14.88 &
  0.09 &
  6.22 &
  $W_{217}$\\
$b_{17}$ &
  -2.42 &
  3.95 &
  -2.39 &
  3.78 &
  -2.39 &
  3.27 &
  -2.30 &
  3.65 &
  $b_{17}$\\
$W_{118}$ &
  2.32 &
  49.95 &
  2.50 &
  42.99 &
  2.33 &
  23.95 &
  2.38 &
  30.36 &
  $W_{118}$\\
$W_{218}$ &
  0.05 &
  89.45 &
  0.05 &
  58.03 &
  0.06 &
  39.79 &
  0.05 &
  37.18 &
  $W_{218}$\\
$b_{18}$ &
  3.61 &
  44.08 &
  3.85 &
  38.82 &
  3.86 &
  20.79 &
  3.98 &
  26.31 &
  $b_{18}$ \\
$W_{211}$ &
  -0.40 &
  88.00 &
  -0.43 &
  41.50 &
  -0.43 &
  61.54 &
  -0.47 &
  18.86 &
  $W_{211}$\\
$W_{221}$ &
  \cellcolor[HTML]{96FFFB}-0.70 &
  \cellcolor[HTML]{96FFFB}$1.70\times10^{2}$ &
  \cellcolor[HTML]{96FFFB}-0.85 &
  \cellcolor[HTML]{96FFFB}$1.27\times10^{2}$ &
  -1.01 &
  47.57 &
  -1.26 &
  51.41 &
  $W_{221}$\\
$W_{231}$ &
  4.56 &
  69.14 &
  \cellcolor[HTML]{96FFFB}4.19 &
  \cellcolor[HTML]{96FFFB}$1.08\times10^{2}$ &
  5.40 &
  63.74 &
  6.30 &
  34.98 &
  $W_{231}$\\
$W_{241}$ &
  0.46 &
  89.95 &
  0.53 &
  37.09 &
  0.47 &
  59.37 &
  0.43 &
  30.41 &
  $W_{241}$\\
$W_{251}$ &
  -0.61 &
  91.52 &
  -0.67 &
  50.46 &
  -0.58 &
  47.23 &
  -0.55 &
  37.15 &
  $W_{251}$\\
$W_{261}$ &
  0.25 &
  56.29 &
  0.27 &
  29.48 &
  0.25 &
  26.36 &
  0.24 &
  21.55 &
  $W_{261}$\\
$W_{271}$ &
  -5.08 &
  29.18 &
  -5.79 &
  33.57 &
  -5.55 &
  22.28 &
  -5.84 &
  5.65 &
  $W_{271}$ \\
$W_{281}$ &
  -6.18 &
  56.22 &
  -5.70 &
  48.94 &
  -6.49 &
  28.18 &
  -6.44 &
  36.92 &
  $W_{281}$\\
$b_{21}$ &
  \cellcolor[HTML]{96FFFB}-2.17 &
  \cellcolor[HTML]{96FFFB}$2.27\times10^{2}$ &
  \cellcolor[HTML]{96FFFB}-2.73 &
  \cellcolor[HTML]{96FFFB}$2.34\times10^{2}$ &
  -2.83 &
  99.72 &
  -3.77 &
  84.23 &
  $b_{21}$\\ \hline
\end{tabular}
\end{table}

\begin{figure}[H]
    \centering
    \includegraphics[width=0.80\textwidth]{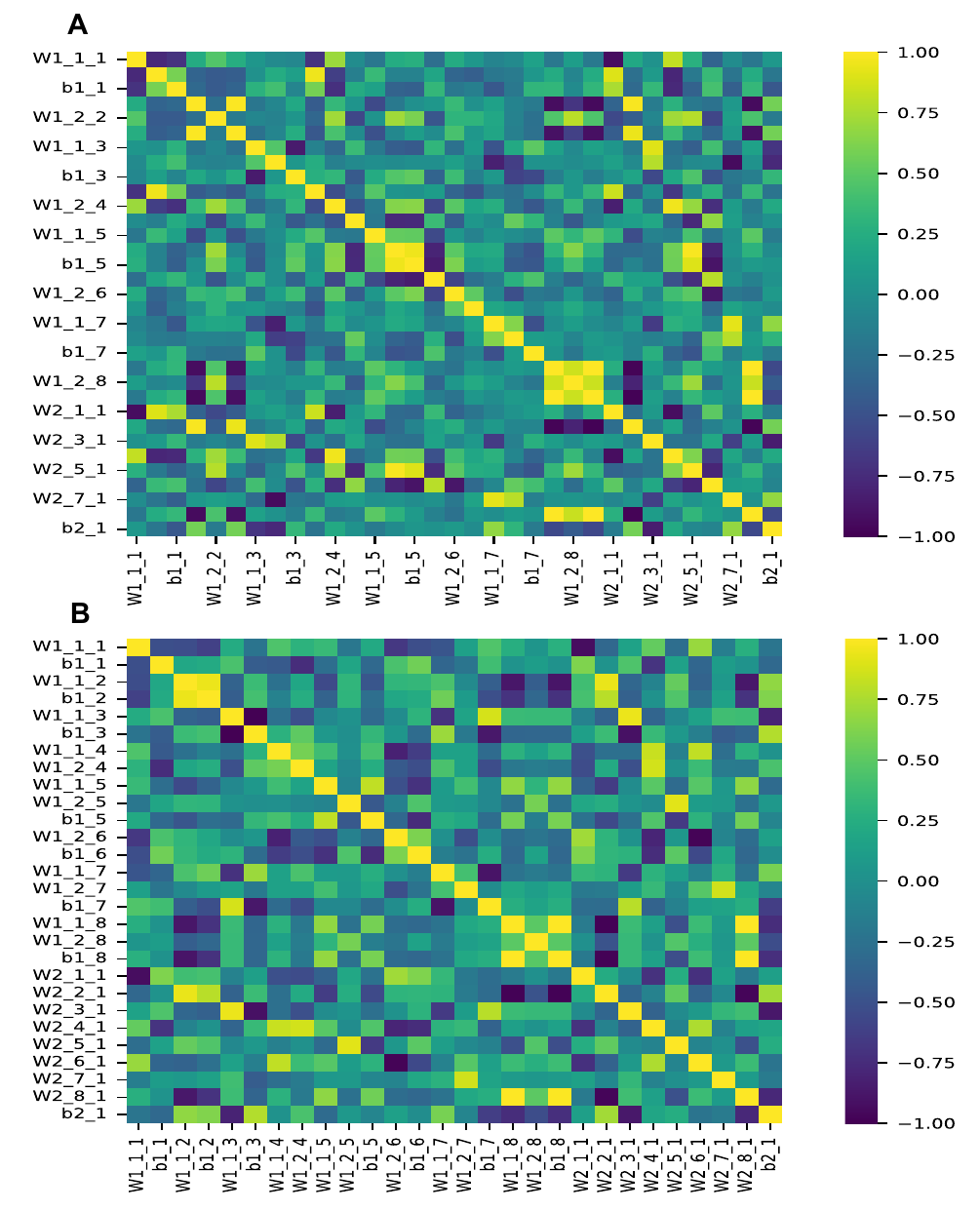}
\caption[\texttt{Hybrid-CDO: }Parameter correlation matrices for the selected hybrid iNODE models.]{
\textbf{\texttt{Hybrid-CDO: }Parameter correlation matrices for the selected hybrid iNODE models.}
\textbf{A,} Correlation matrix for the selected architecture obtained under the initial experimental design, comprising 33 parameters.
\textbf{B,} Correlation matrix for the selected architecture after calibration with two experiments and subsequent pruning, comprising 28 parameters.
Rows and columns correspond to the weights and biases of the neural network embedded in the hybrid NODE model.
The colour scale shows pairwise Pearson correlation coefficients; values close to \(\pm 1\) indicate strong linear dependencies, whereas values near zero indicate weak correlations.
For readability, only a subset of parameter labels is displayed on each axis.
Compared with the initial design, the two-experiment pruned model exhibits a reduced prevalence of strong parameter correlations, consistent with improved practical identifiability.
}
\label{fig:Hybrid_CDO_iNODEs_corr}
\end{figure}

The correlation matrices in Figure~\ref{fig:Hybrid_CDO_iNODEs_corr} provide a structural explanation for the improvement in confidence intervals. The single-experiment model shows blocks of strong positive and negative correlations, indicating compensatory relationships among neural weights and biases. After adding a second experiment and pruning weakly informed parameters, these correlations are substantially reduced and the matrix becomes closer to diagonal dominance. This decorrelation confirms that the retained parameters are more independently constrained by the data, without loss of predictive accuracy.

This parameter decorrelation provides a structural explanation for the observed reduction in confidence intervals and for the elimination of parameters with relative uncertainty exceeding 100\%. Importantly, these gains are achieved without compromising predictive performance, confirming that the removed parameters were not essential to accurately represent the system dynamics. Overall, the correlation analysis reinforces the conclusion that combining experimental diversity with identifiability-guided network pruning yields a well-conditioned hybrid iNODE model with statistically robust and dynamically meaningful parameter estimates.

\subsection*{S.4.2 \texttt{Hybrid-CDO: }Implementation with the conventional HNODE framework}
\addcontentsline{toc}{subsection}{S.4.2 \texttt{Hybrid-CDO: }Implementation with the conventional HNODE framework} 
To benchmark the performance of the proposed identifiability‑aware hybrid modelling pipeline, a conventional Hybrid Neural Ordinary Differential Equations (HNODE) model was implemented in Python and trained following standard practices commonly adopted in the literature. This reference model was calibrated using exactly the same experimental data as the \texttt{CDO\_2exp} configuration, ensuring a fair and controlled comparison between the two approaches.

The conventional HNODE baseline was trained in Python using the same \texttt{CDO\_2exp} dataset. Training was formulated as a continuous-time trajectory-level optimisation problem, with model trajectories obtained by numerical integration from the initial conditions of both experiments. Parameters were estimated by minimising a time-weighted mean squared error using automatic differentiation, Adam optimisation, gradient clipping and scheduled learning-rate reduction. Trajectory subsampling and multiple shooting were used to improve training stability. Candidate architectures were evaluated by the average training NRMSE across both experiments, and the best-performing architecture was then trained to convergence.

Neural network architectures of varying depth and width are systematically explored through a grid‑search procedure, where each candidate model is trained for a limited number of epochs to assess its ability to capture the underlying dynamics (Table~\ref{tab:Hybrid_CDO_HNODEs}). The model is selected by evaluating the normalised root‑mean‑square error (NRMSE) averaged across both experiments, and the architecture yielding the lowest NRMSE is used for full convergence. The final optimised HNODE is then used to generate continuous‑time predictions, which are compared directly against the experimental data. This training strategy enables reliable identification of the hybrid dynamical system while maintaining stability of the learnt vector field and robust generalisation under experimental conditions.

\begin{table}[h!]
\centering
\caption[\texttt{Hybrid-CDO: }Comparison of conventional HNODE architectures for the cubic damped oscillator.]{
\textbf{\texttt{Hybrid-CDO: }Comparison of conventional HNODE architectures for the cubic damped oscillator.}
The table reports the total number of trainable parameters and the training normalised root-mean-square error (NRMSE) for each evaluated architecture.
All neural networks use hyperbolic tangent activation functions in the hidden layers and a linear activation function in the output layer.
The highlighted architecture corresponds to the selected conventional HNODE baseline used for comparison with the hybrid iNODE models.
}
\label{tab:Hybrid_CDO_HNODEs}
\vspace{6pt}
\begin{tabular}{ccc}
\hline
Architecture      & Parameters & $ \text{NRMSE}_{\text{train}} (\%) $ \\ \hline
MLP(8)           & 33         & 18.08 \\
MLP(16)          & 65         & 18.29 \\
MLP(8,8)         & 105        & 15.82 \\
MLP(8,8,8)       & 177        & 15.89  \\
MLP(8,8,8,8)     & 249        & 14.24 \\
MLP(16,16)       & 337        & 15.95 \\
MLP(16,16,16)    & 609        &  14.15\\
\rowcolor[HTML]{FFCE93} 
MLP(16,16,16,16) & 881        & 5.91  \\ \hline
\end{tabular}
\end{table}

Table~\ref{tab:Hybrid_CDO_HNODEs} shows that shallow conventional HNODE architectures underfit the hybrid CDO dynamics, whereas deeper and wider models reduce the training error. The selected MLP(16,16,16,16) baseline contains 881 parameters and achieves the lowest training NRMSE among the evaluated architectures, indicating that the conventional workflow requires substantially higher capacity to approximate the same dynamics.

The difference between the conventional HNODE and hybrid iNODE results cannot be explained by architecture size alone. The same nominal MLP(8) architecture performs substantially worse in the conventional workflow than in the analytically embedded iNODE formulation. This suggests that model formulation and calibration strategy, not only network capacity, affect the ability to recover predictive and identifiable parameterisations in this nonlinear system. The iNODE formulation benefits from analytic embedding and global optimisation, which appear to facilitate the identification of parameter sets that are both more predictive and more consistent with the underlying dynamics. 

\begin{figure}[h!]
    \centering
    \includegraphics[width=0.92\textwidth]{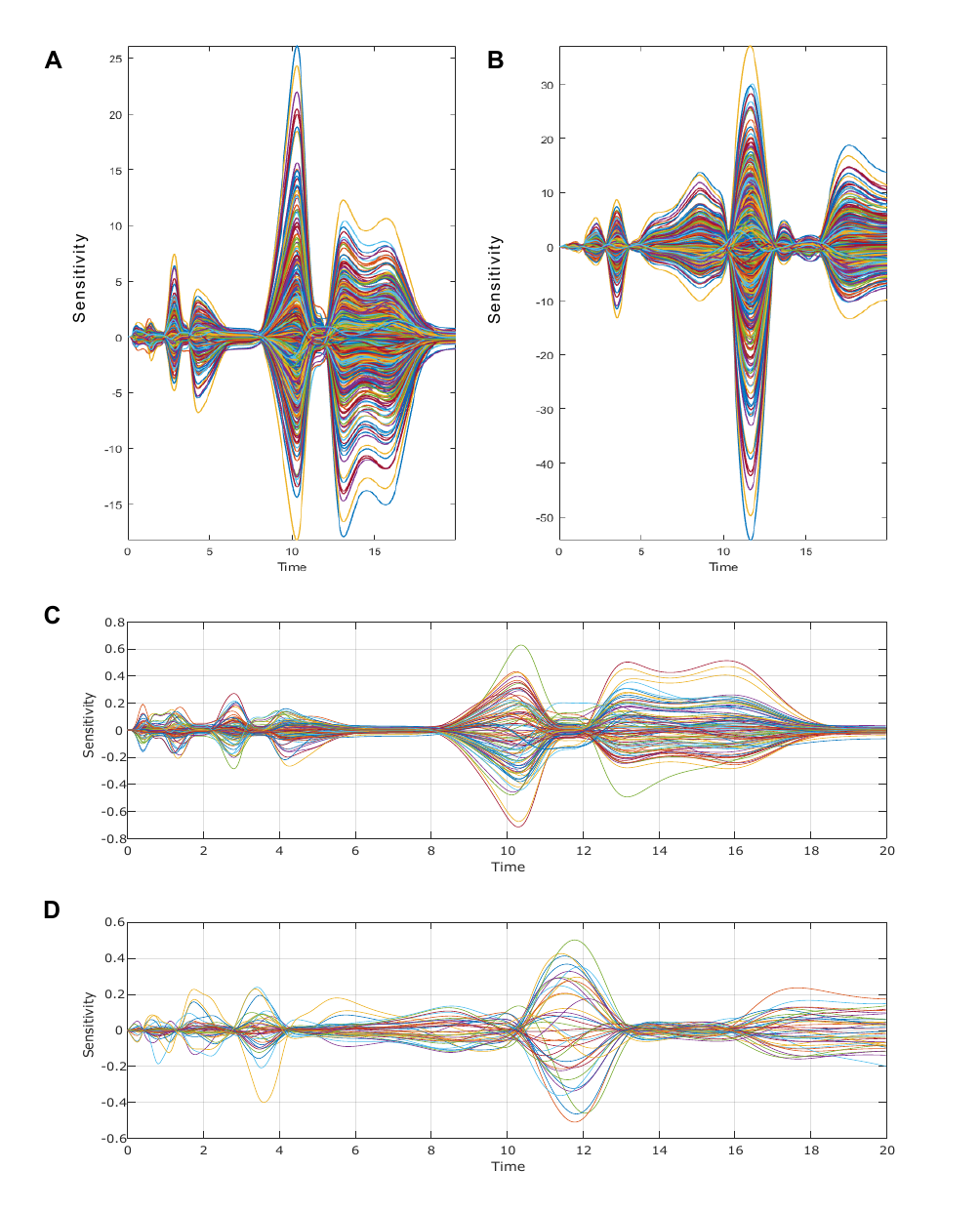}
\caption[\texttt{Hybrid-CDO: }Sensitivity profiles for the conventional HNODE model of the cubic damped oscillator.]{
\textbf{\texttt{Hybrid-CDO: }Sensitivity profiles for the conventional HNODE model of the cubic damped oscillator.}
\textbf{A--B,} Time-dependent sensitivities of the HNODE outputs \(x_1\) and \(x_2\) with respect to the estimated model parameters, computed using the initial conditions of Experiment~1.
\textbf{C--D,} Zoomed views of the sensitivities shown in \textbf{A} and \textbf{B}.
The prevalence of near-zero and overlapping sensitivity profiles is consistent with a non-invertible FIM, a rank-deficient sensitivity matrix and practical non-identifiability of the conventional HNODE model.
}   \label{fig:Hybrid_CDO_sens_HNODEs}
\end{figure}

After training, we used AMIGO2 to compute confidence intervals for the parameters in the optimal HNODE architectures, MLP(16,16,16,16). However, the resulting FIM was found to be non-invertible, which precluded the computation of confidence intervals and suggested potential identifiability issues or parameter redundancy in the model. The parameter sensitivity analysis of the HNODE model reveals pronounced heterogeneity between parameters, with only a small subset exerting a meaningful influence on the observable states (Figure~\ref{fig:Hybrid_CDO_sens_HNODEs}). Most parameters exhibit near-zero or highly correlated sensitivities, indicating that multiple parameter directions contribute redundantly to the model output. Less influential parameters include $W_{3,10,12}$, $W_{3,11,13}$, $W_{2,2,12}$, $W_{4,4,7}$, $W_{4,4,10}$, $W_{3,5,12}$, $W_{4,9,17}$, $W_{2,7,4}$, $W_{4,16,6}$, $W_{2,7,11}$, $W_{2,7,9}$ and $W_{2,13,4}$.

Furthermore, sensitivities are largely confined to a short initial transient and rapidly decay as the system approaches an overdamped regime, where perturbations have minimal impact on the predicted dynamics. This interpretation is supported by the sensitivity ranking based on the $\delta^{msqr}$ metric, which spans several orders of magnitude (from $5 \times 10^{-1}$ to $2 \times 10^{-4}$). Such dispersion highlights a strong imbalance in parameter influence: while some parameters significantly affect the outputs, others contribute negligibly and are therefore difficult to estimate reliably, leading to practical non-identifiability. In addition, the sensitivity matrix is rank-deficient (618 instead of 881), thus reflecting strong correlations among certain parameters.

\begin{figure}[h!]
    \centering
    \includegraphics[width=0.85\textwidth]{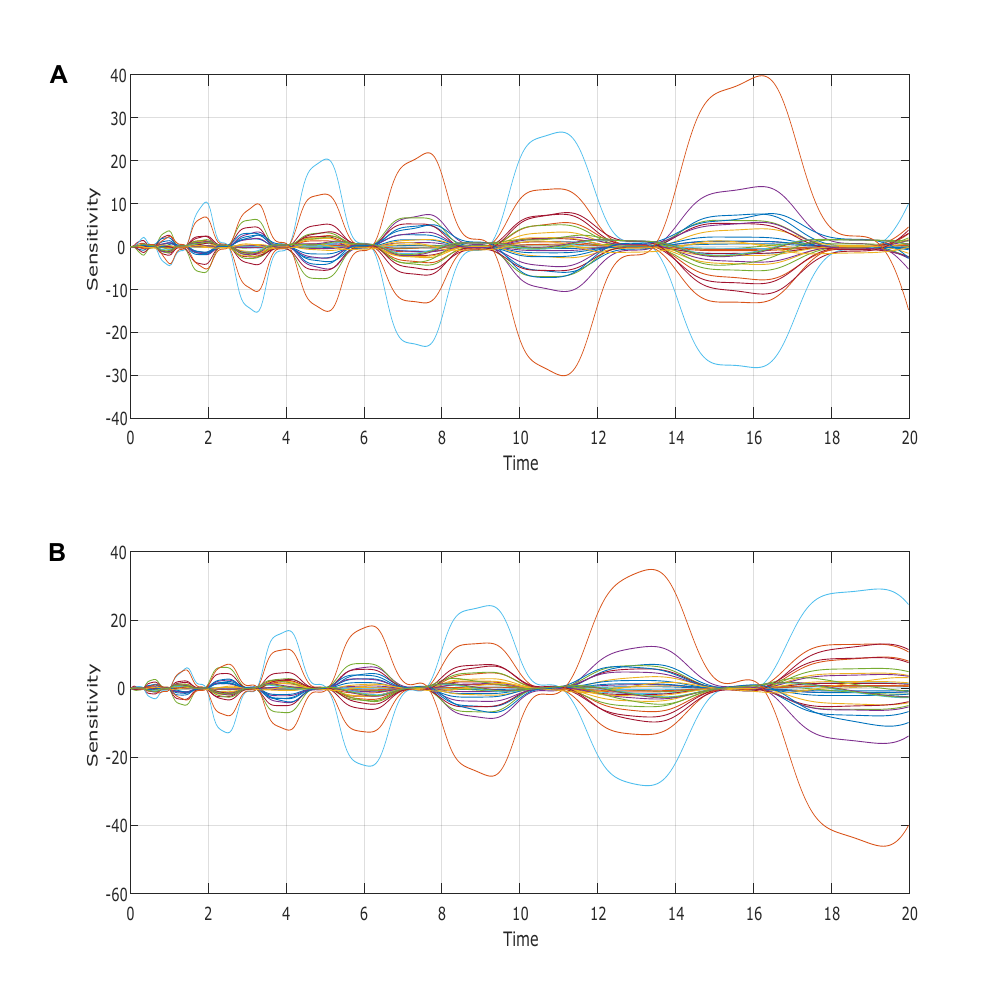}
\caption[\texttt{Hybrid-CDO: }Sensitivity profiles for the selected hybrid iNODE model of the cubic damped oscillator.]{
\textbf{\texttt{Hybrid-CDO: }Sensitivity profiles for the selected hybrid iNODE model of the cubic damped oscillator.}
\textbf{A--B,} Time-dependent sensitivities of the hybrid iNODE outputs \(x_1\) and \(x_2\) with respect to the estimated model parameters, computed using the initial conditions of Experiment~1.
}
    \label{fig:hybrid_CDO_sens_iNODEs}
\end{figure}

A complementary perspective is provided by the sensitivity analysis of the corresponding Hybrid iNODE model, shown in Figure~\ref{fig:hybrid_CDO_sens_iNODEs}. The sensitivities of hybrid iNODE are more evenly distributed across the parameters and display richer temporal variation, with substantially fewer trajectories collapsing toward near-zero values in the integration window (the values of $\delta^{msqr}$ range from $6.89\times10^{-1}$ to $8.15\times10^{-3}$).  This more homogeneous distribution reduces the presence of practically unidentifiable parameters and limits the flat directions in the parameter space. In addition, the hybrid iNODE formulation yields less collinear and more diverse temporal patterns. This reflects a more effective excitation of parameter directions by the experimental data.

The \(\delta_{\mathrm{msqr}}\) values range from \(2\times10^{-1}\) to \(5\times10^{-3}\), with most parameters concentrated between \(10^{-2}\) and \(10^{-1}\). This more balanced sensitivity distribution indicates that the selected iNODE parameterisation is better informed by the available data than the conventional NODE baseline.

\section*{S.5~\texttt{tLV: }Lotka-Volterra with time-dependent interactions}
\addcontentsline{toc}{section}{S.5~\texttt{tLV: }Lotka-Volterra with time-dependent interactions} 
The Lotka--Volterra (LV) framework is a classical and widely used class of mechanistic models in theoretical ecology and systems biology, providing a parsimonious description of interacting populations through coupled nonlinear differential equations \cite{Lotka1925,Volterra1926}. In microbial and cellular systems, these variables are naturally interpreted as population sizes rather than species abundances, and the model is often extended with logistic self-limitation terms to account for finite resources and carrying capacities. For a system of two interacting cell populations, the logistic Lotka--Volterra equations can be written as

\begin{equation}
\left\{
\begin{aligned}
\frac{dN_1(t)}{dt} &= \mu_1 N_1 - \frac{\mu_1 N_1^2}{k_1} - \frac{\mu_1 a_{12}, N_1 N_2}{k_1}\\
\frac{dN_2(t)}{dt} &= \mu_2 N_2 - \frac{\mu_2 N_2^2}{k_2} - \frac{\mu_2 a_{21}, N_1 N_2}{k_2}
\end{aligned}
\right.
\end{equation}

where ($N_i(t)$) denotes the size of population (i), ($\mu_i$) is its intrinsic growth rate, and ($k_i$) is the carrying capacity that accounts for limitations in resources, space, or environmental conditions. The interaction coefficients ($a_{ij}$) quantify the effect of population (j) on the growth of population (i), with positive values indicating inhibitory or competitive effects and negative values indicating cooperative interactions. Despite its simplicity, this framework has been extensively used to analyse interacting biological populations and to infer effective interaction structures from time-series data. At the same time, several recent studies have noted that constant interaction coefficients may be overly restrictive in microbial and cellular systems, where interaction strengths can change over time due to physiological adaptation, environmental variation, or higher-order effects \cite{Hosoda2021Umibato}. This motivates the use of a hybrid formulation in which the mechanistic growth structure is retained, while the interaction terms vary dynamically and are inferred from data. This generalisation allows both the magnitude and the qualitative nature of the interactions between populations to evolve dynamically, enabling the representation of nonstationary biological processes that cannot be captured by models with constant interaction parameters. The time-varying Lotka-Volterra system considered in this study is defined as follows:

\begin{equation}
\left\{
\begin{aligned}
\frac{dN_1(t)}{dt} &= \mu_1 N_1 - \frac{\mu_1 N_1^2}{k_1} - \frac{\mu_1 a_{12}(t)\, N_1 N_2}{k_1}\\
\frac{dN_2(t)}{dt} &= \mu_2 N_2 - \frac{\mu_2 N_2^2}{k_2} - \frac{\mu_2 a_{21}(t)\, N_1 N_2}{k_2}\\
a_{12}(t) &= 2\,\mathrm{e}^{-0.025 t}\\
a_{21}(t) &= -2 + \frac{0.0004\, t^2}{1 + 0.0001\, t^2}
\end{aligned}
\right.
\end{equation}

In this formulation, both interaction terms evolve over time, allowing the balance between inhibition and promotion to change dynamically. Specifically, $a_{12}(t)$ remains positive throughout the time interval considered, but gradually decreases in magnitude, indicating a weakening of the inhibitory effect of the second population on the first. In contrast, $a_{21}(t)$ transitions from negative to positive values, reflecting a qualitative shift in the interaction between the first and second populations. As a result, population dynamics are governed not only by intrinsic growth and self-limitation, but also by time-dependent interspecific influences that cannot be captured by classical Lotka-Volterra models with constant interaction parameters.

To investigate the temporal evolution of species interactions, we adopt a hybrid modelling approach in which the mechanistic Lotka-Volterra framework is retained, while the time dependence of the interaction coefficients is learnt directly from data using neural networks. This corresponds to a latent-parameter dynamics formulation, in which the neural component does not replace the population-state dynamics but instead represents latent, time-varying interaction parameters embedded within an otherwise mechanistic model. The resulting hybrid system is given by:

\begin{equation}
\left\{
\begin{aligned}
\frac{dN_1}{dt} &= \mu_1 N_1 - \frac{\mu_1 N_1^2}{k_1} - \frac{\mu_1 a_{12}(t)\, N_1 N_2}{k_1}\\
\frac{dN_2}{dt} &= \mu_2 N_2 - \frac{\mu_2 N_2^2}{k_2} - \frac{\mu_2 a_{21}(t)\, N_1 N_2}{k_2}\\
a_{12}(t) &= \mathcal{NN}(t;\boldsymbol{\phi}_{\mathrm{NN}})\\
a_{21}(t) &= \mathcal{NN}(t;\boldsymbol{\phi}_{\mathrm{NN}})
\end{aligned}
\right.
\end{equation}

where the neural network $\mathcal{NN}(\cdot;\boldsymbol{\phi}_{\mathrm{NN}})$ parametrizes the latent time-varying interaction coefficients.

This formulation preserves the mechanistic interpretability and biological insight of classical Lotka-Volterra models, while providing the flexibility required to infer nonstationary interaction dynamics directly from data. 

To study this latent-parameter setting, we considered a realistic sparse and noisy experimental scenario and implemented two complementary modelling strategies on the same dataset. First, we applied the proposed identifiability-aware hybrid Neural ODEs framework (Hybrid iNODE). Second, we implemented a conventional hybrid NODE model in Python to provide a direct baseline for comparison.

\subsection*{S.5.1 \texttt{tLV: }Implementation with the hybrid iNODE framework}
\addcontentsline{toc}{subsection}{S.5.1 \texttt{tLV: }Implementation with the hybrid iNODE framework}

In standard Lotka--Volterra models, calibration from coculture data alone often leads to substantial practical identifiability problems because intrinsic growth parameters and interaction coefficients can compensate for one another, making it difficult to disentangle self-limitation from interspecific effects. 

As discussed in \cite{BalsaCanto2025Uncertainty}, this type of mechanistic model benefits from sequential or multiexperiment calibration strategies that improve the conditioning of the estimation problem. When available, monoculture experiments provide direct information on the intrinsic growth rates and carrying capacities of each population in the absence of interactions, thus constraining ($\mu_i$) and ($k_i$) before introducing latent time-varying interactions. This improves the practical identifiability of the remaining mechanistic and neural parameters in the coculture model and supports a more robust recovery of the time-varying interaction functions. 

\begin{table}[h!] 
\centering 
\small 
\caption[\texttt{tLV: }Training and validation scenarios for the Lotka--Volterra model with time-dependent interactions.]{
\textbf{\texttt{tLV: }Training and validation scenarios for the Lotka--Volterra model with time-dependent interactions.}
The scenarios define the training and validation trajectories used to assess recovery of time-varying interaction functions, state prediction and practical identifiability under different sampling and noise conditions.
}\vspace{6pt} 
\label{tab:tLV_noisy_norm_scenarios} 
\begin{tabular}{cccccc} \hline \textbf{Dataset} & \multicolumn{2}{c}{\textbf{Scenario}} & \textbf{Experiment(s)} & \textbf{\begin{tabular}[c]{@{}c@{}}Measurements\\ per state\end{tabular}} & \textbf{\begin{tabular}[c]{@{}c@{}}Total data \\ points\end{tabular}} \\ \hline \multirow{4}{*}{Train} & Monoculture & 1exp & $(N_{1}(0), N_{2}(0)) = (1\times10^{6}, 5\times10^{6})$ & 10 & 20 \\[3mm] & \multirow{3}{*}{Coculture} & \multirow{3}{*}{3exp} & $(N^{(1)}_{1}(0), N^{(1)}_{2}(0)) = (1 \times 10^{8}, 1 \times 10^{8})$ & \multirow{3}{*}{\begin{tabular}[c]{@{}c@{}}10 in each\\ experiment\end{tabular}} & \multirow{3}{*}{60} \\ & & & $(N^{(2)}_{1}(0), N^{(2)}_{2}(0)) = (5 \times 10^{7}, 1 \times 10^{5})$ & & \\ & & & $(N^{(3)}_{1}(0), N^{(3)}_{2}(0)) = (1 \times 10^{5}, 1 \times 10^{7})$ & & \\ [3mm] \multirow{2}{*}{Validation} & \multirow{2}{*}{Coculture} & \multicolumn{1}{l}{val1} & $(N_{1}(0), N_{2}(0)) = (1\times10^{1}, 1\times10^{7})$ & \multicolumn{1}{l}{} & \multicolumn{1}{l}{} \\[3mm] & & val2 & $(N_{1}(0), N_{2}(0)) = (1\times10^{6}, 1\times10^{6})$ & 10 & 20 \\ \hline 
\end{tabular} 
\end{table}

The experimental design adopted in this work is summarised in Table~\ref{tab:tLV_noisy_norm_scenarios}. Experimental data were generated assuming realistic conditions with 10 sampling times and 10\% Gaussian noise. For convenience the Lotka-Volterra latent parameter dynamics were reformulated in logarithmic coordinates. Defining $N_1(t)=10^{y_1(t)}$ and $N_2(t)=10^{y_2(t)}$, the resulting system takes the form:

\begin{equation}
\left\{
\begin{aligned}
N_1(t) &= 10^{y_1(t)}\\
N_2(t) &= 10^{y_2(t)}\\
\frac{dy_1(t)}{dt} &= \frac{\mu_1}{\log(10)}\left( \frac{k_1 - N_1(t) - a_{12}(t)\, N_2(t)}{k_1} \right)\\
\frac{dy_2(t)}{dt} &= \frac{\mu_2}{\log(10)}\left( \frac{k_2 - N_2(t) - a_{21}(t)\, N_1(t)}{k_2} \right)\\
a_{12}(t) &= \mathcal{NN}(t;\boldsymbol{\phi}_{\mathrm{NN}})\\
a_{21}(t) &= \mathcal{NN}(t;\boldsymbol{\phi}_{\mathrm{NN}})
\end{aligned}
\right.
\end{equation}

Tables~\ref{tab:tLV_noisy_norm_iNODEs_gridsearch} and~\ref{tab:tLV_noisy_norm_iNODEs_parameters} compare the performance and parameter estimates of different hybrid iNODE architectures.

\begin{table}[h!]
\centering
\small
\caption[\texttt{tLV: }Candidate hybrid iNODE architectures for the Lotka--Volterra model with time-dependent interactions.]{
\textbf{\texttt{tLV: }Candidate hybrid iNODE architectures for the Lotka--Volterra model with time-dependent interactions.}
For each architecture, the table reports the best-fit objective value, the total number of estimated parameters (\(p\)), the Akaike information criterion (AIC), the training normalised root-mean-square error (NRMSE), and the NRMSE for the recovered time-dependent interaction terms.
The maximum and median confidence intervals across all estimated parameters are also reported, together with the Model Quality Ranking Index (MQRI).
All neural networks use hyperbolic tangent activation functions in the hidden layers and a linear activation function in the output layer.
The highlighted architecture has the lowest MQRI and is selected as the best compromise between state prediction, latent-interaction recovery, parsimony and practical identifiability.
}
\vspace{6pt}
\label{tab:tLV_noisy_norm_iNODEs_gridsearch}
\begin{tabular}{ccccccccc}
\hline
Architecture &
Best fit &
$p$ &
AIC &
\makecell{$\text{NRMSE}_{\text{train}}$ \\ (\%)} &
\makecell{NRMSE $\beta$ \\ (\%)} &
\makecell{Max CI \\ (\%)} &
\makecell{Median CI \\ (\%)} &
MQRI \\
\hline
MLP(1)     & 0.0012 & 10 & -630 & 1.90 & 10.53 & 16.83                 & 2.08  & 7  \\
MLP(1,1)   & 0.0011 & 12 & -626 & 1.91 & 9.75 & 87.44                 & 12.75 & 14 \\
\rowcolor[HTML]{FFCE93} MLP(2) &  0.0009 &  14 &   -634 &   1.89 &   2.89 &   83.10 &
  7.69 &   6 \\
MLP(1,1,1) & 0.0012 & 14 & -621 & 1.92 & 16.57 & -                     & -     & -  \\
MLP(2,1)   & 0.0010 & 15 & -628 & 1.89 & 19.98 & $3.15 \times 10^{12}$ & 36.34 & 13 \\ \hline
\end{tabular}
\end{table} 

All architectures result in an excellent fit quality. However, deeper architectures, such as MLP(1,1) and MLP(2,1), exhibit pronounced identifiability issues with large confidence intervals. Notably, this poor identifiability arises despite excellent reconstruction accuracy, underscoring that goodness of fit alone is insufficient for reliable model selection in hybrid neural differential equation frameworks.

In contrast, single-hidden-layer architectures provide a substantially more robust balance between expressive power and identifiability. Among these, the MLP(2) architecture achieves the lowest best-fit value ($9\times10^{-4}$), the most favourable AIC, and comparatively moderate confidence intervals, resulting in the lowest MQRI in all configurations tested. Compared to deeper alternatives, this architecture avoids excessive parameter uncertainty while maintaining near-optimal predictive accuracy. Table~\ref{tab:tLV_noisy_norm_iNODEs_parameters} presents the parameter values and the corresponding confidence intervals as obtained for the selected architecture.

\begin{table}[h!]
\centering
\caption[\texttt{tLV: }Estimated parameters and confidence intervals for the selected hybrid iNODE model.]{
\textbf{\texttt{tLV: }Estimated parameters and confidence intervals for the selected hybrid iNODE model.}
The table reports the calibrated parameter values \((\theta^*)\) for the Lotka--Volterra model with time-dependent interactions, together with their relative confidence intervals expressed as percentages.
Parameters with large confidence intervals indicate weakly constrained directions and reduced practical identifiability.
}
\vspace{6pt}
\label{tab:tLV_noisy_norm_iNODEs_parameters}
\begin{tabular}{ccc}
\hline
Parameters & $\theta^*$ & CI (\%) \\ \hline
$\mu_1$        & 0.08                 & $2.94 \times 10^{-4}$  \\
$\mu_2$        & 0.02                 & 9.45                   \\
$k_1$          & $1.04 \times 10^{8}$ & 4.17                   \\
$k_2$          & $1.02 \times 10^{7}$ & 7.53                   \\
$W_{111}$    & -0.03                & 83.10                  \\
$b_{11}$      & 0.45                 & 1.93                   \\
$W_{112}$    & -0.02                & 9.33                  \\
$b_{12}$      & 1.42                 & 0.21                   \\
$W_{211}$    & 0.83                 & 41.90                  \\
$W_{221}$    & 0.42                 & 64.80                  \\
$b_{21}$      & 1.13                 & 1.22                   \\
$W_{212}$    & -0.46                & 66.70                 \\
$W_{222}$    & -1.34                & 6.04                  \\
$b_{22}$      & -0.71                & $5.65 \times 10^{-5}$ \\ \hline
\end{tabular}
\end{table}

In contrast, neural-network weights associated with latent interaction functions display larger—but still bounded—uncertainty, reflecting the intrinsic difficulty of inferring time-varying interactions from sparse and noisy population measurements. Importantly, no parameter has a confidence interval greater than 100\%, indicating that all parameters retained are practically identifiable in this final configuration. Bias terms are generally well constrained, suggesting that the overall scale of the learnt interaction dynamics is robustly identified.

\subsection*{S.5.2 \texttt{tLV: }Implementation with the conventional HNODE framework}
\addcontentsline{toc}{subsection}{S.5.2 \texttt{tLV: }Implementation with the conventional HNODE framework} 
Finally, this section presents the results obtained with the conventional HNODE framework. The model was trained using the same data and normalisation as used for the iNODE implementation. The neural component is implemented as a multilayer perceptron with hyperbolic tangent activation functions in the hidden layers and linear activation in the output layer, parameterising the time-dependent interactions. The mechanistic parameters are trained jointly with the neural network and are constrained to remain positive via appropriate reparameterisations. Training is formulated as a trajectory-level optimisation problem across multiple experimental conditions. For each experiment, the hybrid system is numerically integrated forward in time using a differentiable fourth‑order Runge–Kutta scheme that operates directly on the continuous‑time formulation. To improve computational efficiency, neural evaluations of the time‑dependent interactions are batched and reused across all Runge–Kutta substeps within each time interval.

The model parameters are optimised by minimising the mean squared error between the simulated and experimental trajectories in log‑space, aggregated across all experiments. Gradients are propagated through the numerical solver using automatic differentiation, enabling end‑to‑end training with backpropagation. Optimisation is performed using the Adam optimiser with separate learning rates for the mechanistic and neural parameters, gradient clipping to prevent numerical instabilities, and an adaptive learning‑rate scheduler to stabilise convergence. 

For the conventional \texttt{tLV} HNODE baseline, deterministic full-trajectory training was sufficient because the mechanistic structure and the use of multiple experiments constrained the optimisation problem. Candidate architectures were evaluated using NRMSE for both the observed trajectories and the reconstructed time-dependent interaction functions. 

Table~\ref{tab:tLV_noisy_norm_HNODEs_grid} shows that shallow networks fit the observed states but recover the latent interactions less accurately, whereas intermediate architectures improve interaction recovery. Very large networks do not provide further gains, consistent with overparameterisation under sparse noisy data. The MLP(32,16) architecture was selected as the optimal configuration. It achieves the lowest interaction error while maintaining a low training error, offering the best compromise between model complexity, robustness to noise, and accurate recovery of the underlying interaction dynamics.

\begin{table}[h!]
\centering
\caption[\texttt{tLV: }Comparison of conventional HNODE architectures for the Lotka--Volterra model with time-dependent interactions.]{
\textbf{\texttt{tLV: }Comparison of conventional HNODE architectures for the Lotka--Volterra model with time-dependent interactions.}
The table reports the total number of model parameters, the training normalised root-mean-square error (NRMSE), and the NRMSE for the recovered interaction-parameter dynamics for each evaluated architecture.
All neural networks use hyperbolic tangent activation functions in the hidden layers and a linear activation function in the output layer.
The highlighted architecture corresponds to the selected conventional HNODE baseline used for comparison with the hybrid iNODE model.
}
\vspace{6pt}
\label{tab:tLV_noisy_norm_HNODEs_grid}
\begin{tabular}{cccc}
\hline
NN Architecture & Parameters & $ \text{NRMSE}_{\text{train}} (\%) $ & $ \text{NRMSE}_{\text{interaction}} (\%) $ \\ \hline
MLP(4)         & 22         & 0.80      & 12.88              \\
MLP(8)         & 38         & 0.79      & 13.94              \\
MLP(16)        & 70         & 0.78      & 10.68              \\
MLP(4,4)       & 42         & 0.74      & 9.94              \\
MLP(4,8)       & 70         & 0.75      & 11.40              \\
MLP(8,8)       & 110        & 0.76      & 10.75              \\
MLP(16,8)      & 190        & 0.76      & 12.18              \\
MLP(16, 8, 4)  & 218        & 6.94      & 6.83              \\
MLP(16,16)     & 342        & 0.79      & 10.54              \\
\rowcolor[HTML]{FFCE93}
MLP(32,16)     & 630        & 0.75      & 4.14              \\
MLP(32,16,8)   & 750        & 0.81      & 5.21              \\
MLP(32,32)     & 1190       & 0.75      & 6.57              \\
MLP(64,32)     & 2278       & 0.67      & 5.82              \\
MLP(64,64)     & 4422       & 0.76      & 13.60              \\ \hline
\end{tabular}
\end{table}

After training and selecting the best conventional HNODE model, we assessed the uncertainty associated with its estimated parameters. This analysis revealed that the FIM is non-invertible, thereby precluding the direct computation of individual parameter confidence intervals. To investigate the origin of this rank deficiency, a local sensitivity rank analysis was conducted using AMIGO\_LRank. 

\begin{figure}[h!]
    \centering
    \includegraphics[width=0.95\textwidth]{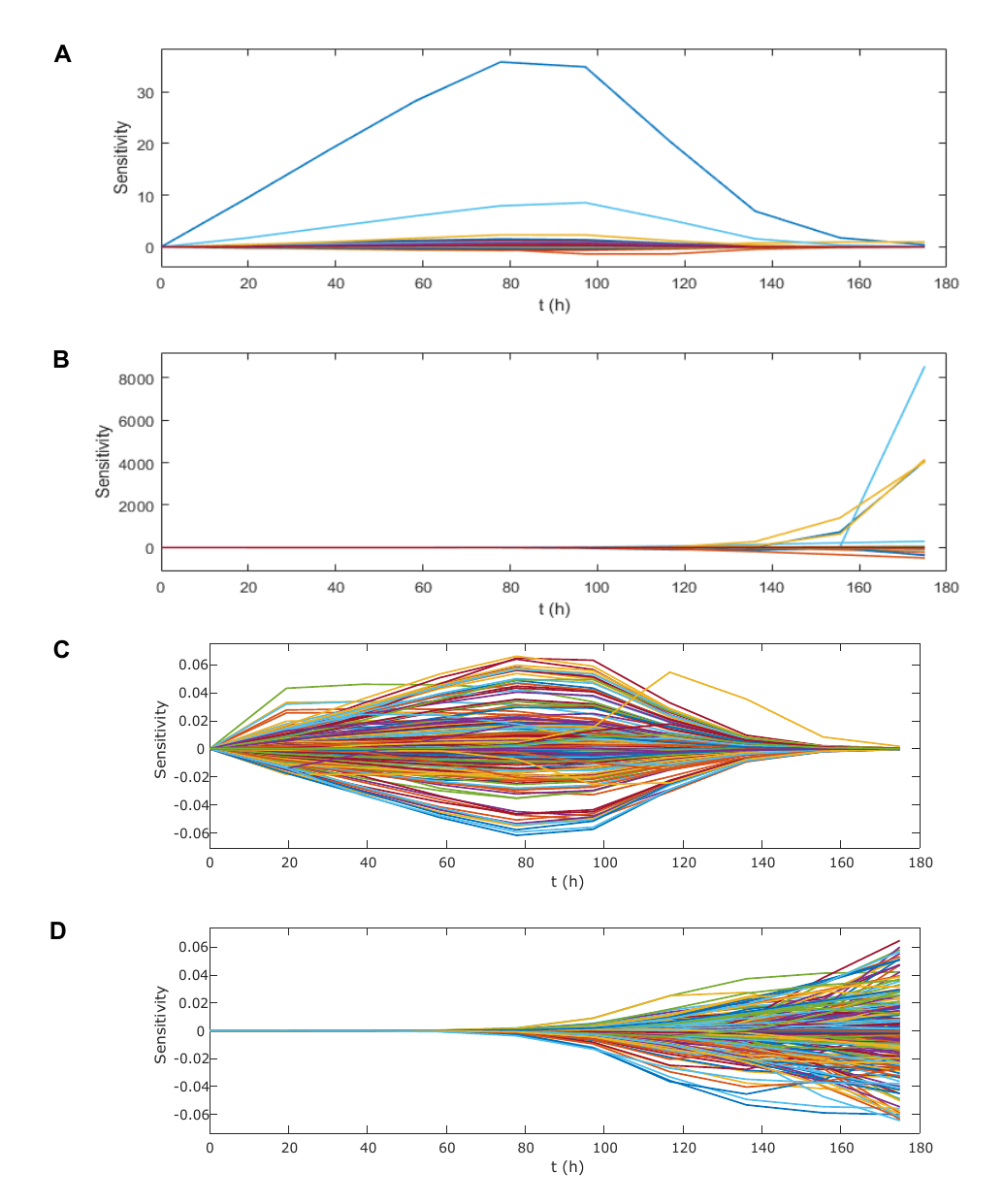}
\caption[\texttt{tLV: }Sensitivity profiles for the conventional HNODE model of the Lotka--Volterra system.]{
\textbf{\texttt{tLV: }Sensitivity profiles for the conventional HNODE model of the Lotka--Volterra system with time-dependent interactions.}
\textbf{A--B,} Time-dependent sensitivities of the HNODE outputs \(y_1\) and \(y_2\) with respect to the estimated model parameters, computed using the initial conditions of Experiment~1.
\textbf{C--D,} Zoomed views of the sensitivities shown in \textbf{A} and \textbf{B}.
The prevalence of near-zero and overlapping sensitivity profiles indicates weakly informed and correlated parameter directions, consistent with limited practical identifiability.
}
    \label{fig:tLV_HNODEs_sens}
\end{figure}

Figure~\ref{fig:tLV_HNODEs_sens} presents the sensitivities with respect to all parameters. The different panels reveal sensitivities spanning several orders of magnitude, with $\delta^{msqr}$ ranging from $10^{-6}$ to $10^{5}$. Representative examples of weakly influential parameters, as identified by low $\delta^{msqr}$ values, include $W_{2,21,15}$, $W_{2,18,15}$, $W_{2,14,15}$, $W_{2,20,15}$, $W_{2,19,15}$, and $W_{2,9,15}$.

This pronounced heterogeneity indicates that only a relatively small subset of parameters exerts a strong influence on the system dynamics, whereas the majority contribute weakly or only over restricted temporal intervals. Consequently, many parameters cannot be reliably estimated from the available data, leading to practical non-identifiability. This limitation is further evidenced by the severe rank deficiency of the sensitivity matrix (40 instead of 630), which highlights the presence of strong linear dependencies among the sensitivity vectors.

\begin{figure}[h!]
    \centering
    \includegraphics[width=0.95\textwidth]{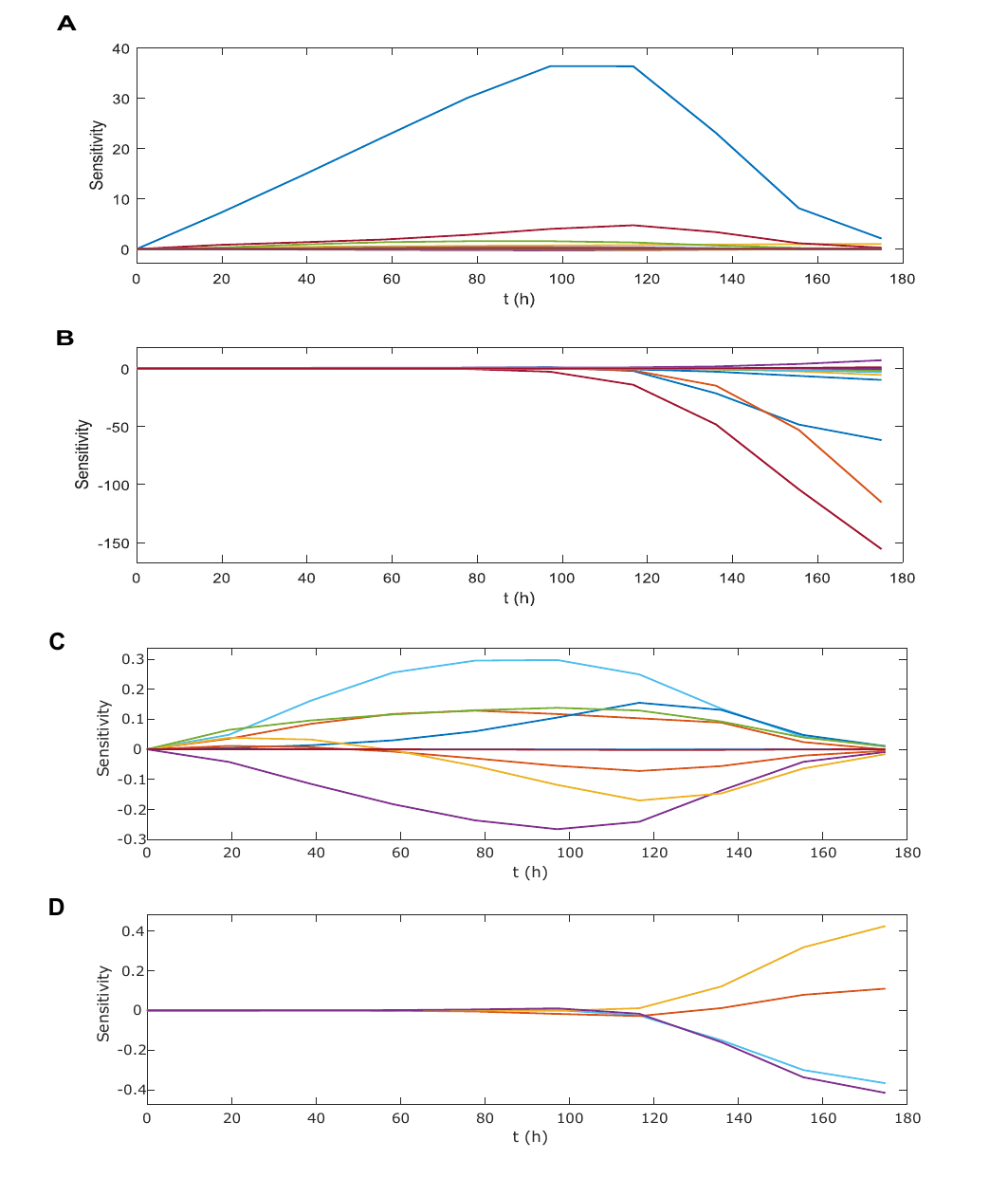}
\caption[\texttt{tLV: }Sensitivity profiles for the selected hybrid iNODE model of the Lotka--Volterra system.]{
\textbf{\texttt{tLV: }Sensitivity profiles for the selected hybrid iNODE model of the Lotka--Volterra system with time-dependent interactions.}
\textbf{A--B,} Time-dependent sensitivities of the hybrid iNODE model outputs \(y_1\) and \(y_2\) with respect to the estimated model parameters, computed using the initial conditions of Experiment~1.
\textbf{C--D,} Zoomed views of the sensitivities shown in \textbf{A} and \textbf{B}.
The sensitivity profiles are more clearly resolved than those of the conventional HNODE baseline, supporting the improved practical identifiability of the selected hybrid iNODE architecture.
}
    \label{fig:tLV_iNODEs_sens}
\end{figure}

The sensitivities of the hybrid iNODE model exhibit a more informative structure, with well-distributed magnitudes and clearly distinguishable temporal patterns (Figure~\ref{fig:tLV_iNODEs_sens}, Panels~\textbf{A--B}). Sensitivity trajectories remain active throughout the observation window, with minimal collapse towards near-zero values, as is clearly observed in Panels~\textbf{C--D}, where a closer inspection of the zoomed-in parameter sensitivities reveals the absence of parameters with negligibly small sensitivities, in contrast to the behaviour observed in Figure~\ref{fig:tLV_HNODEs_sens}. This improvement is consistent with the sensitivity ranking based on the $\delta^{msqr}$ metric, whose values span from $2\times10^{1}$ to $1\times10^{-2}$. This indicates that multiple parameter directions are effectively informed by the data, in contrast to the pronounced imbalance observed in the HNODE formulation.

\pagebreak

Distinct parameters induce separable perturbations in both state variables, reflecting a substantial reduction in redundancy and sensitivity collinearity. Consequently, the hybrid formulation yields a compact yet informative sensitivity structure with greater diversity over time and across parameters. This behaviour directly stems from the identifiability-aware design of the hybrid iNODE framework, which combines mechanistic structure with a constrained neural representation.

By limiting redundant degrees of freedom and aligning model complexity with the information content of the data, the hybrid formulation increases the effective rank of the sensitivity matrix and enables reliable inversion of the Fisher Information Matrix (FIM). As a result, well-defined confidence intervals can be computed for all model parameters, enabling uncertainty quantification.

\pagebreak

Finally, after evaluating the parameter sensitivities of both the conventional HNODE and the hybrid iNODE formulations, an additional analysis was carried out to further investigate the consequences of the identifiability limitations observed in the HNODE model and their impact on predictive uncertainty.

The severe rank deficiency observed in the HNODE sensitivity matrix has direct implications for uncertainty quantification. Since the Fisher Information Matrix (FIM) is singular, the classical asymptotic covariance approximation cannot be applied and individual parameter confidence intervals cannot be obtained through inversion of the FIM. To further characterise parameter uncertainty despite these limitations, a practical identifiability analysis was performed using the resampling bootstrap approach implemented in \texttt{AMIGO\_RIdent}. 

Two independent bootstrap analyses were conducted, each consisting of 500 resamplings. After completing both analyses, the confidence intervals estimated for each parameter were examined and, for every parameter, the widest interval obtained across the two bootstrap procedures was selected. This ensured a conservative and robust quantification of parameter uncertainty, consistent with the identifiability limitations revealed by the sensitivity analysis.

The resulting confidence intervals exhibited an extremely broad distribution, ranging from $2.5\times10^{-13}\%$ to $1590\%$, with a median value of $9.1\%$. This remarkable variability indicates that the available experimental data provide highly uneven information across the parameter space: while some parameters are tightly constrained, others display very large uncertainty regions.

Such behaviour is characteristic of practical non-identifiability, whereby distinct parameter combinations can generate nearly indistinguishable model predictions. The coexistence of extremely small and extremely large confidence intervals reveals that only a subset of parameter directions is effectively informed by the data, whereas many others remain weakly constrained owing to redundancy and strong parameter correlations. These findings are fully consistent with the heterogeneous sensitivity structure observed in Figure~\ref{fig:tLV_HNODEs_sens} and with the severe rank deficiency of the sensitivity matrix. Taken together, they provide independent evidence that the conventional HNODE formulation suffers from substantial practical non-identifiability despite achieving an accurate fit to the experimental observations.

The bootstrap analysis further demonstrates that multiple parameter combinations can reproduce the observed dynamics with comparable accuracy, implying that the predictive behaviour of the model is supported by a highly non-unique parameterisation. This result highlights a fundamental limitation of high-dimensional neural representations and reinforces the advantages of the identifiability-aware hybrid iNODE formulation.

Once confidence intervals had been obtained for both model formulations, the corresponding prediction uncertainty bands were generated following the ensemble-based uncertainty propagation methodology described in detail in the Methods section.

\section*{S.6~\texttt{tSIR: }Susceptible-Infectious-Recovered model with time-dependent transmission rate}
\addcontentsline{toc}{section}{S.6~\texttt{tSIR: }Susceptible-Infectious-Recovered model with time-dependent transmission rate} 

The Susceptible-Infectious-Recovered (SIR) model is one of the most widely used mechanistic frameworks in epidemiology, owing to its conceptual simplicity, interpretability, and low computational cost \citep{Kermack1927,Brauer2008}. By partitioning the population into susceptible (S), infectious (I), and recovered (R) compartments, the model provides an intuitive description of the spread of infection and recovery processes, with parameters directly linked to epidemiologically meaningful quantities such as transmission and recovery rates. These features make SIR-type models particularly well suited for scenario analysis, policy assessment, and qualitative understanding of epidemic dynamics \citep{Brauer2008}.

Despite these advantages, the classical SIR formulation relies on strong simplifying assumptions, including homogeneous mixing of the population, constant model parameters, and the absence of behavioural or environmental feedback mechanisms. In real epidemics, contact patterns, intervention measures, and population behaviour typically change over time, limiting the ability of standard SIR models to capture time-dependent transmission rate. From an inference perspective, these limitations are further compounded by identifiability problems, as distinct parameter configurations can yield similar epidemic trajectories, particularly when only partial observations are available \citep{Tuncer2018,Massonis2020,Kharazmi2021}. This makes SIR-type models a useful benchmark for the present work: they are simple enough to remain interpretable, yet sufficiently challenging to test whether hybrid neural differential equations can recover latent time-dependent transmission rate under partial observability.

In its classical form, the SIR model is defined as follows:

\begin{equation}
\left\{
\begin{aligned}
\frac{dS}{dt} &= -\dfrac{\beta\, S\, I}{N}\\
\frac{dI}{dt} &= \dfrac{\beta\, S\, I}{N} - \gamma\, I\\
\frac{dR}{dt} &= \gamma\, I\\
N &= S + I + R
\end{aligned}
\right.
\end{equation}

where $\beta$ denotes the time-dependent transmission rate, $\gamma$ is the recovery rate, and $N$ is the total population size, assumed to remain constant.

To better capture nonstationary epidemic dynamics, the SIR framework can be extended to explicitly allow the time-dependent transmission rate to vary over time. In this study, we consider a mechanistic SIR model with a time-dependent transmission rate of the following form:

\begin{equation}
\left\{
\begin{aligned}
\frac{dS}{dt} &= -\dfrac{\beta(t)\, S\, I}{N}\\
\frac{dI}{dt} &= \dfrac{\beta(t)\, S\, I}{N} - \gamma\, I\\
\frac{dR}{dt} &= \gamma\, I\\
N &= S + I + R\\
\beta(t) &= b_{\min} + (b_{\max} - b_{\min})
\left[
\frac{1}{1 + \mathrm{e}^{-k_1 (t - t_1)}}
\left(
1 - \frac{1}{1 + \mathrm{e}^{-k_2 (t - t_2)}}
\right)
\right]
\end{aligned}
\right.
\end{equation}

where $\beta(t)$ denotes the transmission rate that varies over time and $\gamma$ is assumed to remain constant. This functional form allows the transmission rate to increase and decrease smoothly over time, enabling the representation of epidemic phases such as the onset and relaxation of containment measures.

Although this formulation preserves the interpretability of the SIR framework and enables the modelling of time-varying transmission dynamics, prescribing a fixed functional form for $\beta(t)$ can introduce modelling bias and limit flexibility when the true temporal evolution of transmission is unknown. These considerations motivate the exploration of data-driven and hybrid modelling approaches that infer time-varying epidemiological parameters directly from data, while retaining the mechanistic structure of compartmental epidemic models.

To address the limitations of classical SIR formulations with fixed transmission parameters, we introduce a hybrid SIR model that retains the mechanistic compartmental structure while inferring the temporal evolution of the transmission rate directly from the data using neural networks. This corresponds to a latent parameter dynamics formulation in which the neural component does not replace the epidemiological state equations but instead represents a latent, time-varying parameter embedded within an otherwise mechanistic model. The resulting hybrid SIR system is given by:

\begin{equation}
\left\{
\begin{aligned}
\frac{dS}{dt} &= -\dfrac{\beta(t)\, S\, I}{N}\\
\frac{dI}{dt} &= \dfrac{\beta(t)\, S\, I}{N} - \gamma\, I\\
\frac{dR}{dt} &= \gamma\, I\\
N &= S + I + R\\
\beta(t) &= \mathcal{NN}\big(t;\boldsymbol{\phi}_{\mathrm{NN}}\big)
\end{aligned}
\right.
\end{equation}

where $\mathcal{NN}(\cdot;\boldsymbol{\phi}_{\mathrm{NN}})$ denotes a neural network parameterized by weights and biases $\boldsymbol{\phi}_{\mathrm{NN}}$ that models the latent time-varying transmission rate.

In addition to the time-dependent transmission coefficient, this case study is characterised by partial observability: only the infected compartment is observed, reflecting common epidemiological scenarios in which incidence or prevalence data are recorded while the susceptible and recovered compartments remain unobserved.

The different hybrid models were calibrated and evaluated using synthetic data generated from the SIR system under partial observability. All datasets were generated using fixed epidemiological parameters, namely a recovery rate $\gamma = 0.1~\mathrm{days}^{-1}$ and a transmission rate varying over time characterised by $b_{\min} = 0.09~\mathrm{days}^{-1}$, $b_{\max} = 0.28~\mathrm{days}^{-1}$, $t_1 = 15~\mathrm{days}$, $k_1 = 0.4~\mathrm{days}^{-1}$, $t_2 = 70~\mathrm{days}$ and $k_2 = 0.20~\mathrm{days}^{-1}$. In total, 60 synthetic measurements of the infected population were used, representing a data-limited epidemiological setting.

To assess the influence of experimental design and data richness on parameter identifiability and predictive performance, three different training configurations were considered as summarised in Table \ref{tab:tSIR_scenarios}. Note that \texttt{1exp\_ext} regards the same experiment as \texttt{1exp} but longer in time and with 20 additional sampling times, while the case \texttt{2exp} regards a scenario in which two different initial states are used.

\begin{table}[H]
\centering
\small
\caption[\texttt{tSIR: }Training and validation scenarios for the SIR model with time-dependent transmission rate.]{
\textbf{\texttt{tSIR: }Training and validation scenarios for the SIR model with time-dependent transmission rate.}
The scenarios define the training and validation trajectories used to assess recovery of the latent transmission-rate dynamics, state prediction and practical identifiability under partial observability.
}
\vspace{6pt}
\label{tab:tSIR_scenarios}
\begin{tabular}{ccccc}
\hline
\textbf{Dataset} &
  \textbf{Scenario} &
  \textbf{Experiment(s)} &
  \textbf{\begin{tabular}[c]{@{}c@{}}Measurements\\ (only I obsevable)\end{tabular}} &
  \textbf{\begin{tabular}[c]{@{}c@{}}Total data \\ points\end{tabular}} \\ \hline
\multirow{4}{*}{Train} & 1exp       & $(S^{(1)}(0), I^{(1)}(0), R^{(1)}(0)) = (100\,000, 100, 0)$ & 60 & 60 \\[3mm]
                       & 1exp\_ ext & $(S^{(1)}(0), I^{(1)}(0), R^{(1)}(0)) = (100\,000, 100, 0)$ & 80 & 80 \\[3mm]
 &
  \multirow{2}{*}{2exp} &
  $(S^{(1)}(0), I^{(1)}(0), R^{(1)}(0)) = (100\,000, 100, 0)$ &
  \multirow{2}{*}{\begin{tabular}[c]{@{}c@{}}40 in each \\ experiment\end{tabular}} &
  \multirow{2}{*}{80} \\
                       &            & $(S^{(2)}(0), I^{(2)}(0), R^{(2)}(0)) = (50\,000, 500, 0)$  &    &    \\[3mm]
Validation             & val        & $(S(0), I(0), R(0)) = (80\,000, 80, 0)$                     & 60 & 60 \\ \hline
\end{tabular}
\end{table}

\subsubsection*{S.6.1 \texttt{tSIR: }Implementation with the hybrid iNODE framework}
\addcontentsline{toc}{subsection}{S.6.1 \texttt{tSIR: }Implementation with the hybrid iNODE framework} 
In this section, we present the results obtained with the proposed hybrid iNODE formulation applied to the SIR system. Tables~\ref{tab:tSIR_iNODES_gridsearch} and~\ref{tab:tSIR_iNODES_parameters} summarise the results of the calibration and model selection procedures, highlighting the performance of the method in terms of parameter identifiability, reconstruction accuracy, and predictive capability. The model design and selection process was carried out using data from a single experiment (\texttt{1exp}). Subsequently, the parameter estimation procedure for the selected model was repeated using the extended dataset (\texttt{1exp\_ext}) as well as the multi-experiment dataset (\texttt{2exp}).

In particular, Table~\ref{tab:tSIR_iNODES_gridsearch} compares multiple hybrid iNODE architectures based on a combination of fit accuracy, model complexity, and practical identifiability criteria. This integrated evaluation enables a balanced assessment of both predictive performance and parameter reliability.

\begin{table}[h!]
\centering
\small
\caption[\texttt{tSIR: }Candidate hybrid iNODE architectures for the SIR model with time-dependent transmission rate.]{
\textbf{\texttt{tSIR: }Candidate hybrid iNODE architectures for the SIR model with time-dependent transmission rate.}
For each architecture, the table reports the best-fit objective value, the total number of estimated parameters (\(p\)), the Akaike information criterion (AIC), the training normalised root-mean-square error (NRMSE), and the maximum and median confidence intervals across all estimated parameters.
The Model Quality Ranking Index (MQRI) is also reported for each configuration.
All neural networks use hyperbolic tangent activation functions in the hidden layers and a linear activation function in the output layer.
The highlighted architecture has the lowest MQRI and is selected as the best compromise between state prediction, parsimony and practical identifiability.
}
\vspace{6pt}
\label{tab:tSIR_iNODES_gridsearch}
\begin{tabular}{ccccccccc}
\hline
Architecture &
Best fit &
$p$ &
AIC &
\makecell{$\text{NRMSE}_{\text{train}}$ \\ (\%)} &
\makecell{NRMSE $\beta$ \\ (\%)} &
\makecell{Max CI \\ (\%)} &
\makecell{Median CI \\ (\%)} &
MQRI\\
\hline
MLP(1)   & $1.60 \times 10^{-3}$ & 5  & -624  & $0.52$ & 65.31 & $1.01 \times 10^{6}$ & $1.17\times 10^{2}$ & 33 \\
MLP(1,1) & $1.58 \times 10^{-3}$ & 7  & -621  & $0.52$ & 65.60 & $3.57 \times 10^{6}$ & $4.15\times 10^{2}$ & 36 \\
\rowcolor[HTML]{FFCE93} MLP(2)   & $4.58 \times 10^{-7}$ & 8  & -1107 & $8.77 \times 10^{-3}$ & 6.20 & $1.76 \times 10^{3}$ & 0.96   & 5  \\
MLP(2,1) & $1.29 \times 10^{-6}$ & 10 & -1041 & $1.47 \times 10^{-2}$ & 4.99 & $3.35 \times 10^{2}$ & 12.79  & 14 \\
MLP(1,2) & $3.97 \times 10^{-5}$ & 10 & -836  & $8.16 \times 10^{-2}$ & 18.77 & $2.30\times 10^{4}$  & 98.16  & 27 \\
MLP(3)   & $4.81 \times 10^{-7}$ & 11 & -1099 & $8.98 \times 10^{-3}$ & 8.07 & $2.30\times 10^{9}$  & 17.74  & 18 \\
MLP(3,1) & $1.22 \times 10^{-6}$ & 13 & -1038 & $1.45 \times 10^{-2}$ & 4.00 & $1.18 \times 10^{6}$ & 31.25  & 24 \\
MLP(1,3) & $2.69 \times 10^{-6}$ & 13 & -991  & $2.13 \times 10^{-2}$ & 13.94 & $2.37 \times 10^{3}$ & 13.58  & 20 \\
MLP(4)   & $4.81 \times 10^{-7}$ & 14 & -1093 & $8.99 \times 10^{-3}$ & 6.54 & $7.70\times 10^{10}$ & 16.89  & 20 \\
MLP(2,2) & $8.10 \times 10^{-7}$ & 14 & -1061 & $1.17 \times 10^{-2}$ & 5.29 & $2.52 \times 10^{4}$ & $4.75 \times 10^{2}$ & 23 \\ \hline
\end{tabular}
\end{table}

The comparison of candidate architectures reveals a clear trade-off between model expressiveness and practical identifiability. The MLP(2) architecture yields the lowest best-fit value with a median relative confidence interval below $1\%$. Consistently, this architecture achieves the most favourable AIC value and the lowest Model Quality Ranking Index, reflecting the best overall balance between predictive accuracy, model complexity, and identifiability.

Although more complex architectures achieve comparable fitting performance, they do not systematically improve identifiability and frequently exhibit large confidence intervals, indicative of overparameterisation. These results underscore the importance of identifiability-aware calibration and establish the MLP(2) architecture as the optimal configuration for the SIR system. 

The Table~\ref{tab:tSIR_iNODES_parameters} reports the estimated parameters and relative confidence intervals for this hybrid iNODE-SIR model in three experimental configurations, enabling a direct evaluation of practical identifiability as a function of data richness and experimental design.

\begin{table}[h!]
\centering
\caption[\texttt{tSIR: }Estimated parameters and confidence intervals for the hybrid iNODE configuration of the SIR model.]{
\textbf{\texttt{tSIR: }Estimated parameters and confidence intervals for the hybrid iNODE configuration of the SIR model with time-dependent transmission rate.}
The table reports the calibrated parameter values \((\theta^*)\) for the \texttt{1exp}, \texttt{1exp\_ext} and \texttt{2exp} hybrid iNODE configurations, together with their relative confidence intervals expressed as percentages.
Parameters highlighted in blue have large confidence intervals, indicating weak practical identifiability under the corresponding experimental design.
}
\vspace{6pt}
\label{tab:tSIR_iNODES_parameters}
\begin{tabular}{ccccccc}
\hline
           & \multicolumn{2}{c}{\texttt{1exp}} & \multicolumn{2}{c}{\texttt{1exp\_ext}} & \multicolumn{2}{c}{\texttt{2exp}} \\
Parameter  & $\theta^*$ & CI (\%) & $\theta^*$ & CI (\%) & $\theta^*$ & CI (\%) \\ \hline
$\gamma$      & 0.098          & 0.23         & 0.098             & 0.20           & 0.100          & 0.01         \\
$W_{111}$     & \cellcolor[HTML]{96FFFB}0.193          & \cellcolor[HTML]{96FFFB}$1.76 \times 10^3$      & \cellcolor[HTML]{96FFFB}0.264             & \cellcolor[HTML]{96FFFB}$4.74 \times 10^2$         & 0.099          & 2.46         \\
$b_{11}$      & -2.880         & 15.67       & -4.105            & 8.86          & -6.944         & 0.09        \\
$W_{112}$     & -0.091         & 0.93        & 0.092             & 15.19          & 0.198          & 0.42         \\
$b_{12}$      & 6.323          & 0.31         & -6.391            & 0.29          & 2.975         & $4.13 \times 10^{-4}$ \\
$W_{211}$     & 0.096          & 0.28         & 0.092             & 0.04           & -0.095         & 17.09       \\
$W_{221}$     & 0.105          & 1.61         & -0.104            & 1.02          & 0.095          & 0.04         \\
$b_{21}$      & 0.078          & $3.87 \times 10^{-4}$ & 0.083   & 2.08           & 0.090          & 0.47         \\ \hline
\end{tabular}
\end{table}

In the \texttt{1exp} configuration, which relies on a single epidemic experiment with limited temporal coverage, the recovery rate $\gamma$ is accurately estimated with a narrow confidence interval. In contrast, parameters associated with the neural representation of the time-varying transmission rate display heterogeneous identifiability. Notably, the weight $W_{111}$ exhibits an extremely large confidence interval ($CI>1700\%$), indicating severe practical non-identifiability. This behaviour suggests that a single epidemic trajectory does not provide sufficient information to uniquely constrain all components of the embedded neural dynamics. Extending the observation window in \texttt{1exp\_ext} yields partial improvements. Although the confidence interval of $W_{111}$ is substantially reduced compared to \texttt{1exp}, it remains large, and the overall identifiability structure remains qualitatively similar. These results indicate that increased temporal coverage alone is insufficient to fully resolve the identifiability limitations of the neural interaction component.

\pagebreak

A marked improvement is observed in the \texttt{2exp} configuration, where the inclusion of a second experiment with distinct initial conditions provides richer dynamical information. In this setting, all parameters—including those previously affected by severe uncertainty—exhibit tight confidence intervals. In particular, the confidence interval of $W_{111}$ is reduced by approximately three orders of magnitude, indicating that the time-varying transmission dynamics are now well constrained. Importantly, this improvement is achieved without increasing the total number of observations compared with the single-experiment case, highlighting the critical role of experimental diversity over data volume alone.

\begin{figure}[h!]
    \centering
    \includegraphics[width=0.95\linewidth]{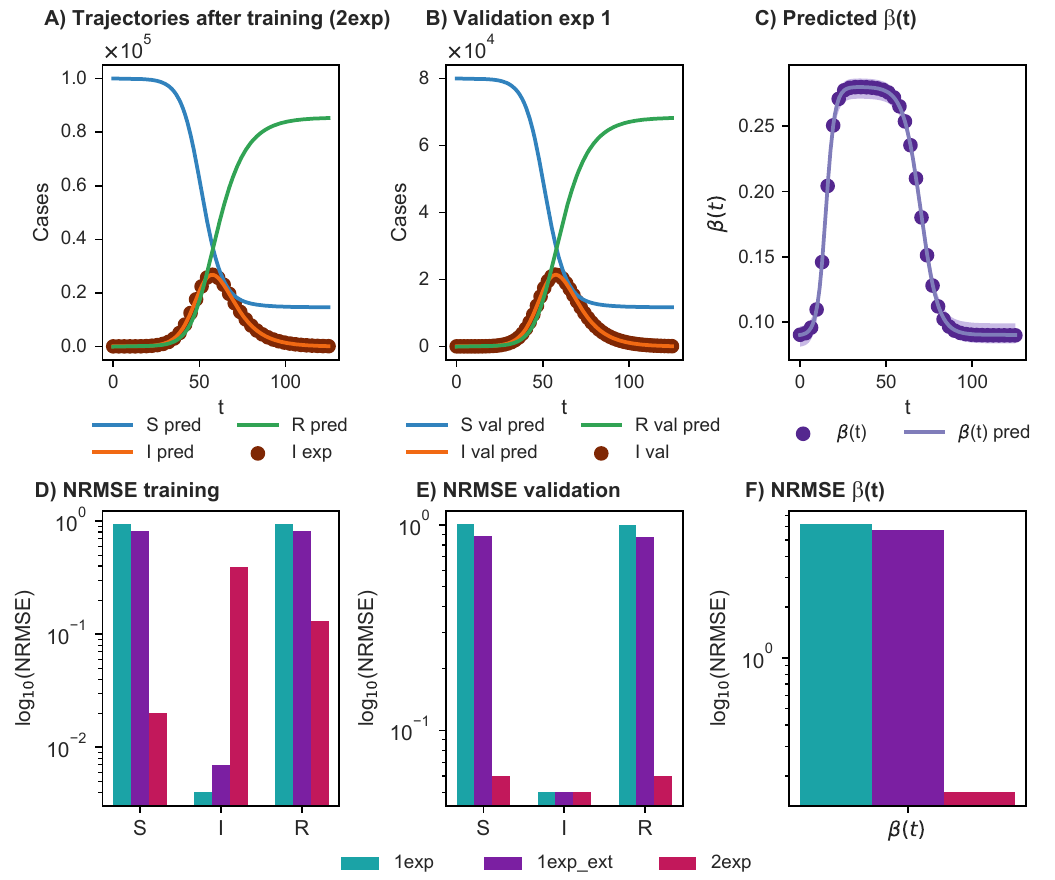}
\caption[\texttt{tSIR: }Predictions, inferred transmission rate and NRMSE comparison for the hybrid iNODE SIR model.]{
\textbf{\texttt{tSIR: }Predictions, inferred transmission rate and NRMSE comparison for the hybrid iNODE SIR model with time-dependent transmission rate.}
\textbf{A--B,} Predicted epidemic dynamics for the first training experiment in the \texttt{2exp} configuration, with initial conditions \((S^{(1)}(0), I^{(1)}(0), R^{(1)}(0)) = (100\,000, 100, 0)\), and for the validation experiment.
\textbf{C,} Time-varying transmission rate inferred by the \texttt{2exp} hybrid iNODE configuration.
In \textbf{A--C}, points denote synthetic observations, solid lines denote model predictions, and shaded regions indicate uncertainty bands.
\textbf{D,} Normalised root-mean-square error (NRMSE) of the three hybrid iNODE configurations across all states in the training experiments.
\textbf{E,} NRMSE for the validation experiment.
\textbf{F,} NRMSE of the inferred transmission rate for the three configurations.
All configurations achieve low prediction errors for the observed states, but the two-experiment configuration yields the most accurate recovery of the latent transmission rate, highlighting the benefit of increased experimental information.
}
    \label{fig:tSIR_iNODEs_models}
\end{figure}

Consistent with these improvements in identifiability, \texttt{2exp} accurately reproduces the training dynamics and achieves the lowest validation NRMSE values (Figure~\ref{fig:tSIR_iNODEs_models}, Panels~\textbf{A}, \textbf{B}, \textbf{D}, \textbf{E}). By contrast, \texttt{1exp} and \texttt{1exp\_ext} exhibit substantially larger validation errors—particularly for susceptible and recovered populations—which worsen under extrapolation due to limited dynamical variability in the training data. 

Although extending the observation window improves the quality of the fitting, it does not translate into improved predictive performance. In addition, predictive uncertainty under \texttt{2exp} is minimal, with narrow uncertainty bands consistent with the tight parameter confidence intervals.

The reconstructed transmission rate $\beta(t)$ under \texttt{2exp} is smooth, accurate, and accompanied by similarly narrow uncertainty bands (Figure~\ref{fig:tSIR_iNODEs_models}, Panel~\textbf{C}). This configuration achieves validation errors more than an order of magnitude lower than those obtained under alternative training designs (Figure~\ref{fig:tSIR_iNODEs_models}, Panel~\textbf{F}).

\begin{figure}[H]
    \centering
    \includegraphics[width=0.75\linewidth]{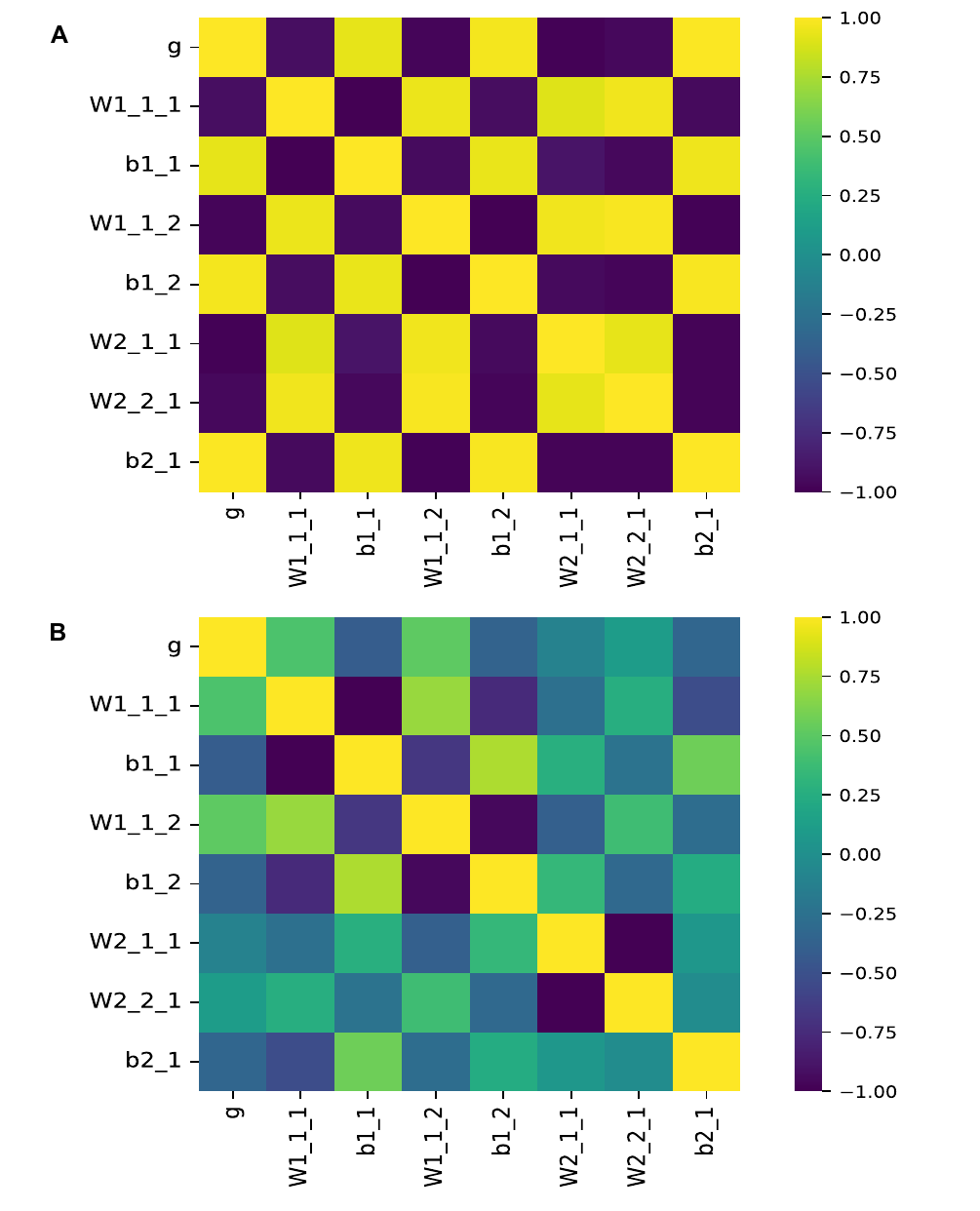}
\caption[\texttt{tSIR: }Parameter correlation matrices for two hybrid iNODE configurations of the SIR model.]{
\textbf{\texttt{tSIR: }Parameter correlation matrices for two hybrid iNODE configurations of the SIR model with time-dependent transmission rate.}
\textbf{A,} Correlation matrix for the \texttt{1exp} hybrid iNODE configuration.
\textbf{B,} Correlation matrix for the \texttt{2exp} hybrid iNODE configuration.
Rows and columns correspond to the neural-network weights and biases, together with the mechanistic parameter \(\gamma\).
The colour scale shows pairwise Pearson correlation coefficients; values close to \(\pm 1\) indicate strong linear dependencies, whereas values near zero indicate weak correlations.
Compared with the \texttt{1exp} configuration, the \texttt{2exp} configuration exhibits fewer strong parameter correlations, consistent with improved practical identifiability.
}
    \label{fig:tSIR_iNODEs_corr}
\end{figure}

The parameter correlation matrices further explain this difference (Figure~\ref{fig:tSIR_iNODEs_corr}). The \texttt{1exp} configuration exhibits strong positive and negative off-diagonal correlations, indicating compensatory effects among parameters. In the \texttt{2exp} configuration, these correlations are markedly reduced, particularly among parameters governing the neural representation of the transmission rate. This decorrelation is consistent with the narrower confidence intervals reported in Table \ref{tab:tSIR_iNODES_parameters} and indicates that the individual parameter effects are more effectively distinguished by the two-experiment design.

\pagebreak

The \texttt{2exp} configuration enables consistent and reliable recovery of both the epidemic trajectories and the underlying time-varying transmission rate. These findings demonstrate that practical identifiability in hybrid mechanistic-neural models is fundamentally governed by experimental design. While extending observation windows can partially reduce parameter uncertainty, the inclusion of heterogeneous experiments probing distinct dynamical regimes is essential to robustly identify embedded neural components. Such experimental diversity is therefore critical to achieving stable inference and reliable uncertainty quantification in partially observable epidemic systems.

\subsubsection*{S.6.2 \texttt{tSIR: }Implementation with the conventional HNODE framework}
\addcontentsline{toc}{subsection}{S.6.2 \texttt{tSIR: }Implementation with the conventional HNODE framework} 
In this section, we report the results obtained using a conventional HNODE formulation of the tSIR model. Model training was performed using the datasets defined in the configuration \texttt{2exp}, which comprises two epidemic trajectories, each providing 40 measurements of the infected population.

The conventional tSIR HNODE baseline was trained using the \texttt{2exp} dataset, comprising two epidemic trajectories with 40 infected-population measurements each. Training was formulated as a continuous-time trajectory-level optimisation problem with shared model parameters across experiments. Positivity of \(\beta(t)\) and \(\gamma\) was enforced through softplus reparameterisation, and the model was calibrated by minimising a weighted trajectory error with regularisation terms promoting smooth \(\beta(t)\) profiles and biologically plausible \(\gamma\) values. Candidate architectures were then compared using NRMSE for both \(I(t)\) and the reconstructed \(\beta(t)\).

The results summarised in Table~\ref{tab:tSIR_HNODEs_grid} reveal a strong dependence of model performance on the neural architecture used to approximate the time‑varying transmission rate $\beta(t)$. Shallow networks with limited expressive capacity, such as the single‑layer architectures MLP(2) and MLP(4), can capture overall epidemic trends but produce relatively large reconstruction errors, particularly in the inferred transmission dynamics. This behaviour indicates insufficient representational capacity to accurately resolve the temporal variability of $\beta(t)$.

\begin{table}[h!]
\centering
\caption[\texttt{tSIR: }Comparison of conventional HNODE architectures for the SIR model with time-dependent transmission rate.]{
\textbf{\texttt{tSIR: }Comparison of conventional HNODE architectures for the SIR model with time-dependent transmission rate.}
For each neural architecture used to represent the time-varying transmission rate \(\beta(t)\), the table reports the total number of trainable parameters and the normalised root-mean-square error (NRMSE) for the infected population trajectory \(I(t)\) and the reconstructed transmission rate \(\beta(t)\).
All neural networks use hyperbolic tangent activation functions in the hidden layers.
The highlighted architecture corresponds to the selected conventional HNODE baseline based on overall predictive performance.
}
\vspace{6pt}
\label{tab:tSIR_HNODEs_grid}
\begin{tabular}{cccc}
\hline
NN Architecture & Parameters & NRMSE $I(t)$ (\%) & NRMSE $\beta(t)$ (\%)\\ \hline
MLP(2)   & 7   & 3.88 & 36.50 \\
MLP(4)   & 13  & 3.51 & 33.02 \\
MLP(8)   & 25  & 2.97 & 29.75 \\
MLP(4,2) & 21  & 2.27 & 25.42 \\
\rowcolor[HTML]{FFCE93}
MLP(4,4) & 34  & 1.40 & 20.35 \\
MLP(8,4) & 57  & 2.22 & 24.68 \\ \hline
\end{tabular}
\end{table}

Increasing network depth and width leads to marked improvements in both state prediction and parameter recovery. Architectures such as MLP(4,2) and MLP(4,4) achieve substantially lower normalised root mean squared errors (NRMSE), reflecting an improved balance between model flexibility and generalisation. Among the tested configurations, the MLP(4,4) architecture yields the lowest errors for both the infected population $I(t)$ and the reconstructed transmission rate $\beta(t)$, while maintaining a moderate number of trainable parameters. This architecture is therefore selected as the optimal HNODE model.

Although the selected conventional HNODE model achieves a reasonable fit to the observed epidemic trajectories, this does not translate into reliable recovery of the latent transmission rate. In particular, the reconstructed ($\beta(t)$) remains substantially less accurate than that obtained with the corresponding hybrid iNODE model, despite the satisfactory agreement with the measured infected population. This indicates that good trajectory-level fit under partial observability is not sufficient to ensure reliable inference of the underlying transmission dynamics. A plausible explanation is that the conventional HNODE architecture contains more trainable parameters than can be robustly constrained by the available data, leading to strong parameter correlations, practical non-identifiability of part of the neural representation, and compensatory parameter combinations that preserve the observed trajectories while distorting the inferred ($\beta(t)$).

The Fisher Information Matrix (FIM) associated with the HNODE SIR model is full rank and therefore invertible. These results indicate that, unlike in more highly overparameterized settings, the available data are sufficiently informative to independently excite nearly all parameter directions. As a consequence, confidence intervals can be computed for all model parameters, including both mechanistic parameters and those associated with the neural approximation of the transmission rate. 

The confidence intervals reported in Table~\ref{tab:tSIR_HNODEs_CI} provide detailed insight into parameter identifiability in the selected HNODE SIR model. Notably, the epidemiologically meaningful recovery rate $\gamma$ is estimated with a narrow confidence interval, indicating strong local identifiability and close agreement with the ground-truth value used to generate the synthetic data. This result demonstrates that the hybrid formulation, combined with training across multiple experiments, enables robust inference of global mechanistic parameters even under partial observability.

\begin{table}[h!]
\centering
\caption[\texttt{tSIR: }Estimated parameters and confidence intervals for the selected conventional HNODE model.]{
\textbf{\texttt{tSIR: }Estimated parameters and confidence intervals for the selected conventional HNODE model of the SIR system with time-dependent transmission rate.}
The table reports the calibrated parameter values \((\theta^*)\), together with their relative confidence intervals expressed as percentages.
Parameters highlighted in blue have large confidence intervals, indicating weak practical identifiability.
}
\vspace{6pt}
\label{tab:tSIR_HNODEs_CI}
\begin{tabular}{ccc}
\hline
Parameters & $\theta^*$ & CI (\%) \\ \hline
$\gamma$        & 0.101  & 2.65 \\
$W_{111}$ & -2.387 & $2.76 \times 10^{-4}$ \\
$W_{112}$ & -0.133 & 1.06 \\
$W_{113}$ & 0.208  & 0.03 \\
$W_{114}$ & 2.360  & $1.97 \times 10^{-4}$ \\
$b_{11}$  & -0.012 & 20.26 \\
$b_{12}$  & 0.793  & 0.03   \\
$b_{13}$  & -0.426 & 0.09  \\
$b_{14}$  & -0.412 & 0.01  \\
$W_{211}$ & 0.519  & 0.04   \\
$W_{221}$ & 0.606  & 0.01   \\
$W_{231}$ & -0.042 & 20.83 \\
$W_{241}$ & -0.385 & 0.04  \\
$W_{212}$ & 0.149  & 5.85   \\
$W_{222}$ & -0.294 & 0.17  \\
$W_{232}$ & $1.80 \times 10^{-5}$ & 1.45 \\
$W_{242}$ & 0.211  & 0.15   \\
\cellcolor[HTML]{96FFFB}$W_{213}$ & \cellcolor[HTML]{96FFFB}-0.003 & \cellcolor[HTML]{96FFFB}$6.13 \times 10^2$\\
$W_{223}$ & -0.108 & 0.17   \\
$W_{233}$ & 0.242  & 0.06    \\
$W_{243}$ & -0.002 & 0.66   \\
$W_{214}$ & 0.001  & 6.99   \\
$W_{224}$ & -0.169 & 0.04  \\
$W_{234}$ & 0.573  & 0.01   \\
$W_{244}$ & 0.145  & 3.77   \\
$b_{21}$  & -0.052 & 0.42  \\
$b_{22}$  & -0.024 & 0.73  \\
$b_{23}$  & -0.167 & 11.15 \\
$b_{24}$  & 0.094  & 0.09   \\
$W_{311}$ & -0.204 & 0.46  \\
$W_{321}$ & 0.499  & 0.13   \\
$W_{331}$ & 0.148  & 20.78  \\
$W_{341}$ & -0.079 & 4.74  \\
$b_{31}$  & -0.039 & 7.02  \\ \hline
\end{tabular}
\end{table}

Parameters associated with the neural representation of the transmission rate $\beta(t)$ exhibit a markedly heterogeneous uncertainty structure. While several first-layer weights and biases are tightly constrained, indicating highly sensitive directions that strongly shape the inferred dynamics, many other neural parameters display broad or highly asymmetric confidence intervals, in some cases spanning several hundred per cent. A similar pattern is observed in deeper layers, where increased parameter correlations across layers lead to substantial redundancy. In multi-layer architectures, different parameter configurations can generate comparable input–output mappings, allowing compensatory adjustments between layers without significantly altering the reconstructed dynamics. As a consequence, large confidence intervals primarily reflect structural overparameterisation rather than a lack of data informativeness, highlighting the fundamental trade-off between representational flexibility and practical identifiability in deep hybrid dynamical models.

\begin{figure}[h!]
    \centering
    \includegraphics[width=0.80\linewidth]{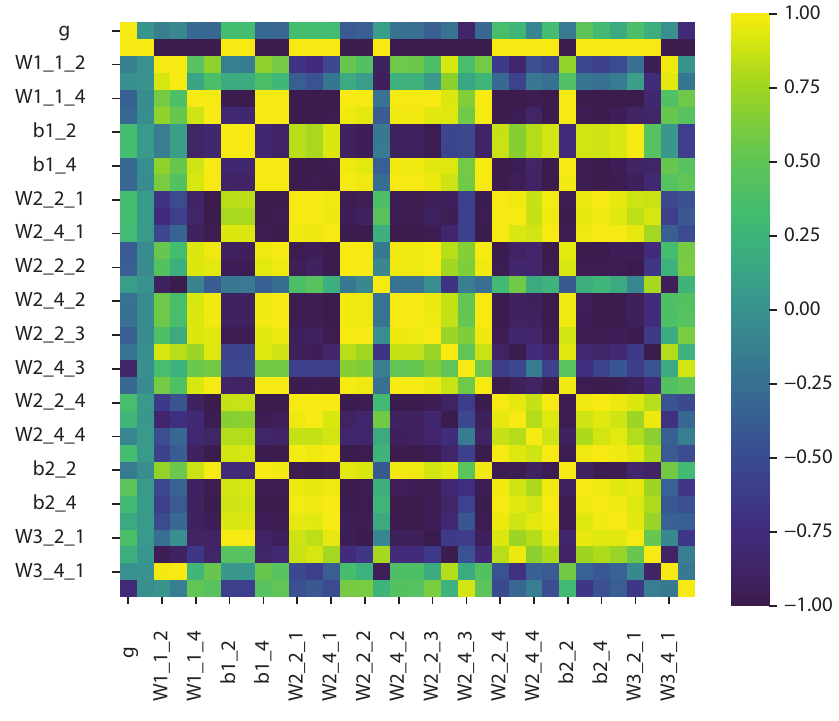}
\caption[\texttt{tSIR: }Parameter correlation matrix for the selected conventional HNODE model.]{
\textbf{\texttt{tSIR: }Parameter correlation matrix for the selected conventional HNODE model of the SIR system with time-dependent transmission rate.}
Rows and columns correspond to the weights and biases of the embedded neural network.
The colour scale shows pairwise Pearson correlation coefficients; values close to \(\pm 1\) indicate strong linear dependencies, whereas values near zero indicate weak correlations.
For readability, only a subset of parameter labels is displayed on each axis.
The presence of strong parameter correlations indicates redundant or weakly informed parameter directions, consistent with limited practical identifiability.
}  \label{fig:tSIR_HNODEs_corr_matrix}
\end{figure}

The parameter correlation matrix of the selected HNODE SIR model (Fig.~\ref{fig:tSIR_HNODEs_corr_matrix}) reveals pronounced positive and negative correlations between parameters belonging to different layers of the neural approximation of the transmission rate $\beta(t)$. This indicates that multiple combinations of weights and biases across layers can generate functionally equivalent dynamics, a characteristic feature of over-parameterised neural architectures embedded within differential equation models.

The strongest correlations are concentrated among neural parameters, particularly between weights in the first and second hidden layers and between weights and biases within the same layer. This correlation structure is consistent with the wide confidence intervals obtained for several neural parameters and reflects the presence of nearly flat directions in the parameter space. In practice, variations in a given layer can be compensated by corresponding adjustments in subsequent layers without significantly affecting the inferred transmission rate $\beta(t)$ or the resulting epidemic trajectories.

In contrast to shallow architectures, deeper networks increase representational flexibility at the cost of amplified parameter redundancy, leading to strong correlations that arise from structural properties rather than limited data informativeness. This trade‑off between expressiveness and identifiability is a known limitation of deep neural models and becomes particularly pronounced when coupled to dynamical systems.

\pagebreak

Across the four benchmark systems, the Supplementary analyses show that predictive accuracy alone does not guarantee practical identifiability. Conventional NODE and HNODE workflows can reproduce observed trajectories while retaining redundant, weakly sensitive or highly correlated parameter directions. In contrast, the iNODE workflow links architecture selection to parameter uncertainty, allowing compact neural differential equations to be selected according to both predictive performance and the information content of the data. The results also show that experimental diversity can improve practical identifiability more effectively than increasing observations along a single trajectory.

\bibliographystyle{unsrt}
\bibliography{biblio}